\documentclass[a4paper, 10pt, twoside, notitlepage]{amsart}
\usepackage{amsmath,amscd}
\usepackage{amssymb}
\usepackage{comment}
\usepackage{graphicx, xcolor}

\usepackage{mathrsfs}
\usepackage[ocgcolorlinks, linkcolor=blue]{hyperref}

\usepackage{bm}
\usepackage{bbm}
\usepackage{url}

\usepackage{geometry}
\usepackage[utf8]{inputenc}
\usepackage{mathtools,amssymb,amsmath}
\usepackage{esint}
\usepackage{tikz}
\usepackage{dsfont}
\usepackage{relsize}
\usepackage{url}
\usepackage{xcolor}
\usepackage{graphicx}
\usepackage{mathrsfs}
\usepackage[shortlabels]{enumitem}
\usepackage{lineno}
\usepackage{amsmath}
\usepackage{enumitem}
\usepackage{amsthm} 
\usepackage{verbatim}
\usepackage{dsfont}
\usepackage{booktabs}
\numberwithin{equation}{section}

\usepackage{slashed}

\allowdisplaybreaks

\mathtoolsset{showonlyrefs}

\graphicspath{{images}}

\newtheorem{theorem}{Theorem}[section]
\newtheorem{lemma}[theorem]{Lemma}
\newtheorem{definition}{Definition}[section]
\newtheorem{corollary}[theorem]{Corollary}

\newtheorem{proposition}[theorem]{Proposition}
\newtheorem{example}{Example}[section]

\newtheorem{remark}{Remark}[section]

\title[Inverse semilinear fourth-order elliptic equations]{Single-point determination of a potential and a nonlinearity in a semilinear biharmonic equation}

\author[S. Fu]{Song-Ren Fu}
\address{School of Mathematics and Statistics, Northwestern Polytechnical University, Xi’an 710129, Shaanxi, PR China.}
\email{songrenfu@nwpu.edu.cn}

\author[L. Li]{Li Li}
\address{Yau Mathematical Sciences Center, Tsinghua University, Beijing, PR China.}
\email{lili19940301@mail.tsinghua.edu.cn}

\author[H. Liu]{Hongyu Liu}
\address{Department of Mathematics, City University of Hong Kong, Hong Kong SAR, China}
\email{hongyliu@cityu.edu.hk}

\author[Z. Zhang]{Zhikun Zhang}
\address{School of Mathematics and Statistics, Northwestern Polytechnical University, Xi’an 710129, Shaanxi, PR China.}
\email{zhikunzhang@nwpu.edu.cn}

\newcommand{\R}{{\mathbb R}}

\newcommand{\Om}{\Omega}

\DeclareMathOperator{\supp}{supp} 

\begin{document}

\begin{abstract}
In this paper, we study the inverse problem of determining the coefficient $q$ and the nonlinearity $F$ in the semilinear biharmonic equation $\Delta^2 u+q(x)u=F(x,u)$ with Navier boundary data. The measurement consists of $\partial_\nu u(x_0)$ at a fixed boundary point $x_0$. We apply the mixed variation approach to linearize the nonlinear problem and construct semiclassical boundary-localized solutions for the linearized equation. Under suitable analyticity, nonresonance, and nondegeneracy assumptions on $q$ and $F$, we prove that this boundary point measurement uniquely determines both $q$ and $F$. We also establish boundary determination of all normal derivatives of $q$ and of the first nonzero Taylor coefficient of $F$ on a nonempty relatively open subset of the boundary. When the linear potential $q$ is known, we prove that either component of the Navier boundary input alone suffices to determine $F$, and that the polynomial nonlinearity can be uniquely determined without requiring any analyticity conditions.  For suitable finite-dimensional classes of nonanalytic potentials, we also obtain, for fixed $F$, the stable determination of $q$ from finitely many boundary inputs and the associated point measurements. In addition, as an important supplement, we provide nontrivial examples verifying the main assumptions, as well as counterexamples concerning the recovery of the linear coefficient. Finally, we present numerical experiments that demonstrate the theoretical results.
\medskip


\end{abstract}
\maketitle

\setlength{\parskip}{0.8ex} 

\section{Introduction}\label{Sec-intro}

The biharmonic operator $\Delta^2$ is a fundamental fourth-order elliptic operator arising, for instance, in the modeling of thin elastic plates (e.g., see \cite{lagnese1988modelling}). A boundary value problem of fourth-order elliptic equation requires two boundary conditions, whose choice reflects the mechanical support of the plate. The homogeneous Navier boundary conditions $u=\Delta u=0$ are  employed as an idealized hinged-plate model. For the unperturbed equation, they also allow a reduction to two successive Dirichlet problems for the Laplacian (see for example \cite{GazzolaGrunauSweers2010}).

Forward problems for fourth-order elliptic equations have drawn significant interest and yielded fruitful results. Classical elliptic boundary regularity is discussed in \cite{Grisvard1985}, while Green functions and positivity for polyharmonic problems are investigated in \cite{GazzolaGrunauSweers2010}. Semilinear problems raise further questions regarding existence and nonexistence, depending on the growth rate of the nonlinear term. The work \cite{GazzolaGrunauSquassina2003} studied the existence and nonexistence of solutions for biharmonic equations with critical growth, using variational methods to overcome the loss of compactness, with particular attention to the Navier boundary case.

Beyond these investigations for forward problems, the inverse problems associated with these equations are a subject of growing interest and rapid development. Inverse boundary problems for polyharmonic operators include the recovery of lower-order perturbations from full Cauchy data and, specifically for the biharmonic operator with Navier boundary conditions, from partial boundary measurements, have been investigated in \cite{KrupchykLassasUhlmann2014,KrupchykLassasUhlmann2012}. Nonlinear tensorial perturbations of the biharmonic operator have also been recovered from full boundary data \cite{BhattacharyyaKrupchykSahooUhlmann2025}.  
Motivated by these works, we consider the inverse problem for semilinear fourth-order elliptic equations. Specifically, we focus on recovering the unknown parameters from suitable boundary measurements at a single fixed point.

Let $\Om\subset\R^n$ ($n\ge2$) be a bounded domain with smooth boundary $\partial\Omega$. In this paper, we consider the following semilinear fourth-order elliptic equation:
\begin{equation}\label{eq:model}
\begin{cases}
\Delta^2u+q(x)u=F(x,u)  & \text{in } \Om,\\
u=f,\, \Delta u=g & \text{on } \partial\Om,
\end{cases}
\end{equation}
where $(u,\Delta u)|_{\partial\Omega}=(f,g)$ specifies the inhomogeneous Navier boundary conditions, with $f$ and $g$ representing the boundary displacement and bending moment, respectively. The function $q(x)$ acts as a spatial potential or a variable stiffness coefficient, and the nonlinear term $F(x,u)$ accounts for nonlinear responses such as large deflections or external reactions.

\subsection{Formulation of the inverse problems} We next formulate the considered inverse problems. Let $\mathbb F\in\{\mathbb R,\mathbb C\}$. We begin with introducing the following function space
\begin{equation}
\mathcal X_{\mathbb F}\vcentcolon=W^{4-1/p,p}(\partial\Omega;\mathbb F)\times W^{2-1/p,p}(\partial\Omega;\mathbb F), \quad  \max\{1,n/3\}<p<\infty.
\end{equation}
For some constant $\varepsilon_0>0$, we define a set for the boundary pair $(h_1,h_2)$ by
\begin{equation}
B_{\mathcal X_{\mathbb F}}(0,\varepsilon_0)\vcentcolon=\{(h_1,h_2)\in \mathcal X_{\mathbb F}: \|(h_1,h_2)\|_{\mathcal X}:=\|h_1\|_{W^{4-1/p,p}(\partial\Omega)}+\|h_2\|_{W^{2-1/p,p}(\partial\Omega)}\le\varepsilon_0\}.
\end{equation}
Fix $x_0\in\partial\Om$. For the moment, we formally define the single-point measurement map by
\begin{equation}
\Lambda_{q,F}: B_{\mathcal X_{\mathbb R}}(0,\varepsilon_0)\rightarrow\mathbb R,\quad \Lambda_{q,F}(f,g)=\partial_\nu u_{f,g}(x_0),
\end{equation}
where $u_{f,g}$ denotes the solution to the nonlinear problem \eqref{eq:model}, and $\nu$ is the unit outward normal vector field along $\partial\Omega$. For sufficiently small $\varepsilon_0>0$, we will establish the (local) well-posedness result for \eqref{eq:model} under some suitable conditions on the linear coefficient $q$ and the nonlinearity $F$ (see Proposition~\ref{pro:well-posedness-nonlinear}). Hence, the map $\Lambda_{q,F}(f,g)$ is well-defined. 

We emphasize that, in what follows, the field $\mathbb F$ will be either $\mathbb R$ or $\mathbb C$ (we will sometimes drop the symbol $\mathbb F$ when no confusion can arise). If necessary, we will also use the complexification of certain real functions and maps (see Remark~\ref{rem:complexification} for more details).

We will also introduce the following one component boundary input spaces
\begin{equation}
\mathcal X^{(0)}:= W^{4-1/p,p}(\partial\Omega)\times\{0\},\quad \mathcal X^{(1)}:=\{0\}\times W^{2-1/p,p}(\partial\Omega).
\end{equation}
The purpose of introducing the above two function spaces is to recover the unknown parameters using only $u|_{\partial\Omega}$ or $\Delta u|_{\partial\Omega}$ as inputs, rather than the full Navier boundary data. Throughout this paper, unless otherwise stated, we assume that the exponent $p$ satisfies
\begin{equation}\label{condition-for-p}
\max\{1,n/3\}<p<\infty.
\end{equation}

The major inverse problem study in this paper is stated as follows.

{\bf (IP)} \emph{Can we recover the coefficient $q$ and the nonlinearity $F$ from the knowledge of the single-point measurement map $\Lambda_{q,F}$}?

For general $q$ and $F$, we point out that the inverse problem possesses a natural coefficient gauge invariance. Indeed, for any sufficiently regular real-valued function $b=b(x)$, set
\begin{equation}
\widetilde q=q+b,\quad  \widetilde F(x,u)=F(x,u)+b(x)u.
\end{equation}
Then, we have 
\begin{equation}
\Delta^2u+\widetilde q(x)u-\widetilde F(x,u)=\Delta^2u+q(x)u-F(x,u),
\end{equation}
and hence the pairs $(q,F)$ and $(\widetilde q,\widetilde F)$ define
exactly the same nonlinear boundary value problem. In particular, we have $\Lambda_{q,F}=\Lambda_{\widetilde q,\widetilde F}$.
This shows that the boundary measurements cannot distinguish $q$ from the linear part $\partial_uF(\cdot,0)$. It is only possible to recover $ q-\partial_uF(\cdot,0)$. Hence, regarding the restricted single-point measurement proposed in this work, we are required to impose some suitable conditions on $q$ and $F$ to break the gauge invariance.

\subsection{Main assumptions and results}\label{sec:setting}
Denote by
\begin{equation}
L_q:=\Delta^2+q(x).
\end{equation}

We introduce the first assumption:

{\bf (A)} \emph{ $q\in L^\infty(\Om;\R)$ and zero is not an eigenvalue of the operator $L_q$. That is, if $ v\in H^4(\Om)$ satisfies the equation $L_qv=0$ with $v=\Delta v=0$ on $\partial\Omega$, then $v\equiv 0$.}

We fix a nonempty connected relatively open subset $\Gamma$ of $\partial\Omega$ such that
\begin{equation}\label{eq:fixed-Gamma}
\Gamma\Subset\partial\Omega\setminus\{x_0\}.
\end{equation}
For the global unique determination of $q,F$, we fix an open connected set $\mathcal O\subset\mathbb R^n$ such that $\overline\Omega\subset\mathcal O$.
For the pair $(q,F)$, we give the following definition:
\begin{definition}\label{def:admissible}
Let $q\in L^\infty(\Omega)$ be a real-valued function and $F:\mathcal{O}\times I\to\mathbb{R}$, where $I$ is an open interval with $0\in I$. We say $(q,F)$ is \emph{an admissible pair} if they satisfy the following conditions:
\begin{enumerate}[label=(\roman*)]
\item $q$ is real analytic in $\Om$ and admits an analytic extension to an open neighborhood  of $\overline\Gamma$;
\item $F\in C^\omega(\mathcal O\times I;\R)$, where $C^\omega$ denotes the analytic function space;
\item there exists an integer $m(F)\ge2$ such that, with $a_F(x):=\partial_u^{m(F)}F(x,0),$
one has, for any $x\in\mathcal O$, $\partial_u^kF(x,0)=0$ for $0\le k<m(F)$, and $a_F|_\Gamma\not\equiv0$;
\item Assumption {\bf (A)} holds for $q$.
\end{enumerate}
\end{definition}

\subsubsection{Infinite-dimensional potentials and uniqueness results}
As for the main theorems of this paper, we first consider the unique determination of $q$ and $F$ (or their boundary normal derivatives), assuming $q$ belongs to certain infinite-dimensional function space.

The following theorem establishes the global determination of $q$ and $F$ from $\Lambda_{q,F}$.
\begin{theorem}\label{thm:main}
For $j=1,2$, let $(q_j,F_j)$ be an admissible pair, with $F_1$ and $F_2$ defined on the same product domain $\mathcal O\times I$. Suppose that there exists a sufficiently small $\varepsilon>0$ such that
\begin{equation}\label{eq:map-equality}
\Lambda_{q_1,F_1}(H)=\Lambda_{q_2,F_2}(H)
\end{equation}
for all $H\in\mathcal X_{\mathbb R}$ with
$\|H\|_{\mathcal X}<\varepsilon$. Then
\begin{equation}
q_1=q_2\quad\text{\rm in}\,\, \Om, \quad F_1=F_2\quad {\rm on}\,\, \Om\times I.
\end{equation}
\end{theorem}

\begin{remark}
We note that, ``single-point measurement'' refers to the location of the scalar observation, not to the number of boundary excitations. The data in Theorem~\ref{thm:main} consist of the map $H=(f,g)\mapsto\Lambda_{q,F}(H) =\partial_\nu u_H(x_0)$
for all $H$ in a neighborhood of zero in the infinite-dimensional space $\mathcal X_{\mathbb R}$. Recovery of unknowns from
finitely many fixed inputs is asserted only in the finite-dimensional coefficient settings stated later.
\end{remark}

\begin{remark}\label{rem:nonanalytic-polynomial}
Section~\ref{sec:conclusion} gives unique recovery results for $F(x,u)=\sum_{k=2}^{M}\frac{a_k(x)}{k!}u^k$
without assuming that the coefficients $a_k$ are spatially analytic. When $q=0$, all the coefficients can be recovered in dimensions $n\geq2$ from sufficiently small inputs $(f,0)$, with $f$ supported in any nonempty relatively open subset of
$\partial\Omega$. For a fixed, known nonresonant $q\in L^\infty(\Omega)$, recovery of coefficients from full Navier boundary input $(f,g)$ holds for $3\leq n\leq6$.  These results do not remove the analyticity assumptions from Theorem~\ref{thm:main} for an unrestricted unknown potential.
\end{remark}

As a direct by-product of the above Theorem~\ref{thm:main}, we have the following corollary:
\begin{corollary}\label{cor:scalar-case}
For $j=1,2$, let $q_j\in L^\infty(\Om;\R)$ satisfy Assumption {\bf (A)}, be real analytic in $\Om$, and admit a real-analytic extension to an open neighborhood of $\overline\Gamma$. Let $F_j:I\to\R$ be real analytic on the same connected interval $I\ni0$, with $F_j(0)=F_j'(0)=0$ and $F_j\not\equiv0$ near $0$. If \eqref{eq:map-equality} holds for all sufficiently small real boundary inputs, then
\begin{equation}
q_1=q_2\quad {\rm in}\,\, \Om, \quad F_1=F_2\quad {\rm on}\,\, I.
\end{equation}
\end{corollary}

The following corollary shows that, for a known potential, each of the boundary inputs $(f,0)$ and $(0,g)$ individually determines the nonlinearity $F$.

\begin{corollary}\label{cor:known-q-one-component}
Let $q\in L^\infty(\Om;\R)$ satisfy the conditions (i) and (iv) in Definition~\ref{def:admissible}. For $j=1,2$, let $F_j\in C^\omega(\mathcal O\times I;\R)$ satisfy
\begin{equation}
F_j(x,0)=\partial_uF_j(x,0)=0.
\end{equation}
Fix $\sigma\in\{0,1\}$. If
\begin{equation}\label{eq:one-component-map-equality}
\Lambda_{q,F_1}(H)=\Lambda_{q,F_2}(H)
\end{equation}
for all sufficiently small real-valued $H\in\mathcal X^{(\sigma)}$, then $F_1=F_2$ on $\Omega\times I$.
\end{corollary}

Let $h$ be a function function that is smooth up to the boundary. We write
\begin{equation}
\partial_\nu^r h(y):=\frac{d^r}{ds^r}h\bigl(y+s\nu(y)\bigr)\big|_{s=0^-},\quad y\in\partial\Om,\quad r\in\mathbb N_0.
\end{equation}
By Proposition~\ref{pro:well-posedness-nonlinear}, the point-measurement map
\begin{equation}
\Lambda_{q,F}:\mathcal X_{\mathbb R}\supset\mathcal U\rightarrow\mathbb R,\quad \Lambda_{q,F}(H)=\partial_\nu u_H(x_0),
\end{equation}
is $C^\infty$ on a neighborhood $\mathcal U$ of zero. For $k\geq1$
and $H_1,\ldots,H_k\in\mathcal X_{\mathbb R}$, its $k$th Fr\'echet
derivative at zero is given by
\begin{equation}\label{eq:measurement-frechet-derivative}
D^k\Lambda_{q,F}(0)[H_1,\ldots,H_k]=\partial_{\varepsilon_1,\cdots,\varepsilon_k}^k
\Lambda_{q,F}\!\big(\varepsilon_1 H_1+\cdots+\varepsilon_kH_k\big)\big|_{\varepsilon_1=\cdots=\varepsilon_k=0}.
\end{equation}
Thus $D^k\Lambda_{q,F}(0)$ is a continuous symmetric $k$-linear
form on $\mathcal X_{\mathbb R}$. When complex-valued solutions of
the linearized equation are used, the same notation denotes its
unique complex multilinear extension to $\mathcal X_{\mathbb C}$.

We have the following boundary determination result for $q,F$.
\begin{theorem}\label{thm:boundary-jets}
For $j=1,2$, let $q_j\in L^\infty(\Om;\R)$ satisfy Assumption~{\bf (A)} and be $C^\infty$ up to
$\partial\Om$ in an open neighborhood of $\overline\Gamma$.  Let $F_j(x,\tau)$ be smooth for $x\in\overline\Om$ and $\tau$ near zero, and set $ a_j(x):=\partial_u^{m_*}F_j(x,0)$ for some $m_*\ge 2$. Suppose that
\begin{equation}
 F_j(x,0)=\partial_uF_j(x,0)=0,\quad \partial_u^kF_j(x,0)=0\ (2\le k<m_*),\quad a_j|_\Gamma\not\equiv0,
\end{equation}
for $j=1,2$.
Assume that 
\begin{equation}\label{cond:DN-m-star}
D\Lambda_{q_1,F_1}(0)=D\Lambda_{q_2,F_2}(0), \quad D^{m_*}\Lambda_{q_1,F_1}(0)=D^{m_*}\Lambda_{q_2,F_2}(0).
\end{equation}
Then there exists a nonempty open set $U\Subset\Gamma$ such that
\begin{equation}\label{eq:boundary-jet-conclusion}
\partial_\nu^rq_1=\partial_\nu^rq_2,\quad \partial_\nu^ra_1=\partial_\nu^ra_2 \quad {\rm on}\,\, U,\quad r=0,1,2,\ldots.
\end{equation}
If both $a_1$ and $a_2$ are nowhere zero on $\Gamma$, the identities hold on all of $\Gamma$.
\end{theorem}

\begin{remark}
The exclusion of $x_0$ in \eqref{eq:fixed-Gamma} is needed for the local smoothness of the point adjoint function $\Psi_q$ (see Lemma~\ref{lem:point-adjoint}).  If $q_1,q_2$ are $C^\infty$ up to all of $\partial\Om$ and $a_1,a_2$ are nowhere zero on $\partial\Omega$, the proof of Theorem~\ref{thm:boundary-jets} applies on every sufficiently small open subset of $\partial\Om\setminus\{x_0\}$.  Hence \eqref{eq:boundary-jet-conclusion} holds on $\partial\Om\setminus\{x_0\}$ and, by continuity of the normal derivatives, also at $x_0$. 
\end{remark}

\begin{remark}\label{rem:boundary-vanishing-profile}
The assumption $a_j|_\Gamma\not\equiv0$ in Theorem~\ref{thm:boundary-jets} is used in the present
proof of simultaneous recovery. It is not claimed to be necessary for the unique determination of the nonlinearity $F$. For example, on
$\Omega=B(0,1)$ the nonlinearity $F(x,u)=\frac{1-|x|^2}{m!}u^m$ has a nonzero leading coefficient in $\Omega$ whose
boundary trace vanishes. The identity \eqref{eq:unknown-profile-leading-coefficient}, obtained from Proposition~\ref{prop:coefficient-extraction} is then unavailable. However, if $q=1$ is known,
this nonlinearity remains covered by Corollary~\ref{cor:known-q-one-component}.
\end{remark}

\subsubsection{Finite-dimensional potentials and stability results}\label{sec:finite-dimensional-statements}
We next consider the recovery of potentials belonging to a certain finite-dimensional space.
Before stating the theorems, we first establish a condition under which a potential within a finite-dimensional class is uniquely determined by its boundary normal derivatives.

\begin{definition}\label{def:jet-visible}
Let $E\subset L^\infty(\Omega;\mathbb R)$ be finite dimensional, and let $U\Subset\partial\Omega\setminus\{x_0\}$ be a relatively open set. Assume that any $e\in E$ is smooth up to $\partial\Omega$ near $\overline U$.  We say that $E$ is \emph{determined by its normal derivatives on $U$} if
\begin{equation}\label{eq:jet-visibility}
\bigl\{e\in E:\partial_\nu^r e|_U=0\,\, {\rm for\,\, every}\,\,\, r\geq0\bigr\}=\{0\}.
\end{equation}
\end{definition}

The finite-dimensional assumption permits global recovery of the potential without requiring it to be real analytic in the interior of $\Omega$.

\begin{theorem}\label{thm:finite-dimensional-uniqueness}
Let $\Gamma_a\Subset\Gamma$ be a nonempty relatively open set.
Let $q_{\mathrm{ref}}\in L^\infty(\Omega;\mathbb R)$, and let
$E\subset L^\infty(\Omega;\mathbb R)$ be finite dimensional.
Assume that $q_{\mathrm{ref}}$ and any $e\in E$ are smooth up to
$\partial\Omega$ near $\overline\Gamma_a$, and that $E$ is
determined by its normal derivatives on $\Gamma_a$ in the sense of
Definition~\ref{def:jet-visible}.
For $j=1,2$, let $q_j\in q_{\mathrm{ref}}+E$ satisfy Assumption {\bf (A)}, and let $F_j$ satisfy conditions (ii)--(iii) in Definition~\ref{def:admissible} for some $m_j=m(F_j)\ge 2$ with $a_{F_j}:=\partial_u^{m_j}F_j(\cdot,0)$ nowhere zero on $\Gamma_a$. If
\begin{equation}
\Lambda_{q_1,F_1}(H)=\Lambda_{q_2,F_2}(H)
\end{equation}
for all sufficiently small real boundary input $H\in\mathcal X_{\mathbb R}$, then
\begin{equation}
q_1=q_2\quad {\rm in}\,\, \Omega,\quad F_1=F_2\quad {\rm on}\,\, \Omega\times I.
\end{equation}
\end{theorem}

For given $\rho>0$, and $\sigma\in\{0,1\}$, we  define
\begin{equation}
\overline B_\rho:=\{H\in\mathcal X_{\mathbb R}:\|H\|_{\mathcal X}\leq\rho\},\quad
\mathscr Y_\rho:=C_b(\overline B_\rho;\mathbb R),
\end{equation}
\begin{equation}\label{eq:restricted-data-spaces}
\overline B_\rho^{(\sigma)}:=\overline B_\rho\cap\mathcal X^{(\sigma)},\quad
\mathscr Y_\rho^{(\sigma)}:=C_b(\overline B_\rho^{(\sigma)};\mathbb R).
\end{equation}
Whenever the small solution is well-defined on $\overline B_\rho$, we write
\begin{equation}\label{eq:finite-data-map}
\mathscr D_F(q)
:=\Lambda_{q,F}|_{\overline B_\rho}
\in\mathscr Y_\rho.
\end{equation}

Assuming the full rank condition on the tensor $D_q\mathscr D_F$ at a fixed function $q_*$, the next result addresses the local stable determination of $q$ from finitely many measurements.

\begin{theorem}\label{thm:local-sensitivity}
Let $F\in C^\omega(\mathcal O\times I;\mathbb R)$ satisfy
$F(x,0)=\partial_uF(x,0)=0$. Let $q_*\in L^\infty(\Omega;\mathbb R)$ satisfy Assumption {\bf (A)}, and let
$E\subset L^\infty(\Omega;\mathbb R)$ be finite dimensional. Suppose that there exist $\rho>0$ and a neighborhood $\mathcal V$ of $q_*$ in
$q_*+E$ such that $\mathscr D_F:\mathcal V\to\mathscr Y_\rho$ is $C^1$. If
\begin{equation}\label{eq:local-sensitivity-rank}
\ker\bigl(D_q\mathscr D_F(q_*)|_E\bigr)=\{0\},
\end{equation}
then there are boundary inputs $H_1,\ldots,H_L\in\overline B_\rho$, a neighborhood $\mathcal U_*\subset\mathcal V$ of $q_*$, and  a constant $C>0$, such that
\begin{equation}\label{eq:local-finite-input-stability}
\|q_1-q_2\|^2_{L^\infty(\Omega)}\leq C\sum_{\ell=1}^L\left|\Lambda_{q_1,F}(H_\ell)-\Lambda_{q_2,F}(H_\ell)\right|^2,\quad q_1,q_2\in\mathcal U_*.
\end{equation}
\end{theorem}

The same argument as in Theorem~\ref{thm:local-sensitivity} yields the following corollary when the Navier boundary inputs are either $(f,0)$ or $(0,g)$.

\begin{corollary}\label{cor:one-component-local-stability}
Under the assumptions of Theorem~\ref{thm:local-sensitivity}, fix
$\sigma\in\{0,1\}$ and set
\begin{equation}
\mathscr D_F^{(\sigma)}(q):=\Lambda_{q,F}|_{\overline B_\rho^{(\sigma)}}\in\mathscr Y_\rho^{(\sigma)}.
\end{equation}
If
\begin{equation}\label{eq:restricted-sensitivity-rank}
\ker\bigl(D_q\mathscr D_F^{(\sigma)}(q_*)|_E\bigr)=\{0\},
\end{equation}
then there are
$H_1,\ldots,H_L\in\overline B_\rho^{(\sigma)}$, a neighborhood $\mathcal U_*$ of $q_*$ in $q_*+E$, and $C>0$ such that
\begin{equation}\label{eq:restricted-local-stability}
\|q_1-q_2\|^2_{L^\infty(\Omega)}\leq C\sum_{\ell=1}^L\left|\Lambda_{q_1,F}(H_\ell)-\Lambda_{q_2,F}(H_\ell)\right|^2,\quad q_1,q_2\in\mathcal U_*
\end{equation}
\end{corollary}

\begin{remark}\label{rem:restricted-rank-necessary}
The rank condition \eqref{eq:restricted-sensitivity-rank} is necessary for a local stability estimate on a neighborhood of $q_*$ in $q_*+E$. Indeed, by applying such an estimate to $q_*$ and $q_*+te$, dividing by $|t|$, and letting $t\to0$, we deduce that $D_q\mathscr D_F^{(\sigma)}(q_*)e=0$ implies $e=0$.
\end{remark}

By analyzing the set $E$ defined in Definition~\ref{def:jet-visible}, we obtain the following corollary, which provides a sufficient condition for the rank assumption \eqref{eq:local-sensitivity-rank}.
\begin{corollary}\label{cor:automatic-sensitivity-rank}
Let $q_*\in L^\infty(\Omega;\mathbb R)$ satisfy Assumption {\bf (A)}, and let $E\subset L^\infty(\Omega;\mathbb R)$ be finite dimensional. Let $\Gamma_a\Subset\Gamma$ be a nonempty relatively open set. Suppose that $q_*$ and every $e\in E$ are smooth up to $\partial\Omega$ in an open neighborhood of $\overline\Gamma_a$, and that $E$ is determined by its normal
derivatives on $\Gamma_a$ in the sense of Definition~\ref{def:jet-visible}. Let $F$ satisfy the conditions of Theorem~\ref{thm:local-sensitivity} and (iii) in Definition~\ref{def:admissible} for some $m_*=m(F)\geq2$, and assume in addition that $a_F$ is nowhere zero on $\Gamma_a$. Then
\begin{equation}\label{eq:automatic-sensitivity-rank}
\ker\bigl(D_q\mathscr D_F(q_*)|_E\bigr)=\{0\}.
\end{equation}
Hence, the local stability at $q_*$ as in Theorem~\ref{thm:local-sensitivity} holds.
\end{corollary}

Global uniqueness yields a stability modulus on compact subsets of the same finite-dimensional class. A uniform rank condition upgrades this conclusion to a stability estimate involving finitely many fixed inputs.

\begin{theorem}\label{thm:finite-dimensional-stability}
Let $q_{\mathrm{ref}}$, $E$, $\Gamma_a$ and the fixed nonlinearity $F$ satisfy the conditions of Theorem~\ref{thm:finite-dimensional-uniqueness}. Suppose that the first nonzero coefficient $a_F=\partial_u^{m_*}F(\cdot,0)$ of $F$ is nowhere zero on $\Gamma_a$. Let $\mathcal K\subset q_{\mathrm{ref}}+E$ be compact, and suppose that every $q\in\mathcal K$ satisfies Assumption {\bf (A)}. Assume further that there exists $\rho>0$, independent of $q\in\mathcal K$, such that
\begin{equation}
\mathscr D_F(q)=\Lambda_{q,F}|_{\overline B_\rho}\in\mathscr Y_\rho
\end{equation}
is well defined for $q\in\mathcal K$. The following statements hold.

\begin{enumerate}[label=(\roman*)]
\item There exists a nondecreasing function $\omega:[0,\infty)\to[0,\infty)$, with $\omega(s)\to0$ as $s\to0$, such that
\begin{equation}\label{eq:finite-modulus-stability}
\|q_1-q_2\|_{L^\infty(\Omega)}\leq\omega\!\left(\|\mathscr D_F(q_1)-\mathscr D_F(q_2)\|_{\mathscr Y_\rho}\right),
\quad q_1,q_2\in\mathcal K.
\end{equation}

\item There exist finitely many inputs $H_1,\ldots,H_L\in\overline B_\rho$ and constants $C,C'>0$ such that
\begin{equation}\label{eq:finite-lipschitz-stability}
\|q_1-q_2\|^2_{L^\infty(\Omega)}\leq C\sum_{\ell=1}^{L}\bigl|\Lambda_{q_1,F}(H_\ell)-\Lambda_{q_2,F}(H_\ell)\bigr|^2
\leq C'\|\mathscr D_F(q_1)-\mathscr D_F(q_2)\|^2_{\mathscr Y_\rho}, \quad q_1,q_2\in\mathcal K.
\end{equation}
\end{enumerate}
\end{theorem}

\begin{remark}\label{rem:finite-stability-scope}
The modulus in \eqref{eq:finite-modulus-stability} need not have an explicit rate.  Both statements in Theorem~\ref{thm:finite-dimensional-stability} concern the stable recovery of $q$ with $F$ fixed.
\end{remark}

\subsection{Related literature and contribution}
The study of inverse problems for fourth-order elliptic (or polyharmonic) equations has attracted considerable attention, including the recovery of unknown sources, coefficients, and nonlinearities from various boundary or scattering measurements. 
The recovery of lower-order coefficients from boundary measurements is more directly related to the present problem. Krupchyk, Lassas and Uhlmann proved that the full Cauchy data determine the first- and zeroth-order perturbations of a polyharmonic operator \cite{KrupchykLassasUhlmann2014}. For the biharmonic operator with Navier boundary conditions, they obtained in \cite{KrupchykLassasUhlmann2012}, a corresponding uniqueness result from partial boundary measurements. Yan studied inverse boundary problems for first order perturbations of the biharmonic operator on a conformally transversally anisotropic Riemannian manifold \cite{Yan2021}. The analysis mainly relies on the construction of complex geometrical optics solutions (CGO) and Carleman estimates, with boundary measurements available on open observation sets.

For inverse source problems, stability estimates for the biharmonic operator with a potential and for the damped plate equation were established in the works \cite{LiYaoZhao2021SIAM,LiYaoZhao2021IP}. Further results concern increasing stability
\cite{Si2026} and the characterization of radiating and nonradiating sources \cite{ChangGuoYinZhao2024}.
The work \cite{ChangGuoZhao2024SIAM} developed a reconstruction method from multifrequency phaseless data using auxiliary sources for phase retrieval. Recovery of random source has also been investigated for linear and semilinear stochastic
biharmonic wave equations \cite{LiWang2022,WangXuZhao2025}. For the recovery of an unknown cavity, \cite{DongLi2024} established uniqueness from far-field patterns and from suitable phaseless near-field data. Sampling and factorization methods have also been developed for clamped cavities and absorbing penetrable scatterers \cite{HarrisLiOzochiawaeze2026,CejaHarrisOzochiawaeze2026}. Related inverse scattering results address first-order
perturbations of the biharmonic Schr\"odinger operator \cite{XuZhao2024} and quasi-linear biharmonic operators
\cite{TyniSerov2019,HarjuKultimaSerovTyni2021}. These studies concern different unknowns and measurement configurations, but provide the broader context for inverse problems governed by fourth-order equations.

Inverse problems for nonlinear elliptic equations provide a mechanism for extracting information of coefficients from interactions between solutions.  In \cite{isakov1993uniqueness}, Isakov demonstrated that a nonlinear term $a(x,z)$ in a nonlinear Schr\"odinger equation $\Delta u+a(x,u)=0$ can be recovered from the DN map, using the first-order linearization technique under suitable assumptions on $a(x,z)$. Many subsequent works concerning inverse second-order elliptic equations arising. Feizmohammadi and Oksanen \cite{FeizmohammadiOksanen2020} and Lassas, Liimatainen, Lin and Salo \cite{LassasLiimatainenLinSalo2021} developed higher-order linearization arguments for recovering coefficients in semilinear equations.  The work \cite{LassasLiimatainenLinSaloPartial2021} also treated partial boundary data and simultaneous recovery of coefficients and an unknown boundary component. For general nonlinearities, based on constructing a good solution map that locally parametrizes solutions of the nonlinear equation by solutions of the linearized equation, the authors in \cite{JohanssonNurminenSalo2026} showed that a general nonlinearity, and not just its Taylor series, can be uniquely determined up to gauge near a fixed solution. 

Related developments have extended inverse boundary value problems to fourth-order equations. In \cite{BhattacharyyaKrupchykSahooUhlmann2025}, the authors proved uniqueness for third-order nonlinear tensorial perturbations of the biharmonic operator from full boundary data. In \cite{BhattacharyyaKumar2024}, the authors proved that local Dirichlet-to-Neumann data uniquely determine all tensorial perturbations up to order \(m\) and the potential in a polyharmonic operator of order \(2m\), assuming the highest-order perturbation is partially isotropic and the inaccessible boundary lies in a hyperplane. Local determination of a nonlinear term depending on $(x,u,\nabla u,\Delta u)$ from Cauchy data on the whole boundary, for solutions close to a fixed solution, was investigated in \cite{Nurminen2025inverse}, where an associated gauge freedom in the general case was also addressed. The main argument relies on first-order linearization, Runge approximation, and uniqueness for the resulting linear inverse problem.

For observations at a single boundary point, the closest second-order result is due to Salo and Tzou \cite{SaloTzou2023}. They showed that the spatial coefficient in a power nonlinearity is determined by the Neumann data at one fixed boundary point, even when the Dirichlet inputs are supported in an arbitrarily small nonempty open subset of the boundary. Later,  \cite{Nurminen2023Unbounded} extended this single-point result to unbounded coefficients in $L^{n+\varepsilon}(\Omega)$ with $\varepsilon>0$. Kian \cite{Kian2024} established uniqueness and stability for spatially independent quasilinear constitutive laws from affine Dirichlet inputs and flux observations at finitely many boundary points, including single-point results under additional structural assumptions. More recently, the work \cite{KianLiuWangZheng2024} studied the recovery and numerical reconstruction of a spatially independent semilinear term $F(u)$ from a one-parameter family of Dirichlet inputs and observations at one boundary point. 

Point observations have also been studied in inverse problems for nonlinear wave equations. Tzou established recovery of a Riemannian manifold from observations at a fixed spatial point \cite{Tzou2023}, while Nursultanov, Oksanen and Tzou determined a Lorentzian conformal structure from observations along a single timelike curve \cite{NursultanovOksanenTzou2025}. In both works, the observations are recorded over a time interval.

The present work addresses simultaneous recovery of a linear potential $q(x)$ and a spatially dependent nonlinearity $F(x,u)$ from a scalar boundary measurement map at a single fixed point. A central difficulty is that the unknown potential enters the first linearized operator and therefore also affects the solutions appearing in the nonlinear interaction identities.
We prove that the point-measurement map determines all boundary normal derivatives of $q$ and of the leading Taylor coefficient
of $F$ on a nonempty open boundary set where the latter does not vanish. Under the stated real-analytic extension assumptions,
these boundary identities lead to global uniqueness. For a known potential $q$, we further establish unique recovery of $F$ using either the Navier boundary inputs $(f,0)$ or the inputs $(0,g)$ alone.

We also consider potentials in certain finite-dimensional classes. Lipschitz stability under such structural restrictions has been established for the conductivity equation by Alessandrini and Vessella \cite{AlessandriniVessella2005} and Harrach \cite{Harrach2019}. For the single-point measurement map considered in our work, analytic continuation is replaced by the requirement that the boundary normal derivatives uniquely determine the potential within its prescribed class. With $F$ fixed, an explicit sensitivity rank condition yields local Lipschitz stability from finitely many Navier boundary inputs. The boundary determination property provides a sufficient condition for this rank condition. On compact admissible classes, we obtain a modulus of stability and, under a uniform rank condition, a stable determination of the potantial $q$.

The main contributions of this paper are as follows.

(1) \emph{Simultaneous recovery from a point measurement.}
Under the stated analyticity and normalization assumptions, the values of $\partial_\nu u_{f,g}(x_0)$ for all sufficiently small Navier boundary inputs $(f,g)$ uniquely determine both the potential $q$ and the nonlinearity $F$. Moreover, when $q=0$, using boundary data of the form $(f,0)$, where $f$ is supported on an arbitrary nonempty relatively open subset of $\partial\Omega$, we can recover the nonlinear coefficients $a_k$ in $F(x,u)=\sum_{k=2}^ma_k(x)u^k$ ($m\ge 2$) by single-point measurements, without imposing any analyticity assumption on $a_k$.

(2) \emph{Boundary determination using two orders of the linearization of the measurement map.}
The first-order Fr\'echet derivative and the first nontrivial higher-order derivative of the point-measurement map determine all boundary normal derivatives of $q$ and of the leading Taylor coefficient of $F$ on a nonempty relatively open subset of the boundary.

(3) \emph{Recovery of the nonlinearity from one Navier boundary input component.}
For a known potential $q$, the point-measurement map uniquely determines $F$ when restricted to either the inputs $(f,0)$ or the inputs $(0,g)$. The result remains valid when the leading Taylor coefficient of $F$ vanishes identically on $\partial\Omega$.

(4) \emph{Finite-dimensional uniqueness and stability without analyticity of $q$.}
For some finite-dimensional potential classes determined by their boundary normal derivatives, the map $\Lambda_{q,F}$ uniquely determines both $q$ and $F$ without requiring $q$ to be real analytic in $\Omega$.
With $F$ fixed, finitely many given Navier boundary inputs yield stable determination of $q$ on compact admissible classes satisfying the stated uniform sensitivity rank condition.

The remainder of this paper is organized as follows. Section~\ref{sec:linearization} is devoted to studying the well-posedness of the forward problem and the higher-order linearization of the nonlinear equation \eqref{eq:model}.  In Section~\ref{sec:adjoint}, we establish a key identity for the single-point measurement and introduce a useful adjoint function together with its properties. In Section~\ref{sec:poisson}, we construct some suitable semiclassical boundary-localized solutions to the linearized equation $L_qu=0$. Section~\ref{sec:mainproof} and Section~\ref{sec:finite-dimensional} focus on proving the main theorems of this paper. In Section~\ref{sec:examples}, some nontrivial examples verifying the main assumptions of each theorem are provided, and counterexamples are given to show the non-uniqueness of recovering the linear potential. 
In Section~\ref{sec:numerical-experiment}, we present some numerical examples illustrating the theory on the recovery of the unknowns.  Some conclusions and possible extensions are given in the final Section~\ref{sec:conclusion}.

\section{Well-posedness of the forward problem and linearization}\label{sec:linearization}
In this section, we study the well-posedness of the nonlinear forward problem \eqref{eq:model} and apply mixed variations within the higher-order linearization framework. 

\subsection{Well-posedness for the forward problem}
For given $r\in (1,\infty),$ we define a map $\mathcal A_q$ by
\begin{equation}
\mathcal A_q: W^{4,r}(\Omega)\to L^r(\Omega)\times W^{4-1/r,r}(\partial\Omega)\times W^{2-1/r,r}(\partial\Omega),\quad  u\mapsto \bigl(L_qu,u|_{\partial\Om},\Delta u|_{\partial\Om}  \bigr).
\end{equation}
We have the following property for $\mathcal A_q$.
\begin{lemma}\label{lem:full-navier-isomorphism}
Let $q\in L^\infty(\Omega)$, and let assumption {\bf (A)} hold. For any $1<r<\infty$, the map $\mathcal A_q$ is an isomorphism, such that
\begin{equation}\label{eq:full-navier-estimate}
\|u\|_{W^{4,r}(\Omega)}\le C_r\bigl( \|L_qu\|_{L^r(\Omega)}+\|u|_{\partial\Om}\|_{W^{4-1/r,r}(\partial\Omega)}+\|\Delta u|_{\partial\Omega}\|_{W^{2-1/r,r}(\partial\Omega)} \bigr).
\end{equation}
\end{lemma}

\begin{proof}
For the case $q=0$, we can apply the standard $W^{k,r}$-regularity theorem for the second-order Dirichlet-Laplacian
problem (e.g., see \cite[Theorems 2.4.2.5, 2.5.1.1 and Remark 2.5.1.2]{Grisvard1985}). More precisely, for 
\begin{equation}
 (f,h_1,h_2)\in L^r(\Omega) \times W^{4-1/r,r}(\partial\Omega) \times W^{2-1/r,r}(\partial\Omega),
\end{equation}
we consider the following equation
\begin{equation}\label{eq:w-Dirichlet}
 \Delta w=f\quad\text{in }\Omega,\quad w=h_2\quad\text{on }\partial\Omega.
\end{equation}
The above inhomogeneous elliptic problem admits a unique solution $w\in W^{2,r}(\Omega)$ satisfying the estimate
\begin{equation}\label{eq:w-estimate}
\|w\|_{W^{2,r}(\Omega)} \le C\left( \|f\|_{L^r(\Omega)}+\|h_2\|_{W^{2-1/r,r}(\partial\Omega)}\right).
\end{equation}
For this $w$, we consider the problem
\begin{equation}\label{eq:u-Dirichlet}
\Delta u=w\quad\text{in }\Omega,\quad u=h_1\quad\text{on }\partial\Omega.
\end{equation}
We can also get $u\in W^{4,r}(\Omega)$, and
\begin{equation}
\|u\|_{W^{4,r}(\Omega)}\le C\left(\|w\|_{W^{2,r}(\Omega)}+\|h_1\|_{W^{4-1/r,r}(\partial\Omega)}\right).
\end{equation}
It then follows that
\begin{equation}\label{eq:A0-estimate}
\|u\|_{W^{4,r}(\Omega)}\le C\left( \|f\|_{L^r(\Omega)}+\|h_1\|_{W^{4-1/r,r}(\partial\Omega)}+\|h_2\|_{W^{2-1/r,r}(\partial\Omega)}
\right).
\end{equation}
Moreover, we have $\Delta^2u=f$ with $u|_{\partial\Omega}=h_1$ and $\Delta u|_{\partial\Omega}=h_2$.
Hence $\mathcal A_0$ is an isomorphism.

We next consider the case $q\ne 0$. We write
\begin{equation}
\mathcal A_q=\mathcal A_0+\mathcal B_q, \quad \mathcal B_qu\vcentcolon=(qu,0,0).
\end{equation}
Since $ W^{4,r}(\Omega)\hookrightarrow L^r(\Omega)$ is compact, and multiplication by $q\in L^\infty(\Omega)$ is bounded on $L^r(\Omega)$, thus $\mathcal B_q$ is compact. Therefore $\mathcal A_q$ is a Fredholm operator and $\operatorname{ind}\mathcal A_q=\operatorname{ind}\mathcal A_0=0.$
It remains only to prove that $\ker\mathcal A_q=\{0\}$.

Suppose that $u\in W^{4,r}(\Omega)$ satisfies
\begin{equation}\label{eq:kernel-equation}
\Delta^2u=-qu\quad\text{in }\Omega, \quad
u=\Delta u=0\quad\text{on }\partial\Omega.
\end{equation}
If $r\geq2$, then the embedding $ W^{4,r}(\Omega)\hookrightarrow W^{4,2}(\Omega)=H^4(\Omega),$ and assumption {\bf (A)} immediately give $u=0$ in $\Omega$.

Assume now that $1<r<2$. We are required to show that $u\in H^4(\Omega)$. Suppose that $ u\in W^{4,r_j}(\Omega)$ for some $r_j<2$. If $4r_j<n$, Sobolev embedding gives $W^{4,r_j}(\Omega)\hookrightarrow L^{r_{j+1}}(\Omega)$
for every $r_{j+1}>r_j$ satisfying $\frac1{r_{j+1}}>\frac1{r_j}-\frac4n.$
Hence $qu\in L^{r_{j+1}}(\Omega)$. Define
\begin{equation}
 z\vcentcolon=\mathcal A_0^{-1}(-qu,0,0)\in W^{4,r_{j+1}}(\Omega)\subset W^{4,r_j}(\Omega).
\end{equation}
By \eqref{eq:kernel-equation}, we have $\mathcal A_0(u-z)=0.$ The injectivity of $\mathcal A_0$ yields $u=z$, and therefore $u\in W^{4,r_{j+1}}(\Omega).$

As long as $r_j<2$ and $4r_j<n$, we can choose
\begin{equation}
1/r_{j+1}=\max\left\{1/2,1/r_j-2/n\right\}.
\end{equation}
Hence, after finitely many iterations, one can reach an exponent $r_j$ such that either $r_j\geq2$ or $4r_j\geq n$. This implies $u\in L^2(\Omega).$ If $4r_j>n$, then
$W^{4,r_j}(\Omega)\hookrightarrow L^\infty(\Omega)$, whereas in the critical case $4r_j=n$ one has
$W^{4,r_j}(\Omega)\hookrightarrow L^s(\Omega)$ for any $s\in [2,\infty)$. It follows that $qu\in L^2(\Omega)$. Let
\begin{equation}
 z_2=\mathcal A_0^{-1}(-qu,0,0)\in W^{4,2}(\Omega).
\end{equation}
Using $\mathcal A_0(u-z_2)=0$, we have $u=z_2$, and thus $u\in H^4(\Omega).$ Again using assumption {\bf (A)} gives $u=0$, and $\ker\mathcal A_q=\{0\}$ for any $r\in (1,\infty)$. 
Since $\mathcal A_q$ is Fredholm of index zero, it is surjective and therefore bijective. The bounded inverse theorem yields the estimate \eqref{eq:full-navier-estimate}.
\end{proof}

Recall the function space $\mathcal X_{\mathbb R}=W^{4-1/p,p}(\partial\Omega;\mathbb R)\times W^{2-1/p,p}(\partial\Omega;\mathbb R)$ for $\max\left\{1,n/3\right\}<p<\infty.$ Choosing a constant $\rho_0>0$, we define a neighborhood $\mathcal U_{\rho_0}\subset W^{4,p}(\Omega)$ by
\begin{equation}
\mathcal U_{\rho_0}\vcentcolon=\left\{w\in W^{4,p}(\Omega): \|w\|_{W^{4,p}(\Omega)}<\rho_0  \right\}.
\end{equation}
Denote a map
\begin{equation}
\mathcal F: W^{4,p}(\Om)\times\mathcal X_{\mathbb R}\to  \mathcal Y_p\vcentcolon=L^p(\Omega)\times \mathcal X_{\mathbb R}
\end{equation}
by
\begin{equation}
\mathcal F(u,f,g)\vcentcolon=\bigl(L_qu-F(\cdot,u), u|_{\partial\Om}-f, \Delta u|_{\partial\Om}-g  \bigr).
\end{equation}
This map is well defined. 

Indeed, it is clear that $L_qu=\Delta^2u+qu\in L^p(\Omega).$ In particular, $4p>n$, and the Sobolev-Morrey embedding theorem yields 
\begin{equation}\label{eq:W4p-C1}
W^{4,p}(\Omega)\hookrightarrow C^{1,\beta}(\overline\Omega)\hookrightarrow L^\infty(\Omega)
\end{equation}
for some $\beta>0$. Hence there exists a constant
$C_{\rm emb}>0$ such that
\begin{equation}\label{eq:W4p-Linfty-estimate}
\|u\|_{L^\infty(\Omega)}\le C_{\rm emb}\|u\|_{W^{4,p}(\Omega)}.
\end{equation}
Let $\delta>0$ be such that $F$ is defined and smooth on $\overline\Omega\times(-\delta,\delta)$. Choose $\rho_0>0$ sufficiently small so that $C_{\rm emb}\rho_0<\delta/2$.
Then $u\in\mathcal U_{\rho_0}$ satisfies $|u(x)|<\delta/2$ on $\overline\Omega$. On the other hand, since $F$ is continuous on the compact set $\overline\Omega\times [-\delta/2,\delta/2]$, there exists $C_F>0$ such that
\begin{equation}
|F(x,s)|\le C_F,\quad x\in\overline\Omega,\quad |s|\le\delta/2.
\end{equation}
Hence, for any $u\in\mathcal U_{\rho_0}$,
\begin{equation}
\|F(\cdot,u(\cdot))\|_{L^p(\Omega)}\le|\Omega|^{1/p}C_F,
\end{equation}
implying that $F(\cdot,u(\cdot))\in L^p(\Omega).$
Thus the composition $x\mapsto F(x,u(x))$ is well defined.
Finally, the trace theorem gives $(u|_{\partial\Omega}, \Delta u|_{\partial\Omega})\in \mathcal X_{\mathbb R}$. 

We have the following (local) well-posedness for the nonlinear problem \eqref{eq:model}.
\begin{proposition}\label{pro:well-posedness-nonlinear}
\label{prop:local-wellposedness}
Let $\max\{1,n/3\}<p<\infty$, and let $q\in L^\infty(\Omega)$. Assume that Assumption {\bf (A)} holds, and that for some $\delta>0$,
\begin{equation}
F\in C^\infty\bigl(\overline\Omega\times(-\delta,\delta);\mathbb R\bigr),\quad
F(\cdot,0)=\partial_uF(\cdot,0)=0.
\end{equation} Then there exist $\varepsilon_0,\rho_0>0$
such that, for any $(f,g)\in B_{\mathcal X_{\mathbb R}}(0,\varepsilon_0),$ the nonlinear problem \eqref{eq:model} admits a unique solution
$u_{f,g}\in\mathcal U_{\rho_0}$. Moreover, the solution map $(f,g)\mapsto u_{f,g}$ is $C^\infty$ near the origin, and is real analytic if $F$ is real analytic.
\end{proposition}

\begin{proof}
We first study the regularity of the Nemytskii operator
\begin{equation}
\mathcal N_F:\mathcal U_{\rho_0}\rightarrow L^p(\Omega),\quad \mathcal N_F(u)\vcentcolon=F(\cdot,u(\cdot)).
\end{equation}
For any integer $m\ge1$, the formal $m$-th derivative is
\begin{equation}\label{eq:Nemytskii-mth-derivative}
D^m\mathcal N_F(u)[v_1,\ldots,v_m] =\partial_u^mF(\cdot,u)v_1\cdots v_m.
\end{equation}
We show that the right-hand side defines a bounded $m$-linear map
from $(W^{4,p}(\Omega))^m$ into $L^p(\Omega)$.

Since $\partial_u^mF$ is continuous on $\overline\Omega\times [-\delta/2,\delta/2]$, there exists $C_m>0$ such that
\begin{equation}
\sup_{x\in\overline\Omega,|s|\le\delta/2}|\partial_u^mF(x,s)|\le C_m.
\end{equation}
It then follows that
\begin{equation}\label{est:mutilinear}
\left\|\partial_u^mF(\cdot,u)v_1\cdots v_m\right\|_{L^p(\Omega)}\le C_m\|v_1\|_{L^p(\Omega)}\prod_{j=2}^m \|v_j\|_{L^\infty(\Omega)}\le C'_m\prod_{j=1}^m \|v_j\|_{W^{4,p}(\Omega)}.
\end{equation}
This gives
\begin{equation}
D^m\mathcal N_F(u)\in\mathcal L^m\bigl( W^{4,p}(\Omega);L^p(\Omega) \bigr).
\end{equation}
Furthermore, if $u_j\rightarrow u$ in $W^{4,p}(\Omega)$, then $u_j\rightarrow u$ uniformly on $\overline\Omega.$
Since $\partial_u^mF$ is uniformly continuous on the compact set
$\overline\Omega\times[-\delta/2,\delta/2]$, it follows that
\begin{equation}
\left\|\partial_u^mF(\cdot,u_j)-\partial_u^mF(\cdot,u)\right\|_{L^\infty(\Omega)}\rightarrow0.
\end{equation}
Combining this convergence with the multilinear estimate \eqref{est:mutilinear} gives
\begin{equation}
\left\| D^m\mathcal N_F(u_j)-D^m\mathcal N_F(u)\right\|_{\mathcal L^m}\rightarrow 0.
\end{equation}
Thus $\mathcal N_F\in C^\infty\bigl(\mathcal U_{\rho_0};L^p(\Omega)\bigr)$, and we conclude that $\mathcal F \in C^\infty\bigl(\mathcal U_{\rho_0}\times\mathcal X_{\mathbb R};\mathcal Y_p\bigr).$

We next consider the real-analytic case. Assume that $F$ is real
analytic in a neighborhood of $\overline\Omega\times\{0\}$. There exist constants
$M,R>0$ such that
\begin{equation}\label{eq:uniform-analytic-bound}
\sup_{x\in\overline\Omega}\frac{|\partial_u^kF(x,0)|}{k!}\le MR^{-k},\quad k=0,1,2,\ldots.
\end{equation}
For $k\ge1$, define the continuous $k$-homogeneous polynomial
\begin{equation}
 P_k(u):=a_k u^k,\quad a_k(x):=\frac{\partial_u^kF(x,0)}{k!}.
\end{equation}
Using again the embedding
$W^{4,p}(\Omega)\hookrightarrow L^\infty(\Omega)$, we can obtain
\begin{equation}
\|P_k(u)\|_{L^p(\Omega)}\le C M R^{-k} \|u\|_{W^{4,p}(\Omega)}^k.
\end{equation}
Thus, after decreasing $\rho_0$ if necessary,
the series
\begin{equation}\label{eq:Nemytskii-power-series}
\sum_{k=0}^{\infty}P_k(u) =\sum_{k=0}^{\infty}\frac{\partial_u^kF(\cdot,0)}{k!}u^k
\end{equation}
converges absolutely and locally uniformly in $L^p(\Omega)$ for $u\in\mathcal U_{\rho_0}$. Noting that the right-hand side of
\eqref{eq:Nemytskii-power-series} equals $F(\cdot,u)$, the map
$\mathcal N_F$ admits a convergent Banach-space power-series
representation. Hence $\mathcal N_F\in C^\omega\bigl(\mathcal U_{\rho_0};L^p(\Omega)\bigr).$
It follows that, in the analytic case, $\mathcal F\in C^\omega\bigl(\mathcal U_{\rho_0}\times\mathcal X_{\mathbb R};\mathcal Y_p\bigr).$

We now apply the implicit-function theorem. Using $F(x,0)=0,$ we have $\mathcal F(0,0,0)=0.$ For $v\in W^{4,p}(\Omega)$, we have
\begin{equation}
D_u\mathcal F(u,f,g)v= \bigl(L_qv-\partial_uF(\cdot,u)v,v|_{\partial\Omega},\Delta v|_{\partial\Omega}\bigr).
\end{equation}
In particular,
\begin{equation}
D_u\mathcal F(0,0,0)v=\bigl( L_qv-\partial_uF(\cdot,0)v, v|_{\partial\Omega},\Delta v|_{\partial\Omega}\bigr).
\end{equation}
By the condition $\partial_uF(x,0)=0,$ we obtain
\begin{equation}\label{eq:linearized-F-Aq}
D_u\mathcal F(0,0,0)v=\bigl( L_qv,v|_{\partial\Omega},\Delta v|_{\partial\Omega}\bigr)=\mathcal A_qv.
\end{equation}
By Lemma~\ref{lem:full-navier-isomorphism}, $\mathcal A_q: W^{4,p}(\Omega)\rightarrow\mathcal Y_p$ is an isomorphism. Therefore the smooth
implicit-function theorem applied to $\mathcal F$ at $(0,0,0)$ yields
neighborhoods $\mathcal V\subset W^{4,p}(\Omega)$ and $\mathcal U\subset \mathcal X_{\mathbb R}$ of the origin and a unique smooth map $\mathcal S_{q,F}: \mathcal U \rightarrow \mathcal V$ such that
\begin{equation}
\mathcal F\bigl(\mathcal S_{q,F}(f,g),f,g\bigr)=0,\quad (f,g)\in\mathcal U.
\end{equation}
This implies $u_{f,g}:=\mathcal S_{q,F}(f,g)$ is the unique solution in $\mathcal V$ of \eqref{eq:model}.

After shrinking the neighborhoods, if necessary, we can write
\begin{equation}
B_{\mathcal X_{\mathbb R}}(0,\varepsilon_0)\subset\mathcal U,\quad \mathcal U_{\rho_0} \subset\mathcal V
\end{equation}
for some $\varepsilon_0,\rho_0>0$. The uniqueness follows from the local uniqueness in the fixed neighborhood $\mathcal V$ of the zero solution.
If $F$ is real analytic, then the analytic implicit-function theorem shows that $\mathcal S_{q,F}:\mathcal U\rightarrow\mathcal V$ is real analytic.

Finally, because $p>n/3$, the Sobolev embedding $W^{4,p}(\Omega)\hookrightarrow C^{1,\beta}(\overline\Omega)$ implies that the map
\begin{equation}
\ell_{x_0}: W^{4,p}(\Omega)\rightarrow\mathbb R,\quad \ell_{x_0}(u):=\partial_\nu u(x_0),
\end{equation}
is a continuous linear functional. 
The nonlinear point-measurement map can therefore be written as $\Lambda_{q,F}=\ell_{x_0}\circ\mathcal S_{q,F}.$
Therefore, $\Lambda_{q,F}$ is $C^\infty$ near the zero boundary data, and it is real analytic whenever $F$ is real analytic.
\end{proof}

\subsection{Higher-order linearization}
For $m\ge 1$, let $H_j=(f_j,g_j)\in\mathcal X_{\mathbb R}$. For $\varepsilon=(\varepsilon_1,\cdots,\varepsilon_m)\in\mathbb R^m$, we define $u_\varepsilon=u_{\varepsilon_1H_1+\cdots\varepsilon_mH_m}$. 
For a nonempty set $\mathcal I\subset\{1,\dots,m\}$, we set the mixed variation 
\begin{equation}
 u_{\mathcal I}:=\Big(\prod_{j\in\mathcal I}\partial_{\varepsilon_j}\Big)u_\varepsilon\Big|_{\varepsilon=0}.
\end{equation}
For $\mathcal I=\{j\}$ we write $u_{\mathcal I}=v_j$.  Differentiating the equation \eqref{eq:model} with respect to $\varepsilon$ gives
\begin{equation}\label{eq:first-linear-vj}
 \begin{cases}
 L_qv_j=0&\text{in }\Om,\\
 v_j=f_j,\, \Delta v_j=g_j&\text{on }\partial\Om.
 \end{cases}
\end{equation}
Whenever $|\mathcal I|\ge2$, we have    
\begin{equation}\label{eq:higher-hom-bc}
u_{\mathcal I}=\Delta u_{\mathcal I}=0\quad\text{on }\partial\Om.
\end{equation}
Moreover, the $m$-th Fr\'echet derivative of $\Lambda_{q,F}$ forms a continuous, symmetric $m$-linear map from $\mathcal X^m_{\mathbb R}\to \mathbb R$ such that 
\begin{equation}
D^m\Lambda_{q,F}(0)[H_1,\cdots, H_m]=\partial^m_{\varepsilon_1,\cdots\varepsilon_m}\Lambda_{q,F}(\varepsilon_1H_1+\cdots+\varepsilon_mH_m)\big|_{\varepsilon=0}.
\end{equation}
In particular, for the linear equation \eqref{eq:first-linear-vj}, we have $D\Lambda_{q,F}(0)[H]=\partial_\nu v_j(x_0)$, $j=1,\cdots,m$.

Let $\Pi(\mathcal I)$ denote the set of all set partitions of $\mathcal I$.  The multivariable Fa\`a di Bruno formula gives
\begin{equation}\label{eq:partition}
L_qu_{\mathcal I}=\sum_{\substack{\pi\in\Pi(\mathcal I), |\pi|\ge2}}\partial_u^{|\pi|}F(x,0)\prod_{\alpha\in\pi}u_\alpha, \quad |\mathcal I|\ge2.
\end{equation}
Suppose that $m_*\ge2$ and
\begin{equation}\label{eq:first-nonlinear-order-assumption}
 \partial_u^kF(x,0)=0\quad(2\le k<m_*),
 \quad \partial_u^{m_*}F(x,0)=\lambda_* A(x),
 \quad \lambda_*\neq0.
\end{equation}
By \eqref{eq:partition} and using assumption {\bf (A)}, we can get
\begin{equation}\label{eq:lower-mixed-vanish}
 u_{\mathcal I}=0, \quad 2\le |\mathcal I|<m_*.
\end{equation}
At order $m_*$, with $\mathcal I_*=\{1,\cdots,m_*\}$, using $\partial_u^{m_*}F(x,0)=\lambda_* A(x),$ and $u_{\{j\}}=v_j$, we obtain
\begin{equation}\label{eq:first-nonlinear-variation}
\begin{cases}
L_qu_{\mathcal I_*}=\lambda_* A(x)v_1\cdots v_{m_*} &\text{in }\Omega,\\
u_{\mathcal I_*}=\Delta u_{\mathcal I_*}=0 &\text{on }\partial\Omega.
\end{cases}
\end{equation}

\begin{remark}\label{rem:complexification}
Although the nonlinear forward problem and the measurement map are defined for real-valued boundary data, we shall use
complex-valued solutions of the linearized equation in next sections. This causes no essential difficulty. For each $m\ge1$, the real continuous multilinear form $T^m=D^m\Lambda_{q,F}(0)$ has a unique complex multilinear extension to the complexification of $\mathcal X^m$ via the formula
\begin{equation}
T^m_{\mathbb C}(H_1,\cdots, H_m)=\sum_{\alpha\in\{0,1\}^m}i^{|\alpha|}T(Z_{1,\alpha_1},\cdots, Z_{m,\alpha_m}),\quad H_j=X_j+i Y_j,
\end{equation}
where $Z_{j,0}=X_j$ and $Z_{j,1}=Y_j$ for $j=1,\cdots,m$.
Hence, after complexification, the higher-order linearization identities, including \eqref{eq:first-linear-vj} and \eqref{eq:first-nonlinear-variation}, remain valid for complex-valued boundary data and the corresponding complex-valued solutions of the linearized equation.
\end{remark}

\section{Adjoint representation of the point measurement}\label{sec:adjoint}
In this section, we first represent the point measurement at $x_0\in\partial\Omega$ by means of a suitable adjoint function $\Psi_q\in L^{p'}(\Omega)$, where $p'$ is the H\"older conjugate of $p$.  The resulting Green identity expresses the point observation in terms of $\Psi_q$ and boundary traces.  Moreover, $\Psi_q$ is regular away from the fixed $x_0$, which allows us to carry out the subsequent boundary analysis on boundary patches not containing the observation point. We also present a (local) unique continuation property of the fourth-order differential operator. Some useful properties of the adjoint function $\Psi_q$ will be introduced in the end of this section.

We emphasize that, in the remainder of this section, unless otherwise stated, the potential $q$ is assumed to satisfy Assumption {\bf (A)}, and the nonlinearity $F$ is smooth near $\overline{\Omega}\times\{0\}$, with $F(x,0)=\partial_uF(x,0)=0$ in $\overline\Omega$.

We introduce a useful integral identity that will be frequently used.
\begin{lemma}
For $u\in W^{4,r}(\Omega)$ and $v\in W^{4,r'}(\Omega)$ with $1<r<\infty$ and $r'=\frac{r}{r-1}$, it holds that
\begin{equation}\label{eq:biharmonic-green}
\int_\Om[(\Delta^2u)v-u\Delta^2v]\,dx=\int_{\partial\Om}\bigl[\partial_\nu\Delta u\,v-\Delta u\,\partial_\nu v+\Delta v\,\partial_\nu u-u\,\partial_\nu\Delta v\bigr]\,dS.
\end{equation}
\end{lemma}

\begin{proof}
Suppose first that $u,v\in C^\infty(\overline\Omega)$.  Applying
Green's formula twice gives
\begin{align} 
\int_\Omega(\Delta^2u)v\,dx&=\int_\Omega\Delta u\,\Delta v\,dx+\int_{\partial\Omega}\bigl((\partial_\nu\Delta u)v
-(\Delta u)\partial_\nu v\bigr)\,dS,\label{eq:green-first-term}
\\
\int_\Omega u\Delta^2v\,dx&=\int_\Omega\Delta u\,\Delta v\,dx+\int_{\partial\Omega}\bigl(u\partial_\nu\Delta v
-(\Delta v)\partial_\nu u\bigr)\,dS.
\label{eq:green-second-term}
\end{align}
Subtracting \eqref{eq:green-second-term} from
\eqref{eq:green-first-term} proves
\eqref{eq:biharmonic-green} for smooth functions.

For general
$u\in W^{4,r}(\Omega)$, the trace theorem gives
\begin{equation}\label{eq:biharmonic-traces}
\bigl(u|_{\partial\Omega},\partial_\nu u|_{\partial\Omega},\Delta u|_{\partial\Omega},\partial_\nu\Delta u|_{\partial\Omega}\bigr)\in W^{4-1/r,r}(\partial\Omega)\times W^{3-1/r,r}(\partial\Omega)\times W^{2-1/r,r}(\partial\Omega)\times W^{1-1/r,r}(\partial\Omega).
\end{equation}
The corresponding traces of $v$ belong to the analogous spaces with
exponent $r'$.  In particular, the traces of $u$ belong to
$L^r(\partial\Omega)$ and those of $v$ belong to
$L^{r'}(\partial\Omega)$.  Hölder's inequality therefore shows that every product in the boundary integral is integrable and that the boundary form is continuous on $W^{4,r}(\Omega)\times W^{4,r'}(\Omega).$ The volume terms are continuous on the same product space.  The result now follows by approximating $u$ and $v$ in their respective Sobolev spaces by functions in $C^\infty(\overline\Omega)$ and passing to the limit.
\end{proof}

\subsection{Some key identities}  We will establish a key identity that connects the point measurement with a suitable adjoint function. Applying Lemma~\ref{lem:full-navier-isomorphism}, we have the following lemma:
\begin{lemma}\label{lem:Green-identity}
For any $1<r<\infty$, the problem
\begin{equation}
 L_qw=f\,\,  {\rm in}\,\, \Om, \quad w=\Delta w=0\,\, {\rm on }\,\,  \partial\Om
\end{equation}
defines a bounded (inverse) map
\begin{equation}\label{eq:Lp-resolvent}
G_q:L^r(\Om)\rightarrow W^{4,r}(\Om), \quad \|G_qf\|_{W^{4,r}(\Omega)}\le C_r\|f\|_{L^r(\Omega)}.
\end{equation}
\end{lemma}

We next establish a local boundary regularity estimate for solutions to $L_qu=0$ that will be applied to the regularized adjoint functions.

\begin{lemma}\label{lem:local-boundary-regularity}
Let $1<s<\infty$, and let
$\mathcal O_1,\mathcal O_2\subset\mathbb R^n$ be open neighborhoods
of a portion of $\partial\Omega$ satisfying $\overline{\mathcal O}_1\subset \mathcal O_2.$ Assume that $q$ is smooth on $\mathcal O_2\cap\overline\Omega$, and that $u\in W^{4,s}(\mathcal O_2\cap\Omega)$ satisfies the problem
\begin{equation}
\begin{cases}
L_qu=0 & {\rm in }\,\, \mathcal O_2\cap\Omega,\\
u=\Delta u=0 & {\rm on }\,\, \mathcal O_2\cap\partial\Omega.
\end{cases}
\end{equation}
Then, for any integer $m\ge0$,
\begin{equation}\label{eq:local-boundary-estimate}
\|u\|_{W^{m+4,s}(\mathcal O_1\cap\Omega)}\le C_{m,s}\|u\|_{L^s(\mathcal O_2\cap\Omega)}, \quad m\ge0.
\end{equation}
The constant $C_{m,s}$ depends on
$\mathcal O_1$, $\mathcal O_2$, $m$, $s$, $\Omega$, and on finitely many
derivatives of $q$ in $\mathcal O_2\cap\overline\Omega$.
\end{lemma}
The proof is deferred to the appendix.

We next establish a key identity associated with the point measurement.
\begin{lemma}\label{lem:point-adjoint}
Let $\max\{1,n/3\}<p<\infty$ and $p'=\frac{p}{p-1}$, and let $f\in L^p(\Omega)$. Assume that $q\in L^\infty(\Omega)$ and that Assumption {\bf (A)} holds.  Then there exists a unique
$\Psi_q\in L^{p'}(\Omega)$ such that
\begin{equation}\label{eq:adjoint-duality}
\int_\Omega\Psi_q f\,dx=\partial_\nu(G_qf)(x_0).
\end{equation}
If $q$ is real valued, then $\Psi_q$ is real valued.  Moreover,
\begin{equation}\label{eq:adjoint-interior}
 L_q\Psi_q=0 \quad {\rm in}\,\,  \mathcal D'(\Omega).
\end{equation}
If $q$ is $C^\infty$ in a boundary neighborhood of a relatively open set $U\Subset\partial\Omega\setminus\{x_0\},$
then, for a sufficiently small open  set $\mathcal O_U$ of $U$, it holds that
\begin{equation}\label{eq:adjoint-local-smooth}
\Psi_q\in C^\infty(\overline\Omega\cap\mathcal O_U),\quad\Psi_q=\Delta\Psi_q=0 \quad {\rm on}\,\, U.
\end{equation}
\end{lemma}

\begin{proof}
We divide the proof into two steps.

\emph{Step 1.} Define a map
\begin{equation}
\ell_q:L^p(\Om;\mathbb R)\rightarrow\mathbb R,\quad \ell_q(f):=\partial_\nu(G_qf)(x_0).
\end{equation}
By \eqref{eq:Lp-resolvent} and \eqref{eq:W4p-C1}, we have
\begin{equation}\label{eq:point-functional-bound}
|\ell_q(f)|=|\partial_\nu(G_qf)(x_0)|\le C\|G_qf\|_{W^{4,p}(\Om)}\le C\|f\|_{L^p(\Om)}.
\end{equation}
Thus $\ell_q$ is a bounded linear functional on
$L^p(\Om;\mathbb R)$.  Since $1<p<\infty$ and $(L^p(\Om;\mathbb R))'=L^{p'}(\Om;\mathbb R),$ there exists a unique real-valued
$\Psi_q\in L^{p'}(\Om)$ such that
\begin{equation}\label{eq:adjoint-duality-real}
\int_\Om\Psi_qf\,dx=\partial_\nu(G_qf)(x_0),\quad f\in L^p(\Om;\mathbb R).
\end{equation}
Extending both sides complex linearly gives
\eqref{eq:adjoint-duality} for complex-valued $f$.  In particular,
the representative $\Psi_q$ remains real valued.

We next prove that $\Psi_q$ satisfies the equation $L_q\Psi_q=0$.  Let $\varphi\in C_c^\infty(\Om).$
Then $\varphi$ has homogeneous Navier boundary data and $G_q(L_q\varphi)=\varphi$.  Using \eqref{eq:adjoint-duality}, we can obtain
\begin{equation}
\int_\Om\Psi_qL_q\varphi\,dx=\partial_\nu\varphi(x_0)=0.
\end{equation}
It then follows that
\begin{equation}
\langle L_q\Psi_q,\varphi\rangle_{\mathcal D'(\Omega),\mathcal D(\Omega)}:=\int_\Om\Psi_qL_q\varphi\,dx=0.
\end{equation}
Since $\varphi\in C_c^\infty(\Om)$ is arbitrary, this implies $L_q\Psi_q=0$ in $\mathcal D'(\Om).$

\emph{Step 2.} It remains to prove smoothness up to the boundary away from $x_0$ of $\Psi_q$.
Let $\varepsilon>0$. We choose a family of smooth boundary functions $\rho_\varepsilon\in C^\infty(\partial\Om),$ satisfying $\rho_\varepsilon\ge0,$ and
\begin{equation}\label{eq:boundary-approximate-identity}
\int_{\partial\Om}\rho_\varepsilon\,dS=1, \quad \operatorname{supp}\rho_\varepsilon \rightarrow\{x_0\},
\quad \rho_\varepsilon \rightarrow\delta_{x_0}\quad\text{in }\mathcal D'(\partial\Om),\quad \text{as } \varepsilon\to 0.
\end{equation}
Here $\delta_{x_0}$ denotes the Dirac distribution on $\partial\Omega$, defined by
\begin{equation}
\langle\delta_{x_0},\phi\rangle=\phi(x_0),\quad \phi\in C^\infty(\partial\Omega).
\end{equation}
By Lemma~\ref{lem:full-navier-isomorphism}, applied with exponent
$p'$, there exists a unique
$\Psi_\varepsilon\in W^{4,p'}(\Om)$ satisfying
\begin{equation}\label{eq:smooth-adjoint-approx}
\begin{cases}
L_q\Psi_\varepsilon=0 &\text{in }\Om,\\
\Psi_\varepsilon=0,\, \Delta\Psi_\varepsilon=\rho_\varepsilon&\text{on }\partial\Om.
\end{cases}
\end{equation}

Let $f\in L^p(\Om)$ and set $W:=G_qf.$
Then
\begin{equation}
L_qW=f \quad\text{in }\Om, \quad W=\Delta W=0 \quad\text{on }\partial\Om.
\end{equation}
Applying the Green formula \eqref{eq:biharmonic-green} for $\Delta^2$ to $\Psi_\varepsilon$ and $W$, we obtain
\begin{equation}\label{eq:approx-duality}
\int_\Om\Psi_\varepsilon f\,dx=\int_{\partial\Om}\rho_\varepsilon\partial_\nu W\,dS.
\end{equation}
The boundary integral in \eqref{eq:approx-duality} is classical because $ W\in W^{4,p}(\Om)\hookrightarrow C^{1,\beta}(\overline\Om)$ for some $\beta>0$. Using $\|\rho_\varepsilon\|_{L^1(\partial\Om)}=1$ and \eqref{eq:W4p-C1}, we find that
\begin{equation}
\Big|\int_{\partial\Om}\rho_\varepsilon\partial_\nu W\,dS\Big|\le\|\partial_\nu W\|_{L^\infty(\partial\Om)}
\le C\|W\|_{W^{4,p}(\Om)} \le C\|f\|_{L^p(\Om)}.
\end{equation}
It follows from \eqref{eq:approx-duality} and
$L^p$--$L^{p'}$ duality that
\begin{equation}\label{eq:adjoint-uniform-Lp}
 \|\Psi_\varepsilon\|_{L^{p'}(\Om)}
 \le C,
\end{equation}
where the constant is independent of $\varepsilon$.

For any fixed $f\in L^p(\Om)$, the function $\partial_\nu(G_qf)$ is continuous on $\partial\Om$.  Hence the defining properties of $\rho_\varepsilon$ imply that
\begin{equation}
\int_{\partial\Om}\rho_\varepsilon\partial_\nu(G_qf)\,dS\rightarrow\partial_\nu(G_qf)(x_0),\quad \varepsilon\to 0.
\end{equation}
Indeed, writing $g_f=\partial_\nu(G_qf)$ and using
$\int_{\partial\Omega}\rho_\varepsilon\,dS=1$, we have
\begin{equation}
\Big|\int_{\partial\Omega}\rho_\varepsilon g_f\,dS-g_f(x_0) \Big|\le\sup_{x\in\operatorname{supp}\rho_\varepsilon}|g_f(x)-g_f(x_0)|.
\end{equation}
Since $g_f\in C^{0,\beta}(\partial\Omega)$ and $\operatorname{supp}\rho_\varepsilon \subset B_{\partial\Omega}(x_0,C_0\varepsilon)=\{x\in\partial\Omega: |x-x_0|<C_0\varepsilon\}$, the right-hand side
tends to zero as $\varepsilon\to0$.
Combining this with \eqref{eq:approx-duality} and
\eqref{eq:adjoint-duality}, we obtain
\begin{equation}
\int_\Om(\Psi_\varepsilon-\Psi_q)f\,dx\rightarrow0\quad \text{for any }f\in L^p(\Om),\quad \varepsilon\to 0.
\end{equation}
Therefore, we can derive
\begin{equation}\label{eq:adjoint-weak-convergence}
\Psi_\varepsilon\rightharpoonup\Psi_q \quad\text{weakly in }L^{p'}(\Om).
\end{equation}

Fix a relatively open set
$U\Subset\partial\Omega\setminus\{x_0\}$.
Choose boundary neighborhoods $\mathcal O_0\Subset\mathcal O_1\Subset\mathcal O_2$ of $\overline U$ such that $\overline U\subset\mathcal O_0\cap\partial\Omega$, $x_0\notin\overline{\mathcal O}_2$ and $q$ is $C^\infty$ up to the boundary in $\mathcal O_2\cap\overline\Omega$. For all sufficiently small $\varepsilon$, $\operatorname{supp}\rho_\varepsilon\cap\bigl(\mathcal O_2\cap\partial\Om\bigr)=\varnothing.$ It then follows from \eqref{eq:smooth-adjoint-approx} that
\begin{equation}\label{eq:adjoint-local-homogeneous}
L_q\Psi_\varepsilon=0 \quad\text{in }\mathcal O_2\cap\Om,\quad \Psi_\varepsilon=\Delta\Psi_\varepsilon=0 
\quad\text{on }\mathcal O_2\cap\partial\Om.
\end{equation}
Applying Lemma~\ref{lem:local-boundary-regularity} with exponent
$s=p'$ and using \eqref{eq:adjoint-uniform-Lp}, we obtain, for every integer $k\ge0$,
\begin{equation}\label{eq:adjoint-local-Wkp}
\|\Psi_\varepsilon\|_{W^{k,p'}(\mathcal O_1\cap\Om)}\le C_k\|\Psi_\varepsilon\|_{L^{p'}(\mathcal O_2\cap\Om)}\le C_k.
\end{equation}

Fix an integer $r\ge0$ and choose $ k>r+\frac{n}{p'}.$ By the local Rellich compactness theorem and the boundary Sobolev
embedding, the restriction map
\begin{equation}
W^{k,p'}(\mathcal O_1\cap\Om)\rightarrow C^r(\overline\Om\cap\mathcal O_0)
\end{equation}
is compact.  It follows that every sequence $\varepsilon_j\to0$ admits a subsequence, denoted again by $\varepsilon_j$, and a function $\Phi\in C^r(\overline\Omega\cap\mathcal O_0)$ such that
\begin{equation}
\Psi_{\varepsilon_j}\rightarrow\Phi\quad\text{in }  C^r(\overline\Omega\cap\mathcal O_0).
\end{equation}
In particular, this subsequence converges to $\Phi$ in $\mathcal D'(\mathcal O_0\cap\Omega)$.  On the other hand, \eqref{eq:adjoint-weak-convergence} implies that the same subsequence converges to $\Psi_q$ in $\mathcal D'(\mathcal O_0\cap\Omega)$.  By uniqueness of distributional limits, we have $\Phi=\Psi_q$ in $\mathcal O_0\cap\Omega.$
Thus every $C^r$-convergent subsequence has the same limit $\Psi_q$.
Since the set of restrictions $\{\Psi_\varepsilon|_{\overline\Omega\cap\mathcal O_0}:0<\varepsilon<\varepsilon_0\}$, where $\varepsilon_0$ is a sufficiently small constant, is relatively compact in
$C^r(\overline\Omega\cap\mathcal O_0)$, a standard subsequence argument yields $\Psi_\varepsilon\rightarrow\Psi_q$ in $C^r(\overline\Omega\cap\mathcal O_0).$ Hence,
\begin{equation}\label{eq:adjoint-local-Cr}
\Psi_\varepsilon\rightarrow \Psi_q \quad\text{in } C^r(\overline\Om\cap\mathcal O_0).
\end{equation}
Since $r$ is arbitrary, we conclude that $\Psi_q\in C^\infty(\overline\Om\cap\mathcal O_0).$

Finally, take $r\ge2$ in \eqref{eq:adjoint-local-Cr}.  For all sufficiently small $\varepsilon$, the support condition on
$\rho_\varepsilon$ and \eqref{eq:smooth-adjoint-approx} give $\Psi_\varepsilon=\Delta\Psi_\varepsilon=0$ on $U$.
Passing to the limit in $C^2$ yields $\Psi_q=\Delta\Psi_q=0$ on $U$.
Taking $\mathcal O_U:=\mathcal O_0$ proves the assertion \eqref{eq:adjoint-local-smooth} as desired.
\end{proof}

Based on Lemma~\ref{lem:point-adjoint}, we can immediately conclude the following corollary:
\begin{corollary}\label{cor:transposition}
Under the conditions of Lemma~\ref{lem:point-adjoint}, if $W\in W^{4,p}(\Om)$ satisfies $ W=\Delta W=0$ on $\partial\Omega$, 
then
\begin{equation}\label{eq:transposition}
\int_\Om\Psi_qL_qW\,dx=\partial_\nu W(x_0).
\end{equation}
Suppose, in addition, that the condition \eqref{eq:first-nonlinear-order-assumption} holds.  For the boundary data $H_j=(f_j,g_j)\in\mathcal X_{\mathbb C}$, let $v_j$ solve the following linear equation
\begin{equation}
 \begin{cases}
L_qv_j=0 & {\rm in}\,\, \Om,\\
v_j=f_j,\, \Delta v_j=g_j & {\rm on}\,\, \partial\Omega.
 \end{cases}
\end{equation}
Then after complexification, it holds that
\begin{equation}\label{eq:first-nonlinear-tensor}
D^{m_*}\Lambda_{q,F}(0)[H_1,\dots,H_{m_*}]=\lambda_*\int_\Om\Psi_q A(x)v_1\cdots v_{m_*}\,dx.
\end{equation}
\end{corollary}

\begin{proof}
Let $ f:=L_qW.$ The uniqueness of solutions gives $G_qf=G_q(L_qW)=W.$ \eqref{eq:transposition} directly follows from  \eqref{eq:adjoint-duality} by letting $f=L_qW$.

We next prove \eqref{eq:first-nonlinear-tensor} for $H_1,\dots,H_{m_*}\in\mathcal X_{\mathbb R}.$
Recall that $u_{\mathcal I_*}$ satisfy the problem \eqref{eq:first-nonlinear-variation}.
Moreover,
\begin{equation}
D^{m_*}\Lambda_{q,F}(0)[H_1,\dots,H_{m_*}]= \partial_\nu u_{\mathcal I_*}(x_0).
\end{equation}
Applying \eqref{eq:transposition} with $W=u_{\mathcal I_*},$ we can get
\begin{equation}
D^{m_*}\Lambda_{q,F}(0)[H_1,\dots,H_{m_*}]=\int_\Om\Psi_qL_qu_{\mathcal I_*}\,dx=\lambda_*
\int_\Om\Psi_q A(x) v_1\cdots v_{m_*}\,dx.
\end{equation}
Both sides are continuous real multilinear functions of $H_1,\dots,H_{m_*}$.  Their unique complex multilinear extensions
therefore agree on $\mathcal X_{\mathbb C}^{m_*}$. Hence \eqref{eq:first-nonlinear-tensor} follows.
\end{proof}

For later use, we denote by $\mathcal Q^\infty_\Gamma(\Omega)$ the
class of potentials $q\in L^\infty(\Omega;\mathbb R)$ for which there
exist an open neighborhood $\mathcal N_\Gamma\subset\mathbb R^n$ of
$\overline\Gamma$ and a function
$\widetilde q\in C^\infty(\mathcal N_\Gamma;\mathbb R)$ such that
\begin{equation}\label{eq:Q-Gamma-infinity}
 \widetilde q=q
 \quad\text{a.e. in }\mathcal N_\Gamma\cap\overline\Omega.
\end{equation}
Thus $\mathcal Q^\infty_\Gamma(\Omega)$ consists of bounded potentials
that admit a smooth extension across the boundary near
$\overline\Gamma$.

The following lemma gives the localized first-order linearization integral identity.
\begin{lemma}\label{lem:first-green}
Assume that $q\in \mathcal Q^\infty_\Gamma(\Omega)$.
Let $U\Subset\Gamma$, and let
$h_0,h_1\in C_c^\infty(U)$. Let $v$ satisfy
\begin{equation}
\begin{cases}
L_qv=0 & {\rm in}\,\, \Omega,\\
v=h_0,\, \Delta v=h_1 & {\rm on}\,\, \partial\Omega.
\end{cases}
\end{equation}
Then it holds that
\begin{equation}\label{eq:first-green}
D\Lambda_{q,F}(0)[(h_0,h_1)]=\int_U h_1\,\partial_\nu\Psi_q\,dS+\int_U h_0\,\partial_\nu\Delta\Psi_q\,dS.
\end{equation}
\end{lemma}

\begin{proof}
Recalling the first-order linearization of $\Lambda_{q,F}$, we have
\begin{equation}\label{eq:first-variation-measurement}
D\Lambda_{q,F}(0)[(h_0,h_1)]=\partial_\nu v(x_0).
\end{equation}
Let $\rho_\varepsilon\in C^\infty(\partial\Omega)$ denote the same function as in \emph{Step 2} of the proof of Lemma~\ref{lem:point-adjoint}, and let $\Psi_\varepsilon$ solve the problem
\begin{equation}
 \begin{cases}
  L_q\Psi_\varepsilon=0 & \text{in }\Omega,\\
  \Psi_\varepsilon=0,\,
  \Delta\Psi_\varepsilon=\rho_\varepsilon
  & \text{on }\partial\Omega.
 \end{cases}
\end{equation}
Applying the Green formula \eqref{eq:biharmonic-green}, we can get
\begin{equation}\label{eq:first-green-epsilon}
\int_{\partial\Omega}\rho_\varepsilon\,\partial_\nu v\,dS=\int_U h_1\,\partial_\nu\Psi_\varepsilon\,dS+
\int_U h_0\,\partial_\nu\Delta\Psi_\varepsilon\,dS.
\end{equation}
Using $\partial_\nu v\in C(\partial\Omega)$, the properties of the approximate identity therefore imply that
\begin{equation}
\int_{\partial\Omega}\rho_\varepsilon\,\partial_\nu v\,dS\rightarrow\partial_\nu v(x_0).
\end{equation}
Moreover, $U\Subset\Gamma$ is separated from $x_0$, and the local boundary regularity established in Lemma~\ref{lem:point-adjoint} gives $\Psi_\varepsilon\rightarrow\Psi_q$ in $C^m(\overline\Omega\cap\mathcal O_U)$, for any $m\geq0$, after choosing a sufficiently small boundary neighborhood $\mathcal O_U$ of $\overline U$.  In particular,
\begin{equation}
\partial_\nu\Psi_\varepsilon\rightarrow\partial_\nu\Psi_q, \quad
\partial_\nu\Delta\Psi_\varepsilon\rightarrow\partial_\nu\Delta\Psi_q \quad\text{uniformly on }\overline U.
\end{equation}
Hence, \eqref{eq:first-green} follows from passing to the limit in \eqref{eq:first-green-epsilon} and using \eqref{eq:first-variation-measurement}.
\end{proof}

\begin{remark}
If the first-order linearization of two measurement maps $\Lambda_{q_1,F_1}$ and $\Lambda_{q_2,F_2}$ coincide, applying \eqref{eq:first-green} first with $h_0=0$ and then with $h_1=0$ gives
\begin{equation}\label{eq:trace-equality}
\partial_\nu\Psi_1=\partial_\nu\Psi_2,\quad\partial_\nu\Delta\Psi_1=\partial_\nu\Delta\Psi_2 \quad{\rm on}\,\, \Gamma.
\end{equation}
\end{remark}

\subsection{A unique continuation property}
We will introduce a boundary unique continuation property (UCP) for the coupled problem
\begin{equation}
\Delta u=w,\quad \Delta w=-cu.
\end{equation}
We refer, for instance, to \cite[Thm.~17.2.1]{HormanderIII} and \cite[Thm.~3.3.1]{isakov2017inverse} for UCP arguments for second-order elliptic equations. It is well known that Carleman estimates play an essential role in deriving the UCP for various PDEs (see, e.g., \cite{leRousseau2022ellipticVol1,leRousseau2022ellipticVol2}).

\begin{lemma}\label{lem:boundary-ucp-system}
Let $U_*\subset\partial\Omega$ be a nonempty relatively open set.
Assume that there exists an open set
$\mathcal N_*\subset\mathbb R^n$ containing $U_*$ such that $c\in L^\infty(\Omega\cap\mathcal N_*).$
Suppose that $u,w\in W^{2,2}_{\rm loc}$ up to $U_*$ and satisfy the equation
\begin{equation}\label{eq:ucp-system}
\Delta u=w,\quad \Delta w=-cu \quad{\rm in}\,\, \Omega\cap\mathcal N_*,\quad u=\partial_\nu u=w=\partial_\nu w=0\quad {\rm on}\,\, U_*.
\end{equation}
Then, for any $x_*\in U_*$, there exists an open neighborhood
$\mathcal V_*\subset\mathbb R^n$ of $x_*$ such that
\begin{equation}\label{local-vanishing-u}
 u=w=0 \quad {\rm in}\,\, \Omega\cap\mathcal V_*.
\end{equation}
If, in addition, \eqref{eq:ucp-system} holds in $\Omega$ with $c\in L^\infty(\Omega)$, then $ u=w=0$ in $\Omega$.
\end{lemma}
The proof of this lemma is provided in the appendix.

\subsection{Properties of the adjoint function $\Psi_q$}\label{subsec:pro-adjoint-Psi}
In the end of this section, we will present some properties of the adjoint function $\Psi_q$, which will be used in the analysis of inverse problems. On $\Gamma$, we set $A_q:=\partial_\nu\Psi_q$ and $B_q:=\partial_\nu\Delta\Psi_q.$ 

We first recall the boundary normal coordinates that are useful for the boundary determination of coefficients.  Let $V\Subset\partial\Omega$ be contained in a boundary coordinate
neighborhood.  Choose an open set
$\mathcal Y\subset\mathbb R^{n-1}$ and a smooth parametrization $\gamma:\mathcal Y\rightarrow V$
of the boundary.  Thus $y\in\mathcal Y$ represents the boundary point
$\gamma(y)\in V$.  For sufficiently small $\varepsilon>0$, the map
\begin{equation}
\Phi(y,t)=\gamma(y)-t\nu(\gamma(y)),\quad 0\leq t<\varepsilon,
\end{equation}
defines inward normal coordinates near $V$. Thus $y$ denotes the tangential variables and $t$ is the distance to the boundary measured in the inward normal direction.  In particular,
\begin{equation}\label{eq:normal-coordinate-sign}
t=0\quad\text{on }\partial\Omega,\quad t>0\quad\text{in }\Omega,\quad
\partial_\nu=-\partial_t\quad\text{at }t=0.
\end{equation}
We write $u(y,t)$ for $u(\Phi(y,t))$. In these coordinates, the Euclidean metric takes the form
$ds^2=dt^2+g_t$, and hence
\begin{equation}\label{eq:laplacian-boundary-normal-form}
\Delta=\partial_t^2+a(y,t)\partial_t+\Delta_{g_t},
\end{equation}
where $a$ is a certain smooth function and $\Delta_{g_t}$ is the Laplace--Beltrami operator in the tangential variables associated with the induced metric $g_t$ on the hypersurface $t=\mathrm{constant}$. More precisely, let $J=\sqrt{\det g_t(y)}$. We have $a(y,t)=\partial_t\log J(y,t)$.

The following lemma reveals the nondegeneracy of the boundary traces $A_q,B_q$.

\begin{lemma}\label{lem:AB}
Assume that $q\in \mathcal Q^\infty_\Gamma(\Omega)$. Then there is no nonempty relatively open subset of $\Gamma$ on which both $A_q$ and $B_q$ vanish identically, and there exist a nonempty connected relatively open set $U\Subset\Gamma$, an integer $\kappa\in\{1,3\}$, a constant $c_0>0$, and a function $\psi_\kappa\in C^\infty(\overline U)$ such that
\begin{equation}\label{eq:nonvanishing-psi-kappa}
\left|\psi_\kappa(y)\right|\geq c_0,\quad y\in\overline U.
\end{equation}
Moreover, in inward normal coordinates $(y,t)$, where $t>0$ in $\Omega$, there exist $\varepsilon>0$ and $R_\kappa\in C^\infty(\overline U\times[0,\varepsilon])$ such that
\begin{equation}\label{eq:kappa-expansion}
\Psi_q(y,t)=\frac{\psi_\kappa(y)}{\kappa!}\,t^\kappa+t^{\kappa+1}R_\kappa(y,t),\quad
(y,t)\in\overline U\times[0,\varepsilon].
\end{equation}
\end{lemma}

\begin{proof}
We divide the proof into two steps.

\emph{Step 1.} Suppose that $A_q=B_q=0$ on a nonempty open subset $U_1\Subset\Gamma$. Set $W:=\Delta\Psi_q.$
By Lemma~\ref{lem:point-adjoint}, the functions $\Psi_q$ and $W$ are smooth up to the boundary near $\overline U_1$ and satisfy
\begin{equation}\label{eq:adjoint-system}
\Delta\Psi_q=W,\quad \Delta W=-q\Psi_q\quad\text{in }\mathcal D'(\Omega).
\end{equation}
Moreover,
\begin{equation}\label{eq:zero-adjoint-cauchy-data}
 \Psi_q=\partial_\nu\Psi_q=W=\partial_\nu W=0 \quad\text{on }U_1.
\end{equation}

We briefly verify the interior regularity required for unique continuation.  Since $\Psi_q\in L^{p'}(\Omega)$, $q\in L^\infty(\Omega)$, and $\Delta^2\Psi_q=-q\Psi_q$ in $\mathcal D'(\Omega),$ the interior $W^{4,r}$ estimates, iterated together with the Sobolev embedding theorem, give $\Psi_q\in H^4_{\mathrm{loc}}(\Omega),\, W\in H^2_{\mathrm{loc}}(\Omega).$
Lemma~\ref{lem:boundary-ucp-system}, applied to
\eqref{eq:adjoint-system}--\eqref{eq:zero-adjoint-cauchy-data}, first
gives $\Psi_q=W=0$ in an interior neighborhood of $U_1$.  Since
$\Omega$ is connected, the interior unique-continuation  then yields 
\begin{equation}\label{eq:adjoint-zero}
\Psi_q=W=0\quad\text{in }\Om.
\end{equation}

We next show that \eqref{eq:adjoint-zero} contradicts the definition of $\Psi_q$ given in \eqref{eq:adjoint-duality-real}. To achieve this, we choose inward normal coordinates near $x_0$, and let $y_0$ be the tangential coordinate corresponding to $x_0$.  Take
$\eta\in C_c^\infty(\mathcal Y)$ and
$\chi\in C_c^\infty([0,\varepsilon))$ such that $\eta(y_0)=1$ and $\chi=1$  near $t=0.$
Denote by
\begin{equation}\label{eq:test-Z}
Z(y,t)=\eta(y)\chi(t)\Big( t-\frac{a(y,0)}{2}t^2\Big)\in C^\infty(\overline\Omega)
\end{equation}
in this coordinate neighborhood and extend it by zero to
$\overline\Omega$. 
At $t=0$,
\begin{equation}
Z=0,\quad \partial_tZ=\eta, \quad \partial_t^2Z=-a(y,0)\eta.
\end{equation}
Since $Z(\,\cdot\,,0)$ vanishes identically, its tangential derivatives
also vanish.  Hence
\begin{equation}
\left.\Delta Z\right|_{t=0} =\left.\bigl( \partial_t^2Z+a(y,0)\partial_tZ+\Delta_{g_0}Z\bigr)\right|_{t=0}=0.
\end{equation}
Outside the support of $\eta$, the function $Z$ vanishes near the
boundary.  Therefore
\begin{equation}\label{eq:Z-navier}
Z=\Delta Z=0\quad\text{on }\partial\Omega.
\end{equation}
Since $\partial_\nu=-\partial_t$ at $t=0$, we have
\begin{equation}\label{eq:Z-normal-derivative}
\partial_\nu Z(x_0)= -\partial_tZ(y_0,0)=-1.
\end{equation}

Let $f:=L_qZ$.  Then $f\in L^p(\Omega)$, and \eqref{eq:Z-navier}, together with uniqueness for the homogeneous Navier problem, implies that $G_qf=Z.$ Using \eqref{eq:adjoint-duality}, \eqref{eq:adjoint-zero}, and
\eqref{eq:Z-normal-derivative}, we obtain
\begin{equation}
0=\int_\Omega\Psi_qf\,dx=\partial_\nu(G_qf)(x_0)=\partial_\nu Z(x_0)=-1,
\end{equation}
which is a contradiction. Thus, $A_q$ and $B_q$ cannot both vanish
identically on any nonempty open subset of $\Gamma$.

\emph{ Step 2.} It remains to prove the boundary expansion formula.  Suppose first that $A_q\not\equiv0$ on $\Gamma$.  By continuity, there exist a nonempty connected open set $U\Subset\Gamma$ and $c_0>0$ such that
\begin{equation}
|A_q(y)|\geq c_0,\quad y\in\overline U.
\end{equation}
We choose $U$ sufficiently small so that it is contained in a boundary
coordinate neighborhood.  Since $\Psi_q(y,0)=0$, Taylor's formula in
the normal variable gives
\begin{equation}
\Psi_q(y,t)=t\,\partial_t\Psi_q(y,0)+t^2R_1(y,t).
\end{equation}
where
\begin{equation}\label{eq:R1-expression}
R_1(y,t)=\int_0^1(1-s)\,\partial_t^2\Psi_q(y,st)\,ds.
\end{equation}
Thus, with $\psi_1(y):=\partial_t\Psi_q(y,0)=-\partial_\nu\Psi_q(y) = -A_q(y),$ the conclusions follow with $\kappa=1$.

Suppose now that $A_q\equiv0$ on $\Gamma$.  We have $B_q\not\equiv0$.  Hence, after choosing $U$ smaller if
necessary, there exist a nonempty connected open set
$U\Subset\Gamma$ and $c_0>0$ such that
\begin{equation}\label{eq:B-lower-bound}
|B_q(y)|\geq c_0,\quad y\in\overline U.
\end{equation}
On $U$ we have
\begin{equation}
\Psi_q(y,0)=0,\quad \partial_t\Psi_q(y,0)=-A_q(y)=0,\quad W(y,0)=\Delta\Psi_q(y,0)=0.
\end{equation}
Using the normal-coordinate expression for the Laplacian and the fact
that $\Psi_q(\,\cdot\,,0)$ vanishes identically, we obtain
\begin{equation}\label{eq:second-normal-vanishing}
 \partial_t^2\Psi_q(y,0)=0.
\end{equation}

Differentiating $W=\Delta\Psi_q$ with respect to $t$ gives
\begin{equation}
\partial_tW=\partial_t^3\Psi_q+(\partial_ta)\partial_t\Psi_q+a\partial_t^2\Psi_q+\partial_t\bigl(\Delta_{g_t}\Psi_q\bigr).
\end{equation}
At $t=0$, the second and third terms vanish by
\eqref{eq:second-normal-vanishing}.  The last term also vanishes because
\begin{equation}
\left.\partial_t\bigl(\Delta_{g_t}\Psi_q\bigr)\right|_{t=0}=\left.(\partial_t\Delta_{g_t})\Psi_q \right|_{t=0} +\Delta_{g_0}\left.\partial_t\Psi_q\right|_{t=0}=0.
\end{equation}
Indeed, both $\Psi_q(\,\cdot\,,0)$ and
$\partial_t\Psi_q(\,\cdot\,,0)$ vanish identically on $U$.  Therefore $\partial_tW(y,0)=\partial_t^3\Psi_q(y,0).$
Since $\partial_\nu=-\partial_t$, $\psi_3(y):=\partial_t^3\Psi_q(y,0)=-\partial_\nu W(y)=-B_q(y).$
It then follows from \eqref{eq:B-lower-bound} that
\begin{equation}
|\psi_3(y)|\geq c_0, \quad y\in\overline U.
\end{equation}
Finally, Taylor's formula and the vanishing of the normal derivatives of orders $0$, $1$, and $2$ give
\begin{equation}
\Psi_q(y,t)=\frac{\psi_3(y)}{3!}\,t^3+t^4R_3(y,t),
\end{equation}
where the remainder term $R_3$ is given by
\begin{equation}\label{eq:R3-expression}
R_3(y,t) =\frac{1}{3!}\int_0^1(1-s)^3\,\partial_t^4\Psi_q(y,st)\,ds.
\end{equation}
Thus the required expansion holds with $\kappa=3$, completing the proof.
\end{proof}

We next show the following lemma.
\begin{lemma}\label{lem:initial-adjoint-jets}
For $j=1,2$, let $q_j\in\mathcal Q^\infty_\Gamma(\Omega)$ satisfy Assumption {\bf (A)} and let $F_j$ be smooth such that $F_j(x,0)=\partial_uF_j(x,0)=0$ in $\overline\Omega$.
Assume that
\begin{equation}\label{eq:equal-first-linearizations}
D\Lambda_{q_1,F_1}(0)[H]=D\Lambda_{q_2,F_2}(0)[H],\quad H\in\mathcal X_{\mathbb R}.
\end{equation}
Then, in every inward normal coordinate neighborhood contained in $\Gamma$, the corresponding point adjoint function
$\Psi_j:=\Psi_{q_j}$ satisfy
\begin{equation}\label{eq:initial-adjoint-jets}
\partial_t^m\Psi_1(y,0)=\partial_t^m\Psi_2(y,0),\quad m=0,1,2,3.
\end{equation}
\end{lemma}

\begin{proof}
Fix an open set $V\Subset\Gamma$ contained in a boundary coordinate
neighborhood, and use the inward normal coordinates $(y,t)$.
For $j=1,2$, set
\begin{equation}
\psi_{j,m}(y):=\partial_t^m\Psi_j(y,0),\quad m\geq0.
\end{equation}
By Lemma~\ref{lem:point-adjoint}, $\Psi_j$ is smooth up to the boundary
near $\overline V$ and satisfies
\begin{equation}
\Psi_j=\Delta\Psi_j=0\quad\text{on }V.
\end{equation}
The first boundary condition gives $ \psi_{1,0}=\psi_{2,0}=0.$ Using $\partial_\nu=-\partial_t$ at $t=0$, the trace identity
\eqref{eq:trace-equality} gives $\psi_{j,1}=-\partial_\nu\Psi_j,$ and hence $\psi_{1,1}=\psi_{2,1}.$

We write $a_0(y):=a(y,0)$ and $a_1(y):=\partial_ta(y,0).$ Recalling that
\begin{equation}
\Delta=\partial_t^2+a(y,t)\partial_t+\Delta_{g_t}
\end{equation}
at $t=0$, using $\Delta\Psi_j(y,0)=0$, we can get 
\begin{equation}
0=\psi_{j,2}+a_0\psi_{j,1}+\Delta_{g_0}\psi_{j,0}.
\end{equation}
Since $\psi_{j,0}=0$ identically on $V$, it follows that $\psi_{j,2}=-a_0\psi_{j,1}.$
The coefficient $a_0$ depends only on the geometry of $\partial\Omega$. Thus, we have $\psi_{1,2}=\psi_{2,2}.$

Finally, differentiating $\Delta\Psi_j$ with respect to $t$ and then setting $t=0$ yields
\begin{equation}
\left.\partial_t(\Delta\Psi_j)\right|_{t=0}=\psi_{j,3}+a_1\psi_{j,1}+a_0\psi_{j,2}+\left(\partial_t\Delta_{g_t})\right|_{t=0}\psi_{j,0}+\Delta_{g_0}\psi_{j,1}.
\end{equation}
The fourth term vanishes because $\psi_{j,0}=0$ identically.  Hence
\begin{equation}\label{eq:psi3}
 \left.\partial_t(\Delta\Psi_j)\right|_{t=0}
 =
 \psi_{j,3}
 +a_1\psi_{j,1}
 +a_0\psi_{j,2}
 +\Delta_{g_0}\psi_{j,1}.
\end{equation}
Since $\partial_t=-\partial_\nu$ at $t=0$, the left-hand side is
\begin{equation}
 \left.\partial_t(\Delta\Psi_j)\right|_{t=0}
 =
 -\partial_\nu\Delta\Psi_j.
\end{equation}
It has the same value for $j=1,2$ by
\eqref{eq:trace-equality}.  All the remaining terms on the right-hand
side of \eqref{eq:psi3}, except $\psi_{j,3}$, are also the same for
$j=1,2$.  Therefore, $\psi_{1,3}=\psi_{2,3}.$
This proves \eqref{eq:initial-adjoint-jets}.
\end{proof}

\section{Construction of boundary-localized semiclassical solutions}\label{sec:poisson}

Throughout this section, we assume that $q\in\mathcal Q_\Gamma^\infty(\Omega)$ satisfies Assumption {\bf (A)}.  We will construct suitable solutions of $L_qv=0$, whose boundary data and leading semiclassical behavior are localized near $\Gamma$, by solving the associated transport equations up to a sufficiently high finite order.

The construction can be viewed as an explicit finite-order realization of the local Poisson operator for an elliptic boundary problem. We refer to Grubb~\cite{Grubb1996} for the general pseudodifferential boundary calculus and Poisson operators.  Boundary normal coordinates and complete symbol expansions of boundary operators also form the basis of the classical
boundary determination method of Lee and Uhlmann~\cite{LeeUhlmann1989}. Here we apply these ideas to the fourth-order operator
$L_q=\Delta^2+q$ with Navier boundary conditions. To recover a given normal derivative $\partial_t^r q|_{\partial\Omega}$, we construct suitable semiclassical boundary solutions only through the finite order at which this derivative first contributes.

We begin with making some preparations.

We choose relatively open sets $\Gamma_1\Subset\Gamma_0\Subset\Gamma$ such that the inward boundary normal coordinates introduced in Subsection~\ref{subsec:pro-adjoint-Psi} are defined on $\Gamma_0\times[0,t_0)$.  When the construction is applied on the nondegenerate set $U$ furnished by Lemma~\ref{lem:AB}, we further choose $\Gamma_0\Subset U$.  Thus, in that application, we can choose $\Gamma_1\Subset\Gamma_0\Subset U\Subset\Gamma.$

Recall that the Laplacian in the coordinates $(y,t)$ is given by
\begin{equation}\label{eq:lap-normal}
\Delta=\partial_t^2+\alpha(y,t)\partial_t+\Delta_{g_t}.
\end{equation}
We also write
\begin{equation}\label{eq:normal-volume-density}
dx=J(y,t)\,dy\,dt,\quad J_0(y):=J(y,0),\quad \alpha(y,t)=\partial_t\log J(y,t).
\end{equation}

Fix $0<t_1<t_0$. To use the rescaled variable $z=t/h$, extend
$g_t$, $\alpha$, and $q$ smoothly to $\Gamma_0\times[0,\infty)$, without changing them for
$0\le t\le t_1$. The extensions can be chosen independent of $t$ for large $t$, with $g_t$ uniformly positive definite.  For any compact set $\mathcal K\Subset\Gamma_0$, all derivatives of the extended coefficients
are bounded on $\mathcal K\times[0,\infty)$.  We note that, these extensions are auxiliary: every approximate solution is multiplied by a suitable normal cutoff function and therefore depends only on the original coefficients.

Let $a\in C_c^\infty(\Gamma_1)$ and let $\phi\in C^\infty(\Gamma_0;\R)$ satisfy
\begin{equation}\label{eq:nonstationary-phase}
 \mu(y):=|d\phi(y)|_{g_0}\ge c_0>0
\end{equation}
on a neighborhood of $\supp a$ for some constant $c_0>0$.  Extend $\phi$ constantly to a new function $\phi(y,t)$ in the normal direction, such that $\partial_t^k\phi(y,t)=0$ for any $k\ge 1$. We also extend $e^{i\phi/h}a$ by zero to
$\partial\Omega\setminus\Gamma_1$.  For $0<h\ll1$, let $v_h$ satisfy the following equation
\begin{equation}\label{eq:poisson-bvp}
 \begin{cases}
 L_qv_h=0&\text{in }\Om,\\
 v_h=e^{i\phi/h}a&\text{on }\partial\Om,\\
 \Delta v_h=0&\text{on }\partial\Om.
 \end{cases}
\end{equation}
Noting that the boundary data is complex-valued, accordingly, the linear problem \eqref{eq:poisson-bvp} is understood after
complexification.  
Moreover, the real Fr\'echet derivatives of the measurement map $\Lambda_{q,F}$ are used only through their unique complex multilinear extensions (see Section~\ref{sec:linearization}). We will construct exact solutions $v_h$ to \eqref{eq:poisson-bvp} in the form
\begin{equation}
v_h=V_h+R_h,\quad V_h= \chi(t)e^{i\phi(y)/h}e^{-\mu(y)t/h}b_h(y,t/h).
\end{equation}
Here $b_h$ denotes the semiclassical amplitude function, and $R_h$ is the remainder term.

\subsection{Normal rescaling and preliminary estimates}

Set $z=t/h$.   A direct calculation for $V=V(y,z)$ gives
\begin{align}
\mathcal P_hV &:={h^2e^{-i\phi/h}\Delta(e^{i\phi/h}V)}\big|_{t=hz}=\partial_z^2V-|d\phi|_{g_{hz}}^2V\notag\\
&\quad+h\bigl[\alpha(y,hz)\partial_zV+2i\langle d\phi,d_yV\rangle_{g_{hz}}+i(\Delta_{g_{hz}}\phi)V\bigr]
+h^2\Delta_{g_{hz}}V.
\label{eq:conjugated-laplacian-explicit}
\end{align}
It is clear that we can write $\mathcal P_h=P_0+hP_1+\cdots+h^NP_N+\cdots$, with $P_0=\partial_z^2-\mu^2(y)$.
Taking $V(y,z)=e^{-\mu(y)z}B(y,z)$ (i.e. $B(y,z)=e^{\mu(y)z}V(y,z)$), we have
\begin{equation}\label{eq:Tmu}
P_0(e^{-\mu z}B)=e^{-\mu z}T_\mu B,\quad T_\mu:=\partial_z^2-2\mu\partial_z.
\end{equation}

For given $\mathcal K\Subset\Gamma_0$, and $N,k\ge0$, we define the seminorm of a function $P$
\begin{equation}\label{eq:weighted-polynomial-seminorm}
|P|_{k,N;\mathcal K}:=\max_{|\gamma|+j\le k}\sup_{(y,z)\in\mathcal K\times[0,\infty)}\langle z\rangle^{-N}
|\partial_y^\gamma\partial_z^j P(y,z)|.
\end{equation}
Here $\langle z\rangle:=(1+z^2)^{1/2}.$
A family $P_h(y,z)$ belongs to $\mathscr S_{\rm pol}(\mathcal K)$ if, for every $k\ge0$, there are an integer $N_k$ and a constant $C_k$ such that
\begin{equation}\label{eq:Spol-definition}
 |P_h|_{k,N_k;\mathcal K}\le C_k,\quad0<h<h_0.
\end{equation}
For $\ell\in\mathbb R$, we write $F_h=O_{\mathscr S_{\rm pol}}(h^\ell)$, if the family $h^{-\ell}F_h$ is bounded in
$\mathscr S_{\rm pol}(\mathcal K)$ uniformly for $0<h<h_0$.

We will prove the following lemma that gives the expansions of $e^{\mu z}\mathcal P^k_he^{-\mu z}$ for $k=1,2$.

\begin{lemma}\label{lem:operator-expansion}
For given integer $N\ge1$, there are differential operators $Q_0,\dots,Q_N$ that are smooth in $y$ and polynomial in $z$, such that
\begin{equation}\label{eq:Q-expansion}
 e^{\mu z}\mathcal P_he^{-\mu z}
 =\sum_{j=0}^Nh^jQ_j+h^{N+1}\mathcal E_{N,h},
 \quad Q_0=T_\mu.
\end{equation}
For given compact set $\mathcal K\Subset\Gamma_0$ and every $k,N_0\in\mathbb N$, there exist $N_1\in\mathbb N$ and a constant
$C>0$, independent of $h$, such that the remainder $\mathcal E_{N,h}$ satisfies
\begin{equation}\label{eq:operator-remainder-seminorm}
\left|\mathcal E_{N,h}w\right|_{k,N_1;\mathcal K}\leq C\left|w\right|_{k+2,N_0;\mathcal K}, \quad 0<h<h_0,
\end{equation}
for any smooth function $w=w(y,z)$ satisfying $|w|_{k+2,N_0;\mathcal K}<\infty.$
Moreover, 
\begin{equation}\label{eq:squared-Q-expansion}
e^{\mu z}\mathcal P_h^2e^{-\mu z}=T_\mu^2+\sum_{m=1}^{N}h^m\mathcal A_m+h^{N+1}\widetilde{\mathcal E}_{N,h},
\end{equation}
where
\begin{equation}\label{eq:Ageom-definition}
\mathcal A_m=\sum_{\substack{\ell+r=m, \ell,r\geq0}}Q_\ell Q_r, \quad Q_0=T_\mu.
\end{equation}
The remainder term $\widetilde{\mathcal E}_{N,h}$ satisfies the analogous estimate
\begin{equation}\label{eq:squared-operator-remainder-seminorm}
|\widetilde{\mathcal E}_{N,h}w|_{k,N_2;\mathcal K}\leq C\left|w\right|_{k+4,N_0;\mathcal K},\quad 0<h<h_0,
\end{equation}
for some $N_2\in\mathbb N$.
\end{lemma}

\begin{proof}
Recalling $t=hz$, we apply Taylor's formula at $t=0$ to each smooth scalar coefficient appearing in \eqref{eq:conjugated-laplacian-explicit}. If $c=c(y,t)$ denotes one of these scalar coefficient functions, then, for every fixed integer
$N\geq1$, we have
\begin{equation}\label{eq:coefficient-taylor-expansion}
c(y,hz)=\sum_{\ell=0}^{N}\frac{h^\ell z^\ell}{\ell!}\partial_t^\ell c(y,0)+h^{N+1}z^{N+1}r_{N,c}(y,z;h),
\end{equation}
where
\begin{equation}\label{eq:coefficient-taylor-remainder}
r_{N,c}(y,z;h)=\frac{1}{N!}\int_0^1(1-\theta)^N\partial_t^{N+1}c(y,\theta hz)\,d\theta.
\end{equation}
By construction, for every extended coefficient $c=c(y,t)$, compact set $\mathcal K\Subset\Gamma_0$, multi-index $\gamma$, and every nonnegative integer $j$, there exists a constant
$C_{\mathcal K,\gamma,j}>0$ such that
\begin{equation}\label{eq:extended-coefficient-bounds}
\sup_{(y,t)\in\mathcal K\times[0,\infty)}|\partial_y^\gamma\partial_t^j c(y,t)| \leq C_{\mathcal K,\gamma,j}.
\end{equation}
Hence, we can get
\begin{equation}\label{eq:taylor-remainder-polynomial-bound}
\max_{|\gamma|+j\leq k}\left|\partial_y^\gamma\partial_z^j\bigl(z^{N+1}r_{N,c}(y,z;h)\bigr)\right|\leq C_k\langle z\rangle^{N+1},
\end{equation}
uniformly for $(y,z)\in\mathcal K\times[0,\infty)$ and $0<h<h_0$.

Conjugation by $e^{-\mu(y)z}$ gives
\begin{equation}\label{eq:conjugated-derivatives}
e^{\mu z}\partial_ze^{-\mu z}=\partial_z-\mu,\quad e^{\mu z}\partial_{y_j}e^{-\mu z}=\partial_{y_j}-z\partial_{y_j}\mu.
\end{equation}
It follows from \eqref{eq:coefficient-taylor-expansion}--\eqref{eq:conjugated-derivatives} that the coefficient of each power $h^\ell$ is a differential operator $Q_\ell$ whose coefficients are smooth in $y$ and polynomial in $z$. Particularly, the leading term is $Q_0=T_\mu.$ Moreover, \eqref{eq:taylor-remainder-polynomial-bound} and Leibniz' rule imply that, for any smooth function $w=w(y,z)$, there exist an integer $N_1\in\mathbb N$ and a constant $C>0$, independent of $h$ and $w$, such that
\begin{equation}
\left|\mathcal E_{N,h}w\right|_{k,N_1;\mathcal K}\leq C\left|w\right|_{k+2,N_0;\mathcal K}, \quad 0<h<h_0,
\end{equation}
Thus \eqref{eq:Q-expansion} and \eqref{eq:operator-remainder-seminorm} follows.

To obtain the expansion of $e^{\mu z}\mathcal P_h^2e^{-\mu z}$, we set
\begin{equation}
\mathcal Q_h:=e^{\mu z}\mathcal P_he^{-\mu z}=\sum_{\ell=0}^{N}h^\ell Q_\ell+h^{N+1}\mathcal E_{N,h}.
\end{equation}
Collecting all terms of order at most $h^N$ gives
\begin{equation}
\mathcal Q_h^2=T_\mu^2+\sum_{m=1}^{N}h^m\sum_{\substack{\ell+r=m,\ell,r\geq0}}Q_\ell Q_r+h^{N+1}\widetilde{\mathcal E}_{N,h}.
\end{equation}
All terms included in $\widetilde{\mathcal E}_{N,h}$ contain either a
Taylor remainder or a product $Q_\ell Q_r$ with $\ell+r\geq N+1$.
Since each $Q_\ell$ is a differential operator of order at most two, \eqref{eq:operator-remainder-seminorm} and Leibniz'
rule show that, for any compact set $\mathcal K\Subset\Gamma_0$ and every pair of nonnegative integers
$k$ and $N_0$, there exist a nonnegative integer $N_2$ and a constant
$C>0$ such that
\begin{equation}
|\widetilde{\mathcal E}_{N,h}w|_{k,N_2;\mathcal K}\leq C|w|_{k+4,N_0;\mathcal K}, \quad 0<h<h_0.
\end{equation}
The constant $C$ is independent of $h$ and $w$. This proves \eqref{eq:squared-Q-expansion}--\eqref{eq:Ageom-definition} as desired. 
\end{proof}

Let $\mathcal K\Subset\Gamma_0$ be a compact set on which $\mu\geq c_0>0$ (see \eqref{eq:nonstationary-phase}), and let $P=P(y,z;h)$ be smooth function satisfying
\begin{equation}
\operatorname{supp}_y P(\cdot,z;h)\subset\mathcal K,\quad z\geq0,\quad 0<h<h_0.
\end{equation}
For a fixed nonnegative integer $s$, suppose that there exist a nonnegative integer $N_s$ and constants $C_{\gamma,j}>0$, independent of $h$, such that
\begin{equation}\label{eq:polynomial-amplitude-bound}
\left|\partial_y^\gamma\partial_z^jP(y,z;h)\right|\le C_{\gamma,j}\langle z\rangle^{N_s}, \quad |\gamma|+j\leq s.
\end{equation}
Moreover, let $c_h\in C_c^\infty(\Gamma_0)$ be supported in a fixed
compact subset of $\Gamma_0$ and assume that
\begin{equation}\label{eq:uniform-boundary-amplitude}
\sup_{0<h<h_0}\|c_h\|_{C^k(\Gamma_0)}<\infty \quad\text{for every }k\geq0.
\end{equation}
We have the following useful estimates.
\begin{lemma}\label{lem:bl-scaling}
For every $m\in\mathbb R$, $1<r<\infty$, and $\chi\in C_c^\infty([0,t_0))$, it holds that
\begin{equation}\label{eq:bl-Wsr-scaling}
\big\|h^m\chi(t)e^{i\phi(y)/h}e^{-\mu(y)t/h}P\left(y,t/h;h\right)\big\|_{W^{s,r}(\Omega)}\leq C h^{m+1/r-s},
\quad 0<h<h_0.
\end{equation}
Here the function in \eqref{eq:bl-Wsr-scaling} is extended by zero
outside the chosen boundary coordinate neighborhood, and the constant
$C$ is independent of $h$.

After extending $c_h$ by zero to $\partial\Omega$, one has, for all $\sigma\geq0$,
\begin{equation}\label{eq:boundary-Wsr-scaling}
\big\|h^m e^{i\phi/h}c_h\big\|_{W^{\sigma,r}(\partial\Omega)}\leq C h^{m-\sigma}, \quad 0<h<h_0,
\end{equation}
where $C$ is also independent of $h$.
\end{lemma}

\begin{proof}
Recall that, in the chosen boundary normal coordinates, the volume element has the form $dx=J(y,t)\,dy\,dt,$
where $J$ is smooth and bounded above and below on
$\mathcal K\times[0,t_0]$.

We first consider the case $s=0$. By
\eqref{eq:polynomial-amplitude-bound}, the support condition in $y$, and the lower bound $\mu\geq c_0$, we have
\begin{align}
\big\|h^m\chi(t)e^{i\phi/h}e^{-\mu t/h}P\left(y,t/h;h\right)\big\|_{L^r(\Omega)}^r&\leq
Ch^{mr}\int_{\mathcal K}\int_0^{t_0}e^{-rc_0t/h}\left\langle t/h\right\rangle^{rN_0}
\,dt\,dy\notag\\
&\leq Ch^{mr+1}\int_0^\infty e^{-rc_0z}\langle z\rangle^{rN_0}\,dz \leq Ch^{mr+1},
\end{align}
Hence, it follows that
\begin{equation}\label{eq:bl-Lr-scaling}
\|h^m\chi(t)e^{i\phi/h}e^{-\mu t/h}P\left(y,t/h;h\right)\|_{L^r(\Omega)}\leq Ch^{m+1/r}.
\end{equation}
Note that
\begin{equation}\label{eq:differentiated-layer-bound}
\left|D^\beta_{t,y}\left[h^m\chi(t)e^{i\phi/h}e^{-\mu t/h}P\left(y,t/h;h\right)\right]\right|
\leq C_\beta h^{m-|\beta|}e^{-c_0t/h} \left\langle t/h\right\rangle^{M_s}
\end{equation}
for some nonnegative integer $M_s$. Thus, \eqref{eq:bl-Lr-scaling} and \eqref{eq:differentiated-layer-bound} gives
\begin{equation}
\left\|D^\beta\left[h^m\chi(t)e^{i\phi/h}e^{-\mu t/h}P\left(y,t/h;h\right)\right]\right\|_{L^r(\Omega)}
\leq Ch^{m+1/r-|\beta|}.
\end{equation}
Since $0<h\leq1$ and $|\beta|\leq s$, the estimate \eqref{eq:bl-Wsr-scaling} follows.

It remains to prove the boundary estimate \eqref{eq:boundary-Wsr-scaling}. For every nonnegative integer $q$,
\eqref{eq:uniform-boundary-amplitude} gives
\begin{equation}\label{eq:integer-boundary-scaling}
\|h^me^{i\phi/h}c_h\|_{W^{q,r}(\partial\Omega)}\leq Ch^{m-q}.
\end{equation}
For a noninteger $\sigma\geq0$, choosing
\begin{equation}
\sigma=q+\theta,\quad q=\lfloor\sigma\rfloor, \quad 0<\theta<1,
\end{equation}
where $\lfloor\sigma\rfloor$ denotes the greatest integer less than or equal to $\sigma$.
The interpolation inequality between the estimates in $W^{q,r}(\partial\Omega)$ and $W^{q+1,r}(\partial\Omega)$ yields
\begin{align}
\|h^me^{i\phi/h}c_h\|_{W^{\sigma,r}(\partial\Omega)}\leq C\left(h^{m-q}\right)^{1-\theta}\left(h^{m-q-1}\right)^\theta= Ch^{m-\sigma}.
\end{align}
Hence, the proof is finished.
\end{proof}

\subsection{Construction of amplitude functions}
In this subsection, we first construct suitable amplitude functions $b_h$, and then construct the approximate solution $V_h$. We begin with observing that the (principal) homogeneous equation $T_\mu^2w=0$ has the general solution
\begin{equation}
w(z)=d_0+d_1z+d_2e^{2\mu z}+d_3ze^{2\mu z}.
\end{equation}
Thus, if
$w=d_0+d_1z$, then
\begin{equation}
\bigl(w(0),T_\mu w(0)\bigr)=(d_0,-2\mu d_1).
\end{equation}
Since $\mu\geq c_0>0$, the map $(d_0,d_1)\mapsto(d_0,-2\mu d_1)$ is uniquely determined by $(w(0),T_\mu w(0))$, and hence by the principal boundary traces of $v_h$.

In the following lemma, we solve the inhomogeneous equation $T_\mu^2w=f$.
\begin{lemma}\label{lem:polynomial-normal}
Let $\mu\geq c_0>0$. For any polynomial $f\in\mathbb C[z]$ and  $\beta_0,\beta_1\in\mathbb C$, there exists a unique polynomial solution $b$ such that
\begin{equation}
T_\mu^2b=f,\quad b(0)=\beta_0,\quad T_\mu b(0)=\beta_1.
\end{equation}
\end{lemma}

\begin{proof}
We first solve the following equation
\begin{equation}\label{eq:Smu}
T_\mu u=f,\quad u(0)=a.
\end{equation}
Since $T_\mu u=u''-2\mu u'$, the unique solution of \eqref{eq:Smu} is
\begin{equation}
\mathcal S_\mu[f;a](z):=a-\int_0^z e^{2\mu\tau}\left(\int_\tau^\infty e^{-2\mu s}f(s)\,ds\right)d\tau.
\end{equation}
It is clear that the solution is a polynomial whenever $f$ is a polynomial.

Now we set $c=\mathcal S_\mu[f;\beta_1],$ $b=\mathcal S_\mu[c;\beta_0].$ Then
\begin{equation}
T_\mu^2b=f,\quad b(0)=\beta_0,\quad T_\mu b(0)=\beta_1.
\end{equation}
For uniqueness, let $b$ and $\tilde b$ be two solutions with the same boundary data, and set $d=b-\tilde b$.
Then $T_\mu^2d=0$ with $d(0)=T_\mu d(0)=0$. Since every bounded polynomial solution of $T_\mu^2d=0$ is of the form $d(z)=a_0+a_1z.$ Thus, we have $a_0=d(0)=0$ and $-2\mu a_1=T_\mu d(0)=0.$ This yields the uniqueness.
Finally, the explicit formula for $\mathcal S_\mu$ shows that, for $\mu\geq c_0>0$, the coefficients of $b$ depend smoothly on $\mu$, $\beta_0$, $\beta_1$, and the coefficients of $f$.
\end{proof}

Choose $\chi\in C_c^\infty([0,t_0))$ such that $\chi=1$ near $t=0$. For some positive integer $M$, in the boundary normal coordinate neighborhood, we construct the approximate solution
\begin{equation}\label{eq:WKB}
V_h^{(M)}=\chi(t)e^{i\phi(y)/h}e^{-\mu(y)t/h}\sum_{m=0}^M h^m b_m\!\left(y,t/h\right).
\end{equation}
The amplitudes are chosen so that
\begin{equation}
\operatorname{supp}_y b_m(\cdot,z) \subset \operatorname{supp}a,\quad z\geq0,\quad 0\leq m\leq M.
\end{equation}
Since the amplitudes $b_m$ are compactly supported in the $y$ variable and $\chi$ is compactly supported in the $t$ variable, $V_h^{(M)}$ can be extended smoothly by zero to $\Omega$. 

For every integer $N\ge0$, Taylor's formula gives
\begin{equation}\label{eq:qTaylor}
h^4q(y,hz) =\sum_{r=0}^{N}h^{r+4}q_r(y)\frac{z^r}{r!}+h^{N+5}z^{N+1}\mathcal R_N(y,hz),\quad q_r=\partial_t^rq(y,0).
\end{equation}
Set $b_h^{(M)}=\sum_{m=0}^{M}h^m b_m$. Using
\begin{equation}
h^4e^{-i\phi/h}L_qe^{i\phi/h}
=\mathcal P_h^2+h^4q(y,hz),
\end{equation}
and the expansion formula in Lemma~\ref{lem:operator-expansion}, considering the powers of $h$, we can get the following equations
\begin{equation}\label{eq:transport}
T_\mu^2b_m=-\sum_{\ell=1}^{m}\mathcal A_\ell b_{m-\ell}-\mathbf 1_{\{m\geq4\}}
\sum_{\substack{r+k=m-4, r,k\geq0}}\frac{q_r(y)z^r}{r!}\,b_k,\quad 0\leq m\leq M.
\end{equation}
Empty sums are understood to be zero.
The operator $\mathcal A_\ell$ is defined in \eqref{eq:Ageom-definition}.  It depends on the geometry, the phase $\phi$, and their derivatives, and is independent of $q$.  It preserves polynomial dependence on $z$.
For instance, let $m=5$, we have
\begin{align*}
 T_\mu^2b_0&=0,\quad T_\mu^2b_1=-\mathcal A_1b_0,\quad  T_\mu^2b_2=-\mathcal A_1b_1-\mathcal A_2b_0,\quad T_\mu^2b_3=-\mathcal A_1b_2-\mathcal A_2b_1-\mathcal A_3b_0,\\
T_\mu^2b_4&=-\sum_{\ell=1}^4\mathcal A_\ell b_{4-\ell}-q_0b_0,\quad T_\mu^2b_5=-\sum_{\ell=1}^5\mathcal A_\ell b_{5-\ell}-q_0b_1-q_1zb_0.
\end{align*}
Clearly, the potential $q$ first enters the equation $T_\mu^2b_4$ for $b_4$.

The Dirichlet condition is imposed exactly by setting
\begin{equation}\label{eq:transport-bc1}
b_0(y,0)=a(y),\quad b_m(y,0)=0\quad(1\leq m\leq M).
\end{equation}
Since $\chi=1$ near $t=0$, \eqref{eq:Q-expansion} gives
\begin{equation}
h^2e^{-i\phi/h}\Delta V_h^{(M)}\big|_{t=0}=\sum_{m=0}^{M}h^m \big[T_\mu b_m+\sum_{j=1}^{m}Q_jb_{m-j}\big]_{z=0} +O(h^{M+1}).
\end{equation}
We therefore impose, for $0\leq m\leq M$, 
\begin{equation}\label{eq:transport-bc2}
T_\mu b_m(y,0)=-\sum_{j=1}^{m}(Q_jb_{m-j})(y,0).
\end{equation}
Together with \eqref{eq:transport} and \eqref{eq:transport-bc1}, this determines $b_m$ successively by Lemma~\ref{lem:polynomial-normal}. In particular, $b_0(y,z)=a(y)$.For $m\geq4$, $b_m$ depends on $q$ only through $q_0,\ldots,q_{m-4}$; the amplitudes $b_0,\ldots,b_3$ are independent of $q$.

We next study the contribution of each boundary normal derivative of $q$ to the amplitudes. 

Fix $r\geq0$ and recall $q_r(y)=\partial_t^r q(y,0)$. By \eqref{eq:transport}, it is clear that $q_r$ first enters the transport equations at order $m=r+4$, through the term
\begin{equation}\label{eq:first-qr-forcing}
-q_r(y)\frac{z^r}{r!}b_0(y,z)
=-q_r(y)a(y)\frac{z^r}{r!}.
\end{equation}
All other terms at the order $m$ depend only on $q_0,\ldots,q_{r-1}$. Moreover, the boundary conditions for $b_{r+4}$ involve only lower-order amplitudes. Thus the part of $b_{r+4}$ arising from $q_r$, denoted by $b_{r+4}^{q_r}$, is the unique polynomial solution of the following problem
\begin{equation}\label{eq:qrcorrODE}
\begin{cases}
T_\mu^2b_{r+4}^{q_r}=-q_r(y)a(y)\dfrac{z^r}{r!},\\
b_{r+4}^{q_r}(y,0)=T_\mu b_{r+4}^{q_r}(y,0)=0.
\end{cases}
\end{equation}

Let $s=\mu(y)z$. We have $T_\mu=\mu^2(y)(\partial_s^2-2\partial_s)$, and
\begin{equation}\label{eq:ar-def}
b_{r+4}^{q_r}(y,z)=q_r(y)a(y)\mu^{-r-4}(y)\mathfrak a_r(\mu(y)z),
\end{equation}
where $\mathfrak a_r$ is the polynomial solution of 
\begin{equation}\label{eq:ar-ODE}
\begin{cases}
(\partial_s^2-2\partial_s)^2\mathfrak a_r=-\dfrac{s^r}{r!},\\
\mathfrak a_r(0)=(\partial_s^2-2\partial_s)\mathfrak a_r(0)=0.
\end{cases}
\end{equation}

The following lemma gives an explicit formula for $\mathfrak a_r$ and determines its sign for $s>0$.
\begin{lemma}\label{lem:ar}
For every integer $r\geq0$, the unique polynomial solution of \eqref{eq:ar-ODE} is
\begin{equation}\label{eq:ar-formula}
\mathfrak a_r(s)=-\sum_{d=1}^{r+2}\gamma_{r,d}\frac{s^d}{d!},
\end{equation}
where
\begin{equation}\label{eq:ar-coefficients}
\gamma_{r,1}=\frac{r+1}{2^{r+3}},\quad \gamma_{r,d}=\frac{r-d+3}{2^{r-d+4}} \quad(2\leq d\leq r+2).
\end{equation}
In particular, $\mathfrak a_r(s)<0$ for $s>0$.
\end{lemma}

\begin{proof}
Put $e_d=s^d/d!$, with $e_d=0$ for $d<0$, and let
$\mathcal T=\partial_s^2-2\partial_s$. Then
\begin{equation}
\mathcal T^2e_d=e_{d-4}-4e_{d-3}+4e_{d-2}.
\end{equation}
Thus, \eqref{eq:ar-formula} satisfies the differential
equation in \eqref{eq:ar-ODE} precisely when
\begin{equation}\label{eq:gamma-local-recurrence}
\gamma_{r,k+4}-4\gamma_{r,k+3}+4\gamma_{r,k+2}=\delta_{kr},\quad 0\leq k\leq r,
\end{equation}
with $\gamma_{r,d}=0$ for $d\notin\{1,\ldots,r+2\}$.

For $0\leq k<r$, \eqref{eq:ar-coefficients} gives
\begin{equation}
\gamma_{r,k+4}-4\gamma_{r,k+3}+4\gamma_{r,k+2}=\frac{(r-k-1)-2(r-k)+(r-k+1)}{2^{r-k}}=0.
\end{equation}
For $k=r$, the same expression equals $4\gamma_{r,r+2}=1$. Moreover,
\begin{equation}
\mathfrak a_r(0)=0,\quad \mathcal T\mathfrak a_r(0) =-\gamma_{r,2}+2\gamma_{r,1}=0.
\end{equation}
Thus \eqref{eq:ar-formula} satisfies the equation \eqref{eq:ar-ODE}. Uniqueness follows from Lemma~\ref{lem:polynomial-normal}
with $\mu=1$. Since $\gamma_{r,d}>0$ for $1\leq d\leq r+2$, it follows that $\mathfrak a_r(s)<0$ for $s>0$.
\end{proof}

\subsection{Construction of the remainder term}
With these preparations, we construct the remainder term for $v_h$ and establish some of its properties.
\begin{lemma}\label{lem:poisson}
Let $L\geq0$ and choose an integer $M\geq L+3-1/p$ with $p$ satisfying $\max\{1,n/3\}<p<\infty$. For fixed $q$, the phase function $\phi$, and  the boundary data $a$,  the solution of
\eqref{eq:poisson-bvp} has the following form
\begin{equation}\label{eq:poisson-exact}
v_h=V_h^{(M)}+R_{M,h},\quad 0<h<h_0.
\end{equation}
where $V_h^{(M)}$ is given by \eqref{eq:WKB}, and uniformly for $0<h<h_0$, the remainder term $R_{M,h}$ satisfies
\begin{equation}\label{eq:remainder-W4p}
\|R_{M,h}\|_{W^{4,p}(\Omega)}
+\|R_{M,h}\|_{C^1(\overline\Omega)}
\leq C h^{M-3+1/p}\leq C h^L.
\end{equation}
The positive constant $C$ is uniform over the fixed admissible triples $(q,\phi,a)$.
\end{lemma}

\begin{proof}
Let us first apply Lemma~\ref{lem:operator-expansion} through order $M$. If $M\geq4$, use \eqref{eq:qTaylor} with $N=M-4$;
if $M<4$, the term $h^4q(y,hz)$ is already of order $h^{M+1}$. We then solve equations \eqref{eq:transport}--\eqref{eq:transport-bc2} to obtain $b_m$ successively for $m=0,\ldots,M$. By \eqref{eq:transport}, the coefficients of
$h^0,\ldots,h^M$ vanish. Since $q$ and its extension agree on $\supp\chi$, we obtain
\begin{equation}\label{eq:residual-structure}
h^4L_qV_h^{(M)}=h^{M+1}\chi(t)e^{i\phi/h}e^{-\mu t/h}r_{M,h}\!\left(y,t/h\right)+O_{L^p}(h^\infty),
\end{equation}
where $r_{M,h}$ has uniformly bounded polynomially weighted seminorms. Indeed, the omitted term is generated by
$[\Delta^2,\chi]$ and is supported where $t\geq c>0$. Moreover, the factor $e^{-\mu t/h}$ makes the commutator term $O(h^N)$ in $L^p$ for every $N$. Applying the estimate \eqref{eq:bl-Wsr-scaling} with $m=M-3$, $s=0$, and $r=p$ yields
\begin{equation}\label{eq:residual-Lp}
\|L_qV_h^{(M)}\|_{L^p(\Omega)}
\leq C h^{M-3+1/p}.
\end{equation}

By \eqref{eq:transport-bc1}, we have $V_h^{(M)}|_{\partial\Omega}=e^{i\phi/h}a$. Using \eqref{eq:transport-bc2}, the coefficients of $h^0,\ldots,h^M$ in 
$$h^2e^{-i\phi/h}\Delta V_h^{(M)}|_{\partial\Omega}$$ vanish. This gives
\begin{equation}
h^2\Delta V_h^{(M)}|_{\partial\Omega}=h^{M+1}e^{i\phi/h}c_{M,h},
\end{equation}
where $c_{M,h}$ has fixed compact support and uniformly bounded derivatives. Applying the estimate \eqref{eq:boundary-Wsr-scaling} with
$m=M-1$, $r=p$, and $\sigma=2-1/p$ yields
\begin{equation}\label{eq:boundary-defect}
\|\Delta V_h^{(M)}|_{\partial\Omega}\|_{W^{2-1/p,p}(\partial\Omega)}\leq C h^{M-3+1/p}.
\end{equation}

Let $\mathscr A_q$ be the Navier isomorphism from
Lemma~\ref{lem:full-navier-isomorphism}, and set
\begin{equation}\label{eq:R-problem}
R_{M,h}=\mathscr A_q^{-1}\bigl(-L_qV_h^{(M)},\,0,\,-\Delta V_h^{(M)}|_{\partial\Omega}\bigr).
\end{equation}
Then $v_h=V_h^{(M)}+R_{M,h}$ solves \eqref{eq:poisson-bvp}. By \eqref{eq:full-navier-estimate}, \eqref{eq:residual-Lp}, and \eqref{eq:boundary-defect}, we can get
\begin{equation}\label{eq:R-W4p}
\|R_{M,h}\|_{W^{4,p}(\Omega)}\leq C h^{M-3+1/p}.
\end{equation}
Since $4-n/p>1$, using the Sobolev embedding theorem, the same bound holds in $C^1(\overline\Omega)$. Using the condition $M-3+1/p\ge L$, the estimate \eqref{eq:remainder-W4p} follows. For fixed $M$, the constant $C>0$ can be chosen uniformly over the triple $(q,\phi,a)$.
\end{proof}

In the remainder of this subsection, we present some auxiliary results that will be used to establish the boundary determination of coefficients in Theorem~\ref{thm:boundary-jets}.

Set $\mathfrak p_0(s)=1$, $\mathfrak p_1(s)=s,$ and define the boundary data
\begin{equation}\label{eq:two-mode-boundary-data}
H_h^{(0)}=\bigl(e^{i\phi/h}a,0\bigr),\quad H_h^{(1)}=\bigl(0,-2\mu^2h^{-2}e^{i\phi/h}a\bigr).
\end{equation}

We have the following lemma:
\begin{lemma}\label{lem:two-navier-modes}
For each $\sigma\in\{0,1\}$, there is a solution $v_h^{(\sigma)}$ of $L_qv_h^{(\sigma)}=0$ with boundary data
$H_h^{(\sigma)}$. It admits the expansion and remainder estimate of
Lemma~\ref{lem:poisson}, with leading amplitude
\begin{equation}\label{eq:two-mode-leading-amplitude}
b_0^{(\sigma)}(y,z)
=
a(y)\mathfrak p_\sigma\bigl(\mu(y)z\bigr).
\end{equation}
For every $r\in\mathbb N_0$, the part of $b_{r+4}^{(\sigma)}$ that is linear in $q_r$, denoted by
$b_{r+4}^{(\sigma;q_r)}$, is
\begin{equation}\label{eq:two-mode-qr-correction}
b_{r+4}^{(\sigma;q_r)}(y,z)=q_r(y)a(y)\mu^{-r-4}(y)\mathfrak a_{r,\sigma}\bigl(\mu(y)z\bigr),
\end{equation}
where
\begin{equation}\label{eq:mode-polynomials}
\mathfrak a_{r,0}=\mathfrak a_r,\quad \mathfrak a_{r,1}=(r+1)\mathfrak a_{r+1}.
\end{equation}
\end{lemma}

\begin{proof}
The case $\sigma=0$ is Lemma~\ref{lem:poisson}. For $\sigma=1$, the
leading amplitude $b_0^{(1)}(y,z)=a(y)\mu(y)z$ belongs to the admissible kernel of $T_\mu^2$ and satisfies
\begin{equation}
b_0^{(1)}(y,0)=0,\quad T_\mu b_0^{(1)}(y,0)=-2\mu^2(y)a(y).
\end{equation}
These are the principal boundary conditions corresponding to $H_h^{(1)}$. For \(m\geq1\), the amplitudes \(b_m^{(1)}\) are determined recursively by \eqref{eq:transport}, with
\begin{equation}
b_m^{(1)}(y,0)=0,\quad T_\mu b_m^{(1)}(y,0)=-\sum_{\ell=1}^m\bigl(Q_\ell b_{m-\ell}^{(1)}\bigr)(y,0).
\end{equation}
Thus, Lemma~\ref{lem:polynomial-normal} yields a unique admissible amplitude at each order. Moreover, the remainder terms and their corresponding estimates can be derived by the same arguments as in Lemma~\ref{lem:poisson}.

At order \(r+4\), the term $b_{r+4}^{(\sigma;q_r)}$ satisfies
\begin{equation}\label{eq:mode-qr-normal-problem}
T_\mu^2 b_{r+4}^{(\sigma;q_r)}=-q_r(y)a(y)\frac{z^r}{r!}\mathfrak p_\sigma\bigl(\mu(y)z\bigr),
\end{equation}
with the homogeneous conditions
\begin{equation}
b_{r+4}^{(\sigma;q_r)}(y,0)=T_\mu b_{r+4}^{(\sigma;q_r)}(y,0)=0.
\end{equation}
Set $s=\mu(y)z$. Since
\begin{equation}
T_\mu=\mu(y)^2(\partial_s^2-2\partial_s),\quad z^r=\mu^{-r}(y)s^r,
\end{equation}
the expression
\begin{equation}
b_{r+4}^{(\sigma;q_r)}(y,z)=q_r(y)a(y)\mu^{-r-4}(y)\mathfrak a_{r,\sigma}\bigl(\mu(y)z\bigr)
\end{equation}
reduces \eqref{eq:mode-qr-normal-problem} to
\begin{equation}\label{eq:mode-correction-equation}
(\partial_s^2-2\partial_s)^2\mathfrak a_{r,\sigma}=-\frac{s^r}{r!}\mathfrak p_\sigma(s),\quad
\mathfrak a_{r,\sigma}(0)=(\partial_s^2-2\partial_s)\mathfrak a_{r,\sigma}(0)=0.
\end{equation}
For $\sigma=0$, this is \eqref{eq:ar-ODE}, so
$\mathfrak a_{r,0}=\mathfrak a_r$. For $\sigma=1$, we have
\begin{equation}
-\frac{s^r}{r!}\mathfrak p_1(s)=-(r+1)\frac{s^{r+1}}{(r+1)!}.
\end{equation}
Uniqueness in Lemma~\ref{lem:polynomial-normal} therefore gives $\mathfrak a_{r,1}=(r+1)\mathfrak a_{r+1}$. This proves
\eqref{eq:two-mode-qr-correction} and \eqref{eq:mode-polynomials}.
\end{proof}


\begin{remark}\label{rem:localization}
Recalling that $\supp a\Subset\Gamma_1\Subset\Gamma$, the approximate solution $V_h^{(M)}$ is supported away from $x_0$. Moreover, for given $N>0$, the approximation order $M$ can be chosen so that
\begin{equation}
\|R_{M,h}\|_{L^p(\Omega)}+\|R_{M,h}\|_{C^0(\overline\Omega)}=O(h^N).
\end{equation}
Hence, for any fixed number of solutions, by H\"older's inequality and $\Psi_q\in L^{p'}(\Omega)$, we have
\begin{equation}
\Big|\int_\Omega\Psi_q\Big(\prod_{j=1}^m v_{j,h}-\prod_{j=1}^m V_{j,h}^{(M)}\Big)dx\Big|=O(h^N),
\end{equation}
Hence, the asymptotic expansion is determined entirely by the approximate solutions in the boundary coordinate neighborhood of $\Gamma_1$, where $\Psi_q$ is smooth. 
\end{remark}

Let $m\geq2$, choose boundary functions $a_1,\ldots,a_m\in C_c^\infty(\Gamma_1),$ and $\phi_1,\ldots,\phi_m\in C^\infty(\Gamma_0;\mathbb R)$. We assume that
\begin{equation}\label{eq:phase-cancellation}
\sum_{\ell=1}^m\phi_\ell=0, \quad \mu_\ell:=|d\phi_\ell|_{g_0}>0
\end{equation}
in a neighborhood of $\bigcup_{\ell=1}^m\operatorname{supp}a_\ell$. Set $S:=\sum_{\ell=1}^m\mu_\ell,$ and let $v_{\ell,h}$ be the solutions to \eqref{eq:poisson-bvp} constructed in Lemma~\ref{lem:poisson} with $(\phi,a)$ replaced by $(\phi_\ell,a_\ell)$.

Suppose that $\Psi$ is smooth in the corresponding boundary coordinate neighborhood and that, for some $\kappa\geq0$, $\partial_t^j\Psi(y,0)=0$ for $0\leq j<\kappa.$
We write $\psi_j(y):=\partial_t^j\Psi(y,0)$ and recall that $dx=J(y,t)\,dy\,dt$, $J_0(y):=J(y,0).$

The following proposition gives a finite-order asymptotic expansion for products of the semiclassical solutions constructed in Lemma~\ref{lem:poisson} and specifies the coefficient of each power of $h$.
\begin{proposition}\label{prop:coefficient-extraction}
For every integer $N\geq0$, it holds that
\begin{equation}\label{eq:full-product-expansion}
\int_\Omega\Psi\prod_{\ell=1}^m v_{\ell,h}\,dx=\sum_{\nu=0}^{N}h^{\kappa+1+\nu}\mathcal I_\nu+O(h^{\kappa+N+2}).
\end{equation}
The coefficient $\mathcal I_\nu$ is obtained by collecting all terms
whose indices satisfy
\begin{equation}\label{eq:coefficient-collection-rule}
\rho+\eta+m_1+\cdots+m_m=\nu,\quad \rho,\eta,m_1,\ldots,m_m\in\mathbb N_0.
\end{equation}
Here $\rho$ and $\eta$ index the normal Taylor expansions of $\Psi$
and the volume density $J$, respectively, whereas $m_\ell$ indexes the
semiclassical amplitude $b_{\ell,m_\ell}$ of the $\ell$th solution. In particular,
\begin{equation}\label{eq:I0-leading}
\mathcal I_0=\int_{\Gamma_1}\frac{\psi_\kappa}{S^{\kappa+1}}\big(\prod_{\ell=1}^m a_\ell\big)J_0\,dy.
\end{equation}
\end{proposition}

\begin{proof}
We choose $L>\kappa+N+2$, and take the order $M$ sufficiently large so that Lemma~\ref{lem:poisson} gives
\begin{equation}\label{eq:product-correction-estimate}
\|R_{\ell,M,h}\|_{W^{4,p}(\Omega)}+\|R_{\ell,M,h}\|_{C^1(\overline\Omega)}\leq C h^L, \quad 1\leq \ell\leq m.
\end{equation}
Recall that, in boundary normal coordinates, $t=hz$, the approximate solutions then take the form
\begin{equation}
V_{\ell,h}^{(M)}(y,hz)=\chi(hz)e^{i\phi_\ell(y)/h}e^{-\mu_\ell(y)z}\sum_{j=0}^{M}h^j b_{\ell,j}(y,z).
\end{equation}
Using $\sum_{\ell=1}^m\phi_\ell=0$, we have
\begin{equation}
\prod_{\ell=1}^m e^{-\mu_\ell(y)z}=e^{-S(y)z},\quad
 S(y)=\sum_{\ell=1}^m\mu_\ell(y).
\end{equation}
The Taylor expansions of $\Psi$ and $J$ give
\begin{align}
\Psi(y,hz) &=\sum_{\rho=0}^{N} h^{\kappa+\rho}\frac{\psi_{\kappa+\rho}(y)}{(\kappa+\rho)!}z^{\kappa+\rho}
+O\!\left(h^{\kappa+N+1}\langle z\rangle^{\kappa+N+1}\right),\label{eq:Psi-rescaled-expansion}\\
J(y,hz)&=\sum_{\eta=0}^{N}h^\eta J_\eta(y)z^\eta+O\!\left(h^{N+1}\langle z\rangle^{N+1}\right),\quad J_\eta(y)=\frac{\partial_t^\eta J(y,0)}{\eta!}.
\label{eq:J-rescaled-expansion}
\end{align}
These expansions are uniform on the compact set containing the
supports of the amplitudes.

Since $dt=h\,dz$, a term involving
$\psi_{\kappa+\rho}$, $J_\eta$, and
$b_{\ell,m_\ell}$, $1\leq\ell\leq m$, has order
\begin{equation}
 h^{\kappa+1+\rho+\eta+m_1+\cdots+m_m}.
\end{equation}
Consequently, the coefficient of $h^{\kappa+1+\nu}$ is obtained
precisely from the indices satisfying
\begin{equation}
 \rho+\eta+m_1+\cdots+m_m=\nu,
\end{equation}
which proves \eqref{eq:coefficient-collection-rule}. Because $S$ is bounded below by a positive constant on the relevant
compact set, all the resulting normal integrals converge uniformly:
\begin{equation}
 \int_0^\infty e^{-S(y)z}\langle z\rangle^K\,dz\leq C_K.
\end{equation}
The Taylor remainders in
\eqref{eq:Psi-rescaled-expansion}--\eqref{eq:J-rescaled-expansion}
therefore contribute $O(h^{\kappa+N+2})$.  Furthermore, since $\chi=1$ for sufficiently small $t$, replacing
$\chi^m(hz)$ by $1$ and extending the integral from $[0,t_0/h)$ to $[0,\infty)$. produces an $O(h^\infty)$ error.

For $\nu=0$, only $\rho=\eta=m_1=\cdots=m_m=0$ appears. Using $b_{\ell,0}(y,z)=a_\ell(y)$, we obtain
\begin{equation}
\mathcal I_0=\int_{\Gamma_1}\frac{\psi_\kappa(y)}{\kappa!}\Big(\prod_{\ell=1}^m a_\ell(y)\Big)J_0(y)
\Big(\int_0^\infty z^\kappa e^{-S(y)z}\,dz\Big)dy=\int_{\Gamma_1}\frac{\psi_\kappa(y)}{S(y)^{\kappa+1}}
\Big(\prod_{\ell=1}^m a_\ell(y)\Big)J_0(y)\,dy.
\end{equation}
Thus \eqref{eq:I0-leading} follows.

Finally, write $v_{\ell,h}=V_{\ell,h}^{(M)}+R_{\ell,M,h},$ $1\leq\ell\leq m.$
The identity
\begin{equation}\label{eq:product-difference}
\prod_{\ell=1}^m v_{\ell,h}-\prod_{\ell=1}^m V_{\ell,h}^{(M)}=\sum_{j=1}^mR_{j,M,h}\Big(\prod_{\ell<j}v_{\ell,h}\Big)
\Big(\prod_{\ell>j}V_{\ell,h}^{(M)}\Big)
\end{equation}
shows that every term on the right-hand side contains a correction factor. The approximate solutions are uniformly bounded in $L^\infty(\Omega)$,
and \eqref{eq:product-correction-estimate} gives
\begin{equation}
\|R_{\ell,M,h}\|_{L^p(\Omega)}+\|R_{\ell,M,h}\|_{L^\infty(\Omega)}\leq Ch^L.
\end{equation}
Hence, the solutions $v_{\ell,h}$ are also uniformly bounded in $L^\infty(\Omega)$.  Hölder's inequality and
$\Psi\in L^{p'}(\Omega)$ therefore yield
\begin{align}
&\Big|\int_\Omega\Psi\Big(\prod_{\ell=1}^m v_{\ell,h}-\prod_{\ell=1}^m V_{\ell,h}^{(M)}\Big)dx
\Big| \notag\\
&\leq\|\Psi\|_{L^{p'}(\Omega)}\sum_{j=1}^m\|R_{j,M,h}\|_{L^p(\Omega)}\prod_{\ell<j}\|v_{\ell,h}\|_{L^\infty(\Omega)}
\prod_{\ell>j}\|V_{\ell,h}^{(M)}\|_{L^\infty(\Omega)}\leq Ch^L.
\end{align}
Since $L>\kappa+N+2$, the remainder terms do not contribute to the coefficients in \eqref{eq:full-product-expansion}.
This completes the proof.
\end{proof}

\begin{remark}
The proof of Proposition~\ref{prop:coefficient-extraction} also applies to the solutions constructed in Lemma~\ref{lem:two-navier-modes}.  For a choice of modes $\sigma_\ell\in\{0,1\}$, the leading amplitudes are $b_{\ell,0}(y,z)=a_\ell(y)\mathfrak p_{\sigma_\ell}\bigl(\mu_\ell(y)z\bigr).$
\end{remark}

Using a suitable multiplier, we can also derive the following asymptotic expansion.
\begin{corollary}\label{cor:weighted-extraction}
Under the conditions of Proposition~\ref{prop:coefficient-extraction}, let \(c\in L^\infty(\Omega)\) be smooth in the boundary coordinate neighborhood.  Suppose that, for some
\(r\in\mathbb N_0\),
\begin{equation}\label{eq:multiplier-vanishing-jets}
\partial_t^s c(y,0)=0\quad (0\leq s<r),\quad c_r(y):=\partial_t^r c(y,0).
\end{equation}
Then, it holds that
\begin{align}
\int_\Omega\Psi c\prod_{j=1}^m v_{j,h}\,dx=h^{r+\kappa+1}\binom{r+\kappa}{r}\int_{\Gamma_1}\frac{\psi_\kappa(y)c_r(y)}
{S(y)^{r+\kappa+1}}\Big(\prod_{j=1}^m a_j(y)\Big)J_0(y)\,dy+O\!\big(h^{r+\kappa+2}\big).
\label{eq:weighted-leading-coefficient}
\end{align}
The remainder has the same uniformity as in
Proposition~\ref{prop:coefficient-extraction}.
\end{corollary}

\begin{proof}
Using the assumptions on $c$ and the expansion of $\Psi$, we can get
\begin{equation}
c(y,t)\Psi(y,t)=\frac{c_r(y)\psi_\kappa(y)}{r!\,\kappa!}t^{r+\kappa}+O\!\left(t^{r+\kappa+1}\right).
\end{equation}
After the change of variables \(t=hz\), the leading term of
\(c(y,hz)\Psi(y,hz)\) is
\[
h^{r+\kappa}
\frac{c_r(y)\psi_\kappa(y)}{r!\,\kappa!}
z^{r+\kappa}.
\]
Since \(dt=h\,dz\), its contribution is of order \(h^{r+\kappa+1}\).  The corresponding integral is
\begin{equation}
\frac{1}{r!\,\kappa!}
\int_0^\infty z^{r+\kappa}e^{-S(y)z}\,dz=\binom{r+\kappa}{r}[S(y)]^{-r-\kappa-1}.
\end{equation}
The terms involving an exact correction are \(O(h^L)\), where \(L\) is chosen
larger than \(r+\kappa+2\).  Hence the total remainder is
\(O(h^{r+\kappa+2})\).
\end{proof}

\section{Proof of the main theorems}\label{sec:mainproof}
In this section, we are concerned with proving the main theorems of this paper.

\subsection{Boundary determination of the potential and the leading coefficient}\label{sec:boundarydetermination}
We next prove Theorem~\ref{thm:boundary-jets} that gives the unique determination of the potential and the leading nonlinear spatial coefficient simultaneously.  

\begin{proof}[Proof of Theorem~\ref{thm:boundary-jets}]
For $j=1,2$, let $\Psi_j=\Psi_{q_j}$ be the adjoint function associated with $q_j$, and write $a_{j,s}(y):=\partial_t^s a_j(y,0),q_{j,s}(y):=\partial_t^s q_j(y,0).$ We divide the proof into six steps.

\emph{Step 1.}
Since $a_{1,0}\not\equiv0$ on $\Gamma$, there is a nonempty relatively open set on which $a_{1,0}$ does not vanish. The identity
\begin{equation}
D\Lambda_{q_1,F_1}(0)=D\Lambda_{q_2,F_2}(0)
\end{equation}
and Lemma~\ref{lem:first-green} imply $\partial_\nu\Psi_1=\partial_\nu\Psi_2$ and $\partial_\nu\Delta\Psi_1=\partial_\nu\Delta\Psi_2$ on $\Gamma$. Hence, by Lemma~\ref{lem:initial-adjoint-jets},
\begin{equation}
\partial_t^s\Psi_1(y,0)=\partial_t^s\Psi_2(y,0), \quad 0\leq s\leq3.
\end{equation}
Lemma~\ref{lem:AB} then yields a connected set
$U\Subset\Gamma$, an index $\kappa\in\{1,3\}$, and a function
$\psi_\kappa$ that does not vanish on $\overline U$ such that
\begin{equation}\label{eq:coupled-common-kappa}
\Psi_j(y,t)=\frac{\psi_\kappa(y)}{\kappa!}t^\kappa+O(t^{\kappa+1}),\quad j=1,2.
\end{equation}
We also choose $U$ so that $a_{1,0}\neq0$ on $\overline U$.

We fix a boundary coordinate neighborhood containing $\overline U$, and choose $\rho\in C^\infty$ such that $d\rho\neq0$ there. Let $\phi_\ell=c_\ell\rho$ with $c_\ell\neq0$ and $\sum_{\ell=1}^{m_*}c_\ell=0$, and let
\begin{equation}
\mu_\ell=|c_\ell||d\rho|_{g_0},\quad S=\sum_{\ell=1}^{m_*}\mu_\ell.
\end{equation}
For $\alpha_\ell\in C_c^\infty(U)$, let $v_{j,\ell,h}$ be the semiclassical solution for $L_{q_j}$ with boundary data $H_{\ell,h}^{(0)}=\bigl(e^{i\phi_\ell/h}\alpha_\ell,0\bigr).$
Applying the complex multilinear extension of
\eqref{eq:first-nonlinear-tensor} to $H_{1,h}^{(0)},\ldots,H_{m_*,h}^{(0)}$, we obtain
\begin{equation}\label{eq:coupled-leading-tensor}
\mathcal T_{j,h}:=D^{m_*}\Lambda_{q_j,F_j}(0)\bigl[H_{1,h}^{(0)},\ldots,H_{m_*,h}^{(0)}\bigr]=\int_\Omega\Psi_j a_j
\prod_{\ell=1}^{m_*}v_{j,\ell,h}\,dx.
\end{equation}
The condition \eqref{cond:DN-m-star} gives $\mathcal T_{1,h}=\mathcal T_{2,h}$. 

\medskip
\emph{Step 2.} We next recover the first four normal derivatives of $a_j$. Observing that the leading coefficient in \eqref{eq:coupled-leading-tensor} is
\begin{equation}\label{eq:unknown-profile-leading-asymptotic}
\mathcal T_{j,h}=h^{\kappa+1}\int_U\frac{\psi_\kappa a_{j,0}}{S^{\kappa+1}}\Big(\prod_{\ell=1}^{m_*}\alpha_\ell\Big)J_0\,dy
+O(h^{\kappa+2}).
\end{equation}
we then have
\begin{equation}\label{eq:unknown-profile-leading-coefficient}
\int_U\frac{\psi_\kappa(a_{1,0}-a_{2,0})}{S^{\kappa+1}}\Big(\prod_{\ell=1}^{m_*}\alpha_\ell\Big)J_0\,dy=0.
\end{equation}
Taking $\alpha_1\in C_c^\infty(U)$ arbitrary and choosing $\alpha_2,\ldots,\alpha_{m_*}$ equal to $1$ near $\operatorname{supp}\alpha_1$, we obtain
\begin{equation}
\int_U\frac{\psi_\kappa J_0}{S^{\kappa+1}}(a_{1,0}-a_{2,0})\alpha_1\,dy=0.
\end{equation}
Since $\psi_\kappa J_0/S^{\kappa+1}$ does not vanish on $U$, it follows that $a_{1,0}=a_{2,0}$ on $U$.
Set $a_0:=a_{1,0}=a_{2,0}$ on $U$. By the choice of $U$, $a_0$ is
nowhere zero on $\overline U$.

We now proceed inductively for $1\leq d\leq3$. Suppose that $a_{1,s}=a_{2,s}$ for $s<d$. For \(d\leq3\), the coefficient of \(h^{\kappa+d+1}\) involves only the amplitudes \(b_{\ell,m}\) with \(m\leq d\) and \(\partial_t^{\kappa+\rho}\Psi_j|_{t=0}\) with \(\rho\leq d\).
The amplitudes $b_{\ell,0},b_{\ell,1},b_{\ell,2},b_{\ell,3}$ are independent of \(q_j\). Moreover, the adjoint equation and the common boundary jets of orders \(0,\ldots,3\) imply recursively that \(\partial_t^s\Psi_1=\partial_t^s\Psi_2\) for \(s\leq\kappa+3\): the first possible contribution of \(q_j\psi_\kappa\) is to the jet of order \(\kappa+4\). 
Together with the induction condition for \(a_j\), this leaves only \(a_{1,d}-a_{2,d}\) at order
\(h^{\kappa+d+1}\). Applying Corollary~\ref{cor:weighted-extraction} with
\(c=a_1-a_2\) and \(r=d\), we can get
\begin{align}
0=h^{\kappa+d+1}\binom{d+\kappa}{d}\int_U\frac{\psi_\kappa(a_{1,d}-a_{2,d})}{S^{d+\kappa+1}}
\Big(\prod_{\ell=1}^{m_*}\alpha_\ell\Big)J_0\,dy+O(h^{\kappa+d+2}).
\label{eq:initial-profile-jet-coefficient}
\end{align}
Dividing by $h^{\kappa+d+1}$ and letting $h\to0$, the same argument as above gives $a_{1,d}=a_{2,d}$ on $U$. Hence, by induction, we have
\begin{equation}\label{eq:first-four-profile-jets}
a_{1,s}=a_{2,s}\quad\text{on }U,\quad 0\leq s\leq3.
\end{equation}

\medskip
\emph{Step 3.} Assume, for some $r\geq0$, that
\begin{equation}\label{eq:coupled-induction-hypothesis}
q_{1,s}=q_{2,s}\quad(0\leq s<r),\quad a_{1,s}=a_{2,s}\quad(0\leq s<r+4)
\end{equation}
on $U$. We write
\begin{equation}
\Delta^2=\partial_t^4+\sum_{\ell=0}^3A_\ell(y,t,D_y)\partial_t^\ell,
\end{equation}
where the tangential operators $A_\ell$ depend only on the boundary
geometry, and set
\begin{equation}
\psi_{j,m}:=\partial_t^m\Psi_j(y,0).
\end{equation}
Differentiating $L_{q_j}\Psi_j=0$ $m$ times in $t$ gives
\begin{align}
\psi_{j,m+4}=-\sum_{\ell=0}^3\sum_{b=0}^m\binom{m}{b}(\partial_t^bA_\ell)(y,0,D_y)\psi_{j,m-b+\ell}-\sum_{s=0}^m
\binom{m}{s}q_{j,s}\psi_{j,m-s}.
\label{eq:adjoint-recurrence-explicit}
\end{align}
By \eqref{eq:coupled-common-kappa}, we have $\psi_{j,m}=0$ for $0\leq m<\kappa$. While Lemma~\ref{lem:initial-adjoint-jets} gives $\psi_{1,m}=\psi_{2,m}$ for $0\leq m\leq3$.
\eqref{eq:coupled-induction-hypothesis} and \eqref{eq:adjoint-recurrence-explicit} imply recursively that
\begin{equation}
\psi_{1,m}=\psi_{2,m},\quad m<r+\kappa+4.
\end{equation}
At order $m=r+\kappa$, the only new term is the one containing
$q_{j,r}\psi_\kappa$. Hence
\begin{equation}\label{eq:coupled-adjoint-new-jet}
\partial_t^{r+\kappa+4}(\Psi_1-\Psi_2)(y,0)=-\binom{r+\kappa}{\kappa}\psi_\kappa(y)(q_{1,r}-q_{2,r}).
\end{equation}

The transport recursion equation \eqref{eq:transport} shows that $b_m^{(\sigma)}$ depends on the potential only through $q_0,\ldots,q_{m-4}$. Hence the induction hypothesis implies
\begin{equation}
b_{1,m}^{(\sigma)}=b_{2,m}^{(\sigma)} \quad 0\leq m\leq r+3.
\end{equation}
At order $r+4$, Lemma~\ref{lem:two-navier-modes} therefore gives
\begin{equation}\label{eq:coupled-poisson-new-jet}
b_{1,r+4}^{(\sigma)}(y,z)-b_{2,r+4}^{(\sigma)}(y,z)=\bigl(q_{1,r}(y)-q_{2,r}(y)\bigr)\alpha(y)\mu^{-r-4}(y)
\mathfrak a_{r,\sigma}\bigl(\mu(y)z\bigr), \quad \sigma\in\{0,1\}.
\end{equation}
It follows that the coefficient of $h^{r+\kappa+5}$ contains only
\begin{equation}
\delta a_{r+4}:=a_{1,r+4}-a_{2,r+4},\quad\delta q_r:=q_{1,r}-q_{2,r}.
\end{equation}
Again by the induction hypothesis, only $\delta a_{r+4}$ and $\delta q_r$ remain in the coefficient of $h^{r+\kappa+5}$. Taking the exponent sufficiently large such that $L>r+\kappa+6$ gives a remainder of $O(h^{r+\kappa+6})$.

\medskip
\emph{Step 4.} For $m_*\geq3$, we consider the  boundary data $H_{\ell,h}^{(0)}$, and define the normalized decay rates
\begin{equation}
x_\ell(y):=\frac{\mu_\ell(y)}{S(y)},\quad \sum_{\ell=1}^{m_*}x_\ell(y)=1.
\end{equation}
Combining \eqref{eq:coupled-adjoint-new-jet},
\eqref{eq:coupled-poisson-new-jet}, and \eqref{eq:ar-formula}, we obtain
\begin{align}
0=\mathcal T_{1,h}-\mathcal T_{2,h}=h^{r+\kappa+5}\int_U\psi_\kappa\Big(\prod_{\ell=1}^{m_*}\alpha_\ell\Big)J_0\bigl[
\Lambda_r(y)\delta a_{r+4}-a_0G_r(y)\delta q_r\bigr]\,dy+O(h^{r+\kappa+6}),
\label{eq:coupled-integrated-symbol-high-m}
\end{align}
where
\begin{equation}
\Lambda_r(y)=\binom{r+\kappa+4}{\kappa}[S(y)]^{-r-\kappa-5},
\end{equation}
and
\begin{equation}
G_r(y)=[S(y)]^{-r-\kappa-5}\Big[\binom{r+\kappa}{\kappa}+\sum_{d=1}^{r+2}\gamma_{r,d}\binom{d+\kappa}{\kappa} \sum_{\ell=1}^{m_*}[x_\ell(y)]^{d-r-4}\Big].
\end{equation}
Dividing by $h^{r+\kappa+5}$, letting $h\to0$, and using that the boundary amplitudes are arbitrary, we have
\begin{equation}\label{eq:coupled-symbol-equation-high-m}
\Lambda_r(y)\delta a_{r+4}(y)-a_0(y)G_r(y)\delta q_r(y)=0\quad\text{on }U.
\end{equation}

We consider the following two phase configurations
\begin{equation}\label{eq:two-phase-configurations}
\boldsymbol c^{(0)}=\bigl(\underbrace{1,\ldots,1}_{m_*-1},-(m_*-1)\bigr),\quad 
\boldsymbol c^{(\varepsilon)}=\bigl(\varepsilon,\underbrace{1,\ldots,1}_{m_*-2},-(m_*-2+\varepsilon)\bigr),\quad \varepsilon>0,\,\, \varepsilon\neq1.
\end{equation}
For $i\in\{0,\varepsilon\}$, set
\begin{equation}
x_\ell^{(i)}:=\frac{|c_\ell^{(i)}|}{\sum_{j=1}^{m_*}|c_j^{(i)}|}.
\end{equation}
The phase condition $\sum_{\ell=1}^{m_*}c_\ell^{(i)}=0$ gives $x_{m_*}^{(i)}=\sum_{\ell=1}^{m_*-1}x_\ell^{(i)}=\frac12.$
Observing that
\begin{equation}
x_1^{(0)}=\cdots=x_{m_*-1}^{(0)}=\frac{1}{2(m_*-1)},
\end{equation}
\begin{equation}
x_1^{(\varepsilon)}=\frac{\varepsilon}{2(m_*-2+\varepsilon)},\quad x_\ell^{(\varepsilon)}=\frac{1}{2(m_*-2+\varepsilon)},
\quad 2\leq\ell\leq m_*-1,
\end{equation}
$x_1^{(\varepsilon)},\ldots,x_{m_*-1}^{(\varepsilon)}$ are not all equal due to $\varepsilon\neq1$. Since
$x\mapsto x^{-j}$ is strictly convex for $j>0$ and every coefficient in the sum in $\frac{G_r}{\Lambda_r}$ is positive, we obtain
\begin{equation}\label{eq:phase-ratio-separation}
R_r^{(\varepsilon)}>R_r^{(0)},\quad R_r^{(i)}:=\frac{G_r^{(i)}(y)}{\Lambda_r^{(i)}(y)}.
\end{equation}
The ratios $R_r^{(i)}$ are independent of $y$. By \eqref{eq:coupled-symbol-equation-high-m}, we have
\begin{equation}
\delta a_{r+4}-a_0R_r^{(i)}\delta q_r=0, \quad i\in\{0,\varepsilon\}.
\end{equation}
Since $a_0$ is nowhere zero on $U$ and $R_r^{(\varepsilon)}\neq R_r^{(0)}$, these equations imply $\delta a_{r+4}=\delta q_r=0$ on $U$.

\medskip

\emph{Step 5.} We now consider the case $m_*=2$. The condition $\phi_1+\phi_2=0$ implies $\mu_1=\mu_2=:\mu$. Thus, we are required to use the two mode pairs $(\sigma,\tau)=(0,0)$ and $(\sigma,\tau)=(1,0)$ introduced in Lemma~\ref{lem:two-navier-modes}. For each $j=1,2$ and $(\sigma,\tau)\in\{(0,0),(1,0)\}$, let $v_{j,1,h}^{(\sigma)}$ and $v_{j,2,h}^{(\tau)}$ be the corresponding solutions with boundary data $H_{1,h}^{(\sigma)}$ and $H_{2,h}^{(\tau)}$, respectively. Set
\begin{equation}\label{eq:two-mode-tensor}
\mathcal T_{j,h}^{\sigma\tau}:=D^2\Lambda_{q_j,F_j}(0)\bigl[H_{1,h}^{(\sigma)},H_{2,h}^{(\tau)}\bigr]=
\int_\Omega\Psi_j a_jv_{j,1,h}^{(\sigma)}v_{j,2,h}^{(\tau)}\,dx.
\end{equation}
Using the condition \eqref{cond:DN-m-star}, we have $\mathcal T_{1,h}^{\sigma\tau}=\mathcal T_{2,h}^{\sigma\tau}.$
Recall that $\mathfrak p_0(s)=1$ and $\mathfrak p_1(s)=s.$ Define
\begin{align}
\Lambda_r^{\sigma\tau}&:=\frac{1}{\kappa!(r+4)!}\int_0^\infty s^{r+\kappa+4}e^{-2s}\mathfrak p_\sigma(s)\mathfrak p_\tau(s)\,ds,\label{eq:mode-Lambda}\\
\Theta_{r,\mathrm{ad}}^{\sigma\tau}&:=-\frac{\binom{r+\kappa}{\kappa}}{(r+\kappa+4)!}\int_0^\infty s^{r+\kappa+4}e^{-2s}
\mathfrak p_\sigma(s)\mathfrak p_\tau(s)\,ds,
\label{eq:mode-Theta-ad}\\
\Theta_{r,\mathrm P}^{\sigma\tau}&:=\frac{1}{\kappa!}\int_0^\infty s^\kappa e^{-2s}\bigl(\mathfrak a_{r,\sigma}\mathfrak p_\tau+\mathfrak p_\sigma\mathfrak a_{r,\tau}\bigr)(s)\,ds,\quad \Theta_r^{\sigma\tau}:=\Theta_{r,\mathrm{ad}}^{\sigma\tau}+\Theta_{r,\mathrm P}^{\sigma\tau}.
\label{eq:mode-Theta-P}
\end{align}

After the change of variables $s=\mu z$, the coefficient of
$h^{r+\kappa+5}$ in $\mathcal T_{1,h}^{\sigma\tau}-\mathcal T_{2,h}^{\sigma\tau}$ is
\begin{align}
0=\mathcal T_{1,h}^{\sigma\tau}-\mathcal T_{2,h}^{\sigma\tau}=h^{r+\kappa+5}\int_U\mu^{-r-\kappa-5}\psi_\kappa\alpha_1\alpha_2J_0\bigl[\Lambda_r^{\sigma\tau}\delta a_{r+4}
+a_0\Theta_r^{\sigma\tau}\delta q_r\bigr]\,dy +O(h^{r+\kappa+6}).
\label{eq:coupled-integrated-symbol-quadratic}
\end{align}
Dividing by $h^{r+\kappa+5}$, letting $h\to0$, we can derive
\begin{equation}\label{eq:coupled-symbol-equation-quadratic}
\Lambda_r^{\sigma\tau}\delta a_{r+4}+a_0\Theta_r^{\sigma\tau}\delta q_r=0\quad\text{on }U.
\end{equation}
For $k\in\mathbb N_0$, it follows from \eqref{eq:ar-formula} that
\begin{equation}\label{eq:ar-moment-identity}
\int_0^\infty s^ke^{-2s}\mathfrak a_r(s)\,ds=-\sum_{d=1}^{r+2}\gamma_{r,d}
\frac{(k+d)!}{d!\,2^{k+d+1}}=-\frac{k!}{2^{r+k+5}}\Big[\binom{r+k+4}{k+2}-(r+k+4)\Big].
\end{equation}
Using $\mathfrak a_{r,0}=\mathfrak a_r$ and $\mathfrak a_{r,1}=(r+1)\mathfrak a_{r+1},$ we obtain
\begin{equation}\label{eq:quadratic-ratio-kappa1}
\frac{\Theta_r^{10}}{\Lambda_r^{10}}-\frac{\Theta_r^{00}}{\Lambda_r^{00}}=\frac{(r+1)(r+7)(r+8)}{12(r+5)}>0,
\quad \kappa=1,
\end{equation}
and
\begin{equation}\label{eq:quadratic-ratio-kappa3}
\frac{\Theta_r^{10}}{\Lambda_r^{10}}-\frac{\Theta_r^{00}}{\Lambda_r^{00}}=\frac{(r+1)(r+9)
\bigl((r+7)^3+11(r+7)+60\bigr)}{60(r+7)(r+6)(r+5)}>0,\quad \kappa=3.
\end{equation}
Moreover, $\Lambda_r^{00}, \Lambda_r^{10}>0.$ Since $a_0\neq0$ on $U$, the determinant of the two equations
\eqref{eq:coupled-symbol-equation-quadratic}, corresponding to $(\sigma,\tau)=(0,0)$ and $(1,0)$, is nonzero on $U$. Therefore, $\delta a_{r+4}=\delta q_r=0$ on $U$.

\medskip
\emph{Step 6.} We complete the induction. In both cases $m_*\geq3$ and $m_*=2$, the two independent equations
give $\delta a_{r+4}=\delta q_r=0$ on $U$. Again using the density argument, together with \eqref{eq:first-four-profile-jets}, induction on $r$
proves
\begin{equation}
\partial_t^r a_1=\partial_t^r a_2,\quad\partial_t^r q_1=\partial_t^r q_2\quad\text{on }U, \quad r\in\mathbb N_0.
\end{equation}
This is exactly the conclusion \eqref{eq:boundary-jet-conclusion}. 

Finally, assume that $a_1$ and $a_2$ are nowhere zero on $\Gamma$. For nonempty set $V\Subset\Gamma$, Lemma~\ref{lem:AB} provides $U_V\Subset V$ on which the previous argument applies. Hence, we can get
\begin{equation}
\partial_t^r(q_1-q_2)=\partial_t^r(a_1-a_2)=0\quad\text{on }U_V,\quad r\in\mathbb N_0.
\end{equation}
For each fixed $r\in\mathbb N_0$, the zero sets of $\partial_t^r(q_1-q_2)|_\Gamma$ and $\partial_t^r(a_1-a_2)|_\Gamma$ meet every nonempty open subset of $\Gamma$ and are therefore dense. Both functions are continuous on $\Gamma$, so we have
\begin{equation}
\partial_t^r(q_1-q_2)=\partial_t^r(a_1-a_2)=0\quad\text{on }\Gamma.
\end{equation}
Therefore, the proof of Theorem~\ref{thm:boundary-jets} is finished.
\end{proof}

\subsection{Proof of the global uniqueness result}
In this subsection, we will give the proof to Theorem~\ref{thm:main} and Corollary~\ref{cor:known-q-one-component}, respectively.

Once the potential is fixed, the Taylor coefficients of the nonlinearity can be recovered using a single linearized operator.
For $\sigma\in\{0,1\}$, let $\mathscr S_q^{(\sigma)}$ denote the set of complexified solutions of $L_qv=0$ arising from smooth boundary data in $\mathcal X^{(\sigma)}$. That is,
\begin{equation}
\Delta v|_{\partial\Om}=0\quad {\rm for}\,\, \sigma=0,\quad v|_{\partial\Om}=0\quad {\rm for}\,\, \sigma=1.
\end{equation}

\begin{proposition}\label{prop:multiplier}
Let $\Gamma_0\Subset\Gamma$ be a nonempty relatively open set. Suppose that $q\in L^\infty(\Om;\R)$ satisfies Assumption {\bf (A)} and is $C^\infty$ up to the boundary in a neighborhood of $\overline\Gamma_0$. Let $\Psi_q$ be the corresponding adjoint function. Suppose that $b\in C^\infty(\overline\Om)$ is real analytic in $\Om$ and admits a real-analytic extension across $\overline\Gamma_0$. Fix $N\geq2$ and $\sigma\in\{0,1\}$. If
\begin{equation}\label{eq:multiplier-identity}
\int_\Om \Psi_q b\,v_1\cdots v_N\,dx=0
\end{equation}
for every $v_1,\ldots,v_N\in\mathscr S_q^{(\sigma)}$, then $b=0$ in $\Om$.
\end{proposition}

\begin{proof}
Applying the local nondegeneracy argument in Lemma~\ref{lem:AB} on $\Gamma_0$, there is a nonempty connected set  $U\Subset\Gamma_0$ and $\kappa\in\{1,3\}$ such that
\begin{equation}
\Psi_q(y,t)=\frac{\psi_\kappa(y)}{\kappa!}t^\kappa+O(t^{\kappa+1}),\quad \psi_\kappa(y)\neq0, \quad y\in U.
\end{equation}
As in the proof of \emph{Step 1} in Theorem~\ref{thm:boundary-jets}, we choose a function $\rho$ in boundary coordinates such that $d\rho\neq0$ on $U$. We set $\phi_1=\cdots=\phi_{N-1}=\rho,\, \phi_N=-(N-1)\rho,$ and write
\begin{equation}
\mu_j=|d\phi_j|_{g_0},\quad S_N=\sum_{j=1}^N\mu_j.
\end{equation}
Given $\alpha_j\in C_c^\infty(U)$, let $v_{j,h}^{(\sigma)}$ be the semiclassical solutions constructed by Lemma~\ref{lem:two-navier-modes}, with phase $\phi_j$ and boundary amplitude $\alpha_j$.

We next show inductively that $b_r:=\partial_t^rb(\cdot,0)$ vanishes on $U$ for $r\geq0$. Suppose that $b_s=0$ on $U$ for $s<r$. Let $K\Subset U$ contain the supports of $\alpha_1,\ldots,\alpha_N$. Taylor's formula therefore gives
\begin{equation}
b(y,t)=\frac{b_r(y)}{r!}t^r+O_K(t^{r+1}),\quad y\in K,\quad t\to0^+,
\end{equation}
where $b_r(y)=\partial_t^rb(y,0)$ and the remainder is uniform for $y\in K$. We insert the solutions $v_{j,h}^{(\sigma)}=V_{j,h}^{(M)}+R_{j,h}^{(M)}$ into \eqref{eq:multiplier-identity}, and choose $M$ sufficiently large so that $\|R_{j,h}^{(M)}\|_{W^{4,p}(\Om)}=O(h^L)$ for some $L>r+\kappa+2$. Since the $V_{j,h}^{(M)}$ are uniformly
bounded in $L^\infty$, using Sobolev embedding and H\"older's inequality, we can obtain
\begin{equation}
\int_\Om\Psi_qb\Big(\prod_{j=1}^Nv_{j,h}^{(\sigma)}-\prod_{j=1}^NV_{j,h}^{(M)}\Big)dx=O(h^L).
\end{equation}
The induction hypothesis and the leading amplitudes $\alpha_j\mathfrak p_\sigma(\mu_jz)$ now yield, by Proposition~\ref{prop:coefficient-extraction},
\begin{align}
0=\int_\Om\Psi_q b\prod_{j=1}^Nv_{j,h}^{(\sigma)}\,dx=h^{r+\kappa+1}\int_U\psi_\kappa b_rC_{r,\kappa,N}^{(\sigma)}
\Big(\prod_{j=1}^N\alpha_j\Big)J_0\,dy+O(h^{r+\kappa+2}),
\label{eq:multiplier-jet-asymptotic}
\end{align}
where
\begin{equation}
C_{r,\kappa,N}^{(\sigma)}=\frac{1}{r!\kappa!}\int_0^\infty z^{r+\kappa}e^{-S_Nz}\prod_{j=1}^N\mathfrak p_\sigma(\mu_jz)\,dz.
\end{equation}
Using $\mathfrak p_0(s)=1$ and $\mathfrak p_1(s)=s$, we have
\begin{equation}\label{eq:one-mode-moment-coefficient}
C_{r,\kappa,N}^{(0)}=\binom{r+\kappa}{r}S_N^{-r-\kappa-1},\quad C_{r,\kappa,N}^{(1)}=\frac{(r+\kappa+N)!}{r!\kappa!}\frac{\prod_{j=1}^N\mu_j}{S_N^{r+\kappa+N+1}}.
\end{equation}              
Both coefficients are strictly positive on $U$.

Dividing \eqref{eq:multiplier-jet-asymptotic} by $h^{r+\kappa+1}$ and letting $h\to0$, taking $\alpha_1$ to be arbitrary and choosing $\alpha_2,\ldots,\alpha_N$ equal to one near its support, we obtain $b_r=0$ on $U$, since $\psi_\kappa$ does not vanish and $J_0>0$. This completes the induction.

For any $y_0\in U$, the identities $b_r=0$ on $U$ imply
\begin{equation}
\partial_y^\gamma\partial_t^r b(y_0,0)=0\quad\text{for all }\gamma\in\mathbb N_0^{n-1}\text{ and }r\geq0.
\end{equation}
Since the coordinate map is a local diffeomorphism, all Euclidean partial derivatives of $b$ vanish in a neighborhood of $(y_0,0)$. The real-analytic extension of $b$ therefore has a zero Taylor series at $(y_0,0)$ and since $\Om$ is connected, the  analyticity  gives $b=0$ in $\Om$.
\end{proof}

With Proposition~\ref{prop:multiplier} at hand, we are now in a position to prove Corollary~\ref{cor:known-q-one-component}.

\begin{proof}[Proof of Corollary~\ref{cor:known-q-one-component}]
Fix $\sigma\in\{0,1\}$. By \eqref{eq:one-component-map-equality}, there exists $\varepsilon_\sigma>0$ such that
\begin{equation}\label{equality:DN1=DN2}
\Lambda_{q,F_1}(H)=\Lambda_{q,F_2}(H), \quad H\in\mathcal X^{(\sigma)},\quad \|H\|_{\mathcal X}<\varepsilon_\sigma.
\end{equation}
Moreover the two maps are $C^\infty$ near zero by Proposition~\ref{pro:well-posedness-nonlinear}. Therefore, for
every $N\geq1$ and $H_1,\ldots,H_N\in\mathcal X^{(\sigma)}$, we have
\begin{equation}
D^N\Lambda_{q,F_1}(0)[H_1,\ldots,H_N]=D^N\Lambda_{q,F_2}(0)[H_1,\ldots,H_N].
\end{equation}
In what follows, the identities for the mixed variations (higher-order linearization) and $D^N\Lambda_{q,F_j}(0)$ are understood as the unique complex multilinear extensions of their real counterparts.

For $j=1,2$ and $N\geq2$, set $a_{j,N}(x):=\partial_u^NF_j(x,0).$
We prove $a_{1,N}=a_{2,N}$ in $\Om$ by induction on $N$. Assume that $a_{1,k}=a_{2,k}$ for $2\leq k<N$.
For fixed boundary data $H_1,\ldots,H_N\in\mathcal X^{(\sigma)}$, set $I_N=\{1,\ldots,N\}$, and let $u_B^{(j)}$ denote the
mixed variation solution associated with the boundary inputs indexed by $\varnothing\neq B\subseteq I_N$.
The first variations coincide because they solve the linear equation $L_qv=0$ with the same Navier boundary data. Thus, we write $v_i:=u_{\{i\}}^{(1)}=u_{\{i\}}^{(2)}$. We claim that
\begin{equation}\label{eq:common-lower-variations}
u_B^{(1)}=u_B^{(2)}, \quad \varnothing\neq B\subsetneq I_N.
\end{equation}
Indeed, if $2\leq |B|<N$, every term in \eqref{eq:partition} contains a coefficient $a_{j,|\pi|}$ with $|\pi|\leq |B|<N$ and variations $u_C^{(j)}$ with $|C|<|B|$. Induction on $|B|$ and uniqueness for the homogeneous difference problem prove
\eqref{eq:common-lower-variations}.

Set $W_N:=u_{I_N}^{(1)}-u_{I_N}^{(2)}$. By \eqref{eq:common-lower-variations}, all terms in \eqref{eq:partition} agree except the term corresponding to the partition $\{\{1\},\ldots,\{N\}\}$, whose difference is $(a_{1,N}-a_{2,N})v_1\cdots v_N$. Thus, we have
\begin{equation}\label{eq:Fm-source}
\begin{cases}
L_qW_N=(a_{1,N}-a_{2,N})v_1\cdots v_N&\text{in }\Om,\\
W_N=\Delta W_N=0&\text{on }\partial\Om.
\end{cases}
\end{equation}
The condition \eqref{equality:DN1=DN2} gives $\partial_\nu W_N(x_0)=0$. Applying \eqref{eq:transposition} to $W_N$ yields
\begin{equation}\label{eq:Fm-pair}
\int_\Om \Psi_q(a_{1,N}-a_{2,N})v_1\cdots v_N\,dx=0
\end{equation}
for all $v_1,\ldots,v_N\in\mathscr S_q^{(\sigma)}$. Since $a_{1,N}-a_{2,N}$ satisfies the conditions of Proposition~\ref{prop:multiplier}, it follows that $a_{1,N}=a_{2,N}$ in $\Om$. The induction is complete.

For $x\in\Om$, the real-analytic function $u\mapsto F_1(x,u)-F_2(x,u)$ has a zero Taylor series at $u=0$. It therefore vanishes throughout the connected interval $I$. Hence, $F_1=F_2$ in $\Om\times I$ as desired.
\end{proof}

We proceed to prove Theorem~\ref{thm:main}.

\begin{proof}[Proof of Theorem~\ref{thm:main}]
By \eqref{eq:map-equality} and Proposition~\ref{prop:local-wellposedness}, the two measurement maps satisfy
\begin{equation}
D^N\Lambda_{q_1,F_1}(0)=D^N\Lambda_{q_2,F_2}(0),\quad N\geq1.
\end{equation} 
Put $m_j=m(F_j)$ and $a_{j,m_j}=\partial_u^{m_j}F_j(\cdot,0)$. By \eqref{eq:lower-mixed-vanish}, we have
\begin{equation}
D^k\Lambda_{q_j,F_j}(0)=0, \quad 2\leq k<m_j.
\end{equation}
On the other hand, $D^{m_j}\Lambda_{q_j,F_j}(0)$ is nonzero. To see this, since $a_{j,m_j}|_\Gamma\not\equiv0$, there is a nonempty open set $V\Subset\Gamma$ on which $a_{j,m_j}$ does not vanish. By Lemma~\ref{lem:AB}, we can choose a connected
$U\Subset V$ and $\kappa_j\in\{1,3\}$ such that $\psi_{j,\kappa_j}\neq0$ on $U$. Choose $\rho$ with $d\rho\neq0$ on $U$, and set $\phi_1=\cdots=\phi_{m_j-1}=\rho$ and $\phi_{m_j}=-(m_j-1)\rho$. Let $\alpha\in C_c^\infty(U)$ be nonzero and nonnegative, and put
\begin{equation}
H_{\ell,h}=(e^{i\phi_\ell/h}\alpha,0),\quad S=\sum_{\ell=1}^{m_j}|d\phi_\ell|_{g_0}.
\end{equation}
The calculation of leading-order in \eqref{eq:unknown-profile-leading-asymptotic}, applied to
$(q_j,F_j)$, gives
\begin{equation}
D^{m_j}\Lambda_{q_j,F_j}(0)
[H_{1,h},\ldots,H_{m_j,h}]
=
h^{\kappa_j+1}C_j+O(h^{\kappa_j+2}),
\end{equation}
where
\begin{equation}
C_j=\int_U\frac{\psi_{j,\kappa_j}(y)a_{j,m_j}(y,0)} {[S(y)]^{\kappa_j+1}}[\alpha(y)]^{m_j}J_0(y)\,dy\neq0.
\end{equation}
Together with \eqref{eq:lower-mixed-vanish}, this shows that
\begin{equation}
m_j=\min\{k\geq2:D^k\Lambda_{q_j,F_j}(0)\neq0\}.
\end{equation}
It then follows that $m_1=m_2=:m_*.$

The conditions of Theorem~\ref{thm:boundary-jets} now hold. It provides a nonempty $U\Subset\Gamma$ such that
\begin{equation}\label{eq:all-normal-jets}
\partial_\nu^r(q_1-q_2)=0,\quad \partial_\nu^r(a_{1,m_*}-a_{2,m_*})=0 \quad\text{on }U,\quad r=0,1,2,\ldots.
\end{equation}
Similarly, in boundary normal coordinates, \eqref{eq:all-normal-jets} and tangential differentiation give
\begin{equation}
\partial_y^\alpha\partial_t^r(q_1-q_2)(y,0)=0,\quad y\in U,\ \alpha\in\mathbb N_0^{n-1},\ r\geq0.
\end{equation}
Thus, the real-analytic extension assumption on $q_1,q_2$ implies that $q_1-q_2$ vanishes on $U$. By analyticity, it follows that $q_1=q_2:=q$ in $\Omega$. By \eqref{eq:map-equality}, for either $\sigma\in\{0,1\}$, we have
\begin{equation}
\Lambda_{q,F_1}(H)=\Lambda_{q,F_2}(H)\quad \text{for all sufficiently small } H\in\mathcal X^{(\sigma)}.
\end{equation}
Applying Corollary~\ref{cor:known-q-one-component}, we have $F_1=F_2$ on $\Om\times I$. Therefore the proof is complete.
\end{proof}

\section{Finite-dimensional potentials and stability results} \label{sec:finite-dimensional}
This section is devoted to proving Theorems~\ref{thm:finite-dimensional-uniqueness}, \ref{thm:local-sensitivity}, \ref{thm:finite-dimensional-stability}, and Corollary~\ref{cor:one-component-local-stability}, under the assumption that the linear potential $q$ belongs to certain finite-dimensional spaces. We begin with the following lemma:
\begin{lemma}\label{lem:finite-jet-truncation}
Let $E\subset L^\infty(\Omega;\mathbb R)$ be an $N$-dimensional space that is determined by its normal derivatives on $\Gamma$ in the sense of Definition~\ref{def:jet-visible}. For $k\geq0$, set
\begin{equation}
E_k:=\bigl\{e\in E:\partial_\nu^r e|_\Gamma=0\,\, {\rm for}\,\, 0\leq r\leq k\bigr\}.
\end{equation}
Then there exists an integer $R\geq0$ such that $E_R=\{0\}.$ Moreover, for any basis $e_1,\ldots,e_N$ of $E$, there
exist an integer $L\geq N$ and pairs $(y_\ell,r_\ell)\in\Gamma\times\{0,\ldots,R\}$, $1\leq\ell\leq L$, such that the matrix $\mathbb J_E=(\partial_\nu^{r_\ell}e_j(y_\ell))_{1\leq\ell\leq L,1\leq j\leq N}$ has column rank $N$.
\end{lemma}

\begin{proof}[Proof of Lemma~\ref{lem:finite-jet-truncation}]
The spaces $E_k$ form a decreasing sequence of subspaces of $E$.  Since $E_{k+1}\subseteq E_k$ and $E$ is finite dimensional, the sequence $(E_k)_{k\geq0}$ eventually stabilizes. Hence, for some $R\geq0$, we have $E_R=\bigcap_{k\geq0}E_k.$
The intersection is $\{0\}$ by \eqref{eq:jet-visibility}.

For $y\in\Gamma$ and $0\leq r\leq R$, define
\begin{equation}
\ell_{y,r}\in E^*,\quad \ell_{y,r}(e):=\partial_\nu^r e(y).
\end{equation}
The kernel of the functionals $\{\ell_{y,r}:y\in\Gamma,\ 0\leq r\leq R\}$ is $E_R=\{0\}$. Let
\begin{equation}
W:=\operatorname{span}\{\ell_{y,r}:y\in\Gamma,\ 0\leq r\leq R\}\subset E^*.
\end{equation}
Noting that the annihilator of $W$ in $E$ is $E_R=\{0\}$, finite dimensionality gives $W=E^*$. Since $\dim E^*<\infty$, there exist $(y_\ell,r_\ell)\in\Gamma\times\{0,\ldots,R\}$, $1\leq\ell\leq L$, such that
\begin{equation}
E^*=\operatorname{span}\{\ell_{y_\ell,r_\ell}:1\leq\ell\leq L\}.
\end{equation}
Their values on the basis $e_1,\ldots,e_N$ form the matrix
$\mathbb J_E$, which therefore has column rank $N$.
\end{proof}

We next prove Theorem~\ref{thm:finite-dimensional-uniqueness}.

\begin{proof}[Proof of Theorem~\ref{thm:finite-dimensional-uniqueness}]
Set $m_j=m(F_j)$. By the same argument as in the proof of Theorem~\ref{thm:main}, upon replacing $\Gamma$ with $\Gamma_a$, it follows that $m_j$ is the smallest integer $k\geq2$ satisfying $D^k\Lambda_{q_j,F_j}(0)\neq0$. Hence, $m_1=m_2=:m_*$.
Moreover, since $a_{F_1}(y)\neq0$ and $a_{F_2}(y)\neq0$ for $y\in\Gamma_a$, Theorem~\ref{thm:boundary-jets}, applied with $\Gamma_a$ in place of $\Gamma$, yields
\begin{equation}\label{eq:finite-all-jets}
\partial_\nu^r(q_1-q_2)|_{\Gamma_a}=0,\quad r\geq0.
\end{equation}
Using $q_1-q_2\in E$, Lemma~\ref{lem:finite-jet-truncation} gives $q_1=q_2=:q$ in $\Omega$.

It remains to recover $F$. Since $q$ is smooth near $\Gamma_a$, the argument in the proof of Corollary~\ref{cor:known-q-one-component} applies with any nonempty $\Gamma_0\Subset\Gamma_a$ in place of $\Gamma$.  Furthermore, the higher-order linearization identities and \eqref{eq:transposition} give, successively for $N\geq2$,
\begin{equation}
\int_\Omega \Psi_q\bigl(\partial_u^NF_1(x,0)-\partial_u^NF_2(x,0)\bigr) v_1\cdots v_N\,dx=0
\end{equation}
for each fixed $\sigma\in\{0,1\}$ and all $v_1,\ldots,v_N\in\mathscr S_q^{(\sigma)}$. Proposition~\ref{prop:multiplier} then yields
\begin{equation}
\partial_u^NF_1(\cdot,0)=\partial_u^NF_2(\cdot,0)\quad\text{in }\Omega,\quad N\geq2.
\end{equation}
Analyticity condition on $F_1,F_2$ gives $F_1=F_2$ on $\Omega\times I$.
\end{proof}

The proof of Theorem~\ref{thm:local-sensitivity} is given below.

\begin{proof}[Proof of Theorem~\ref{thm:local-sensitivity}]
Define
\begin{equation}
\mathscr R(u,H,q):=\bigl(\Delta^2u+qu-F(x,u),\,u|_{\partial\Omega}-H_0,\,\Delta u|_{\partial\Omega}-H_1\bigr).
\end{equation}
Its derivative with respect to $u$ at $(0,0,q_*)$ is
\begin{equation}
D_u\mathscr R(0,0,q_*)v=\bigl(L_{q_*}v,\,v|_{\partial\Omega},\,\Delta v|_{\partial\Omega}\bigr),
\end{equation}
which is an isomorphism by Lemma~\ref{lem:full-navier-isomorphism}. The analytic implicit-function theorem therefore provides a jointly analytic solution map $(q,H)\mapsto u_{q,H}$ in a neighborhood of $(q_*,0)$ in $(q_*+E)\times\mathcal X_{\mathbb R}$.

We choose $\delta,\rho>0$ sufficiently small such that the map $(q,H)\mapsto u_{q,H}$ is defined for $\|q-q_*\|_E<2\delta$ and
$\|H\|_{\mathcal X}\leq2\rho$. After decreasing $\delta$ and $\rho$, the analyticity of the
solution map gives a uniform bound for $D_q^2u_{q,H}$ when
$q\in B_\delta(q_*)$ and $H\in\overline B_\rho$. If
$q,q+e\in B_\delta(q_*)$, Taylor's formula yields
\begin{equation}
u_{q+e,H}-u_{q,H}-D_qu_{q,H}[e]=\int_0^1(1-s)D_q^2u_{q+se,H}[e,e]\,ds.
\end{equation}
Hence, we have
\begin{equation}
\sup_{H\in\overline B_\rho}
\|u_{q+e,H}-u_{q,H}-D_qu_{q,H}[e]\|_{W^{4,p}(\Omega)}
\leq C\|e\|_E^2.
\end{equation}
Define a map 
\begin{equation}
\mathcal U:\mathcal V\longrightarrow C_b\bigl(\overline B_\rho;W^{4,p}(\Omega)\bigr),\quad\bigl(\mathcal U(q)\bigr)(H):=u_{q,H},
\end{equation}
where $\mathcal V$ is a neighborhood of $q_*$ in $q_*+E$. Thus, for each fixed $q$, $\mathcal U(q)$ maps
$H\in\overline B_\rho$ to the corresponding solution $u_{q,H}$. The uniform estimate shows that $\mathcal U(q)$ is Fr\'echet differentiable as a map into $C_b(\overline B_\rho;W^{4,p}(\Omega))$, with derivative $e\mapsto D\mathcal U(q)[e]$, where $\bigl(D\mathcal U(q)[e]\bigr)(H)=D_qu_{q,H}[e]$. The uniform bound for $D_q^2u_{q,H}$ makes this derivative continuous in $q$, hence the map is $C^1$. Since
\begin{equation}
\mathscr D_F(q)(H)=\partial_\nu\bigl(\mathcal U(q)(H)\bigr)(x_0),
\end{equation}
and $u\mapsto\partial_\nu u(x_0)$ is a continuous linear functional on $W^{4,p}(\Omega)$, the $C^1$ regularity of $\mathcal U$ implies that $q\mapsto\mathscr D_F(q)$ is $C^1$ with values in $\mathscr Y_\rho$.

Fix $H\in\overline B_\rho$ and set $u_H=u_{q_*,H}$. Using $u_{q_*,0}=0$ and the solution map is $C^1$, we have
\begin{equation}
\sup_{H\in\overline B_\rho}\|u_H\|_{C^0(\overline\Omega)}\leq C\rho.
\end{equation}
The condition $\partial_uF(x,0)=0$ therefore implies that $\|\partial_uF(\cdot,u_H)\|_{L^\infty(\Omega)}$ tends to zero
uniformly in $H\in\overline B_\rho$ as $\rho\to0$. At $(q_*,H)$, denote by
\begin{equation}
\mathcal L_H:=\Delta^2+q_*-\partial_uF(x,u_H).
\end{equation}
For given $H\in\overline B_\rho$ and $f\in L^p(\Omega)$, the problem
\begin{equation}
\mathcal L_Hw=f\quad\text{in }\Omega,\quad w=\Delta w=0\quad\text{on }\partial\Omega
\end{equation}
has a unique solution $w\in W^{4,p}(\Omega)$, with a $W^{4,p}$ estimate uniform in $H$. Let $\Phi_H$ be the adjoint function associated with $\mathcal L_H$. For $e\in E$, let $z_{e,H}:=D_qu_{q,H}\big|_{q=q_*}[e].$ Differentiating the nonlinear equation for $u_{q,H}$ with respect to $q$ at $q=q_*$ in the direction $e$, while keeping $H$ fixed, yields
\begin{equation}
\begin{cases}
\mathcal L_Hz_{e,H}=-e\,u_H
&\text{in }\Omega,\\
z_{e,H}=\Delta z_{e,H}=0
&\text{on }\partial\Omega.
\end{cases}
\end{equation}
Applying \eqref{eq:transposition} to the potential $q_*-\partial_uF(\cdot,u_H)$ for $\mathcal L_H$ yields
\begin{equation}\label{eq:sensitivity-formula}
\bigl(D_q\mathscr D_F(q_*)e\bigr)(H)=\partial_\nu z_{e,H}(x_0)=-\int_\Omega e\,u_H\Phi_H\,dx.
\end{equation}

Let $N=\dim E$ and set $\mathscr V:=D_q\mathscr D_F(q_*)E\subset\mathscr Y_\rho.$ By \eqref{eq:local-sensitivity-rank}, we have $\dim\mathscr V=N$. Moreover there are $H_1,\ldots,H_L\in\overline B_\rho$ such that
\begin{equation}
\mathcal E:\mathscr V\longrightarrow\mathbb R^L, \quad \mathcal E(\varphi):=\bigl(\varphi(H_1),\ldots,\varphi(H_L)\bigr)
\end{equation}
is injective. Consequently, the map
\begin{equation}
G(q):=\bigl(\Lambda_{q,F}(H_1),\ldots,\Lambda_{q,F}(H_L)\bigr)
\end{equation}
has an injective derivative at $q_*$, since $DG(q_*)|_E=\mathcal E\circ D_q\mathscr D_F(q_*)|_E.$ Thus $DG(q_*)|_E$ has rank $N$. Selecting $N$ components of $G$ gives a map $\widetilde G:q_*+E\to\mathbb R^N$ whose derivative at $q_*$ is an isomorphism. The finite-dimensional inverse function theorem gives, after restricting to a smaller neighborhood $\mathcal U_*$ of $q_*$, a Lipschitz inverse for $G$. The stability result \eqref{eq:local-finite-input-stability} follows from the equivalence of norms on $E$. Therefore, the proof of Theorem~\ref{thm:local-sensitivity} is finished.
\end{proof}

We next prove Corollary~\ref{cor:one-component-local-stability} and Corollary~\ref{cor:automatic-sensitivity-rank}.

\begin{proof}[Proof of Corollary~\ref{cor:one-component-local-stability}]
Restriction from $\overline B_\rho$ to $\overline B_\rho^{(\sigma)}$ defines a bounded linear map $\mathscr Y_\rho\to\mathscr Y_\rho^{(\sigma)}$. Consequently, $\mathscr D_F^{(\sigma)}$ is $C^1$, and its derivative is the restriction of $D_q\mathscr D_F$. The formula~\eqref{eq:sensitivity-formula} therefore holds for $H\in\overline B_\rho^{(\sigma)}$. By the condition \eqref{eq:restricted-sensitivity-rank}, the range of $D_q\mathscr D_F^{(\sigma)}(q_*)|_E$ has dimension $N=\dim E$. By some similar arguments as that in proving Theorem~\ref{thm:local-sensitivity}, we can prove \eqref{eq:restricted-local-stability}.
\end{proof}

\begin{proof}[Proof of Corollary~\ref{cor:automatic-sensitivity-rank}]
Let $e\in E$ belong to the kernel in \eqref{eq:automatic-sensitivity-rank}, and write $\mathcal L(q,H):=\Lambda_{q,F}(H)$. We have $D_q\mathcal L(q_*,H)[e]=0$ for all $H\in\overline B_\rho.$
The map $(q,H)\mapsto\mathcal L(q,H)$ is smooth near $(q_*,0)$.
Differentiating this identity with respect to $H$ at $H=0$
therefore gives
\begin{equation}\label{eq:zero-mixed-parameter-data}
D_qD_H^k\mathcal L(q_*,0) [e;H_1,\ldots,H_k]=0,\quad k=1,m_*,
\end{equation}
for all $H_1,\ldots,H_k\in\mathcal X_{\mathbb R}$. Set $q_\tau=q_*+\tau e$, and let $\Psi_\tau$ be the point
adjoint function associated with $q_\tau$. Write $\Psi=\Psi_0$. Since $L_{q_\tau}=L_{q_*}+\tau e$, we see that
\begin{equation}
\bigl(G_{q_\tau}^{(p)}-G_{q_*}^{(p)}\bigr)f=-\tau G_{q_*}^{(p)}\bigl(e\,G_{q_\tau}^{(p)}f\bigr), \quad f\in L^p(\Omega),
\end{equation}
where $G_{q_\tau}^{(p)}$ and $G_{q_*}^{(p)}$ denote the inverse maps as in Lemma~\ref{lem:Green-identity}.
Thus, for $f\in L^p(\Omega;\mathbb R)$, the point-adjoint identity
\eqref{eq:adjoint-duality-real} gives
\begin{equation}
\int_\Omega(\Psi_\tau-\Psi)f\,dx=-\tau\,\partial_\nu\bigl(G_{q_*}^{(p)}(eG_{q_\tau}^{(p)}f)\bigr)(x_0)
=-\tau\int_\Omega\Psi\,e\,G_{q_\tau}^{(p)}f\,dx=-\tau\int_\Omega G_{q_\tau}^{(p')}(e\Psi)\,f\,dx.
\end{equation}
The last equality follows from the $L^p$--$L^{p'}$ duality argument. Thus,
\begin{equation}
\Psi_\tau-\Psi
=-\tau G_{q_\tau}^{(p')}(e\Psi)
\quad\text{in }L^{p'}(\Omega).
\end{equation}
By continuity of $G_{q_\tau}^{(p')}$ at $\tau=0$, we can get
\begin{equation}\label{eq:linearized-point-adjoint}
\dot\Psi=-G_{q_*}^{(p')}(e\Psi),\quad L_{q_*}\dot\Psi=-e\Psi.
\end{equation}
In particular, $\dot\Psi$ has zero Navier boundary traces. The first Green identity and \eqref{eq:zero-mixed-parameter-data} with $k=1$ give $\partial_\nu\dot\Psi=\partial_\nu\Delta\dot\Psi=0$ on $\Gamma_a$. Hence the first four normal derivatives of $\dot\Psi$ vanish on $\Gamma_a$.

We show inductively that $e_r:=\partial_t^re(\cdot,0)=0$ on $\Gamma_a$ for every $r\geq0$. Suppose this holds for $0\leq s<r$.
By Lemma~\ref{lem:AB}, every nonempty relatively open subset of $\Gamma_a$ contains an open set $U$ on which
\begin{equation}
\Psi(y,t)=\frac{\psi_\kappa(y)}{\kappa!}t^\kappa+O(t^{\kappa+1}),\quad \kappa\in\{1,3\},\quad \psi_\kappa\neq0.
\end{equation}
We use the oscillatory Navier boundary data from the proof of Theorem~\ref{thm:boundary-jets}, with phases whose sum vanishes
and amplitudes supported in $U$. Let $b_{\tau,\ell,m}$ denote the order-$m$ amplitude in the asymptotic construction of the solution to $L_{q_\tau}v=0$ with Navier boundary data $H_{\ell,h}$, and set $\dot b_{\ell,m}:=\left.\partial_\tau b_{\tau,\ell,m}\right|_{\tau=0}.$
Assume inductively that $e_s=0$ on $U$ for $0\leq s<r$. Differentiating the equation $L_{q_*}\dot\Psi=-e\Psi$ successively in the normal variable $t$, and using the vanishing lower-order normal derivatives, gives
\begin{equation}
\partial_t^j\dot\Psi(y,0)=0,\quad 0\leq j<r+\kappa+4,\quad\partial_t^{r+\kappa+4}\dot\Psi(y,0)=-\binom{r+\kappa}{\kappa}\psi_\kappa(y)e_r(y).
\end{equation}
Since the boundary data are independent of $\tau$ and $e_s=0$ for $s<r$, differentiating the transport system
\eqref{eq:transport}--\eqref{eq:transport-bc2} for the $\ell$th solution gives $\dot b_{\ell,m}=0$ for $m<r+4$. At order $r+4$, the only new source is $-e_r\alpha_\ell z^r/r!$. The homogeneous boundary conditions and Lemma~\ref{lem:polynomial-normal} therefore give
\begin{equation}
\dot b_{\ell,r+4}(y,z)=e_r(y)\alpha_\ell(y)[\mu_\ell(y)]^{-r-4}\mathfrak a_r\bigl(\mu_\ell(y)z\bigr).
\end{equation}

For $q_\tau$ and the boundary input $H_{\ell,h}$, denote by $V_{\tau,\ell,h}^{(M)}$ the approximation solution of order $M>0$ as in \eqref{eq:WKB}, and write the corresponding solution as $v_{\tau,\ell,h}=V_{\tau,\ell,h}^{(M)}+R_{\tau,\ell,h}^{(M)}$. Set
\begin{equation}
\mathcal E_{\tau,\ell,h}^{(M)}:=L_{q_\tau}V_{\tau,\ell,h}^{(M)},\quad \dot V_{\ell,h}^{(M)}:=\partial_\tau V_{\tau,\ell,h}^{(M)}\big|_{\tau=0}.
\end{equation}
Recalling $q_\tau=q_*+\tau e$, the derivative of the equation is
\begin{equation}
\dot{\mathcal E}_{\ell,h}^{(M)}:=\partial_\tau\mathcal E_{\tau,\ell,h}^{(M)}\big|_{\tau=0}=L_{q_*}\dot V_{\ell,h}^{(M)}
+eV_{0,\ell,h}^{(M)}.
\end{equation}
For any sufficiently small $\tau$, the transport equations
make the coefficients of $h^0,\ldots,h^M$ vanish in both
$h^4L_{q_\tau}V_{\tau,\ell,h}^{(M)}$ and
$h^2\Delta V_{\tau,\ell,h}^{(M)}|_{\partial\Omega}$.
Since these coefficients depend smoothly on $\tau$ and vanish
identically, their $\tau$-derivatives vanish as well. The scaling estimates
\eqref{eq:bl-Wsr-scaling} and
\eqref{eq:boundary-Wsr-scaling} therefore give
\begin{equation}
\|\dot{\mathcal E}_{\ell,h}^{(M)}\|_{L^p(\Omega)}+\left\|\partial_\tau\bigl(\Delta V_{\tau,\ell,h}^{(M)}
|_{\partial\Omega}\bigr)\big|_{\tau=0}\right\|_{W^{2-1/p,p}(\partial\Omega)}\leq C_Mh^{M-3+1/p}.
\end{equation}
Let $\dot R_{\ell,h}^{(M)}:=\partial_\tau R_{\tau,\ell,h}^{(M)}\big|_{\tau=0}.$ We find that
\begin{equation}
\begin{cases}
L_{q_*}\dot R_{\ell,h}^{(M)}=-\partial_\tau\bigl(L_{q_\tau}V_{\tau,\ell,h}^{(M)}\bigr)\big|_{\tau=0}
-eR_{0,\ell,h}^{(M)}&\text{in }\Omega,\\[2mm]
\dot R_{\ell,h}^{(M)}=0,\,\Delta\dot R_{\ell,h}^{(M)}=-\partial_\tau\bigl(\Delta V_{\tau,\ell,h}^{(M)}\bigr)
\big|_{\tau=0}&\text{on }\partial\Omega.
\end{cases}
\end{equation}
By \eqref{eq:R-W4p}, we have
\begin{equation}
\|eR_{0,\ell,h}^{(M)}\|_{L^p(\Omega)}\leq \|e\|_{L^\infty(\Omega)}\|R_{0,\ell,h}^{(M)}\|_{L^p(\Omega)}\leq C_Mh^{M-3+1/p}.
\end{equation}
Applying Lemma~\ref{lem:full-navier-isomorphism}, we obtain
\begin{equation}\label{eq:linearized-poisson-remainder}
\|\dot R_{\ell,h}^{(M)}\|_{W^{4,p}(\Omega)}\leq C_Mh^{M-3+1/p}.
\end{equation}

Denote by
\begin{equation}
\mathcal T_\tau(h):=D_H^{m_*}\Lambda_{q_\tau,F}(0)[H_{1,h},\ldots,H_{m_*,h}].
\end{equation}
By \eqref{eq:zero-mixed-parameter-data}, $\left.\partial_\tau\mathcal T_\tau(h)\right|_{\tau=0}=0$. As in the derivation of
\eqref{eq:coupled-integrated-symbol-high-m} from \eqref{eq:coupled-adjoint-new-jet} and \eqref{eq:coupled-poisson-new-jet} in Theorem~\ref{thm:boundary-jets}, and taking $M$ such that $M-3+1/p\geq r+\kappa+6$, we obtain
\begin{equation}\label{eq:linearized-rank-coefficient}
0=-h^{r+\kappa+5}\int_U\psi_\kappa(y)a(y,0)e_r(y)J_0(y)\Big(\prod_{\ell=1}^{m_*}\alpha_\ell(y)\Big)\mathcal C_{r,\kappa}(y)\,dy+O(h^{r+\kappa+6}),
\end{equation}
where $S(y)=\sum_{\ell=1}^{m_*}\mu_\ell(y)$ and
\begin{equation}
\mathcal C_{r,\kappa}(y)=\binom{r+\kappa}{\kappa}[S(y)]^{-r-\kappa-5}+\sum_{\ell=1}^{m_*}\sum_{d=1}^{r+2} \gamma_{r,d}\binom{d+\kappa}{\kappa}\mu_\ell(y)^{d-r-4}[S(y)]^{-d-\kappa-1}>0.
\end{equation}
Here $\gamma_{r,d}>0$ are the coefficients in Lemma~\ref{lem:ar}.

Divide \eqref{eq:linearized-rank-coefficient} by $h^{r+\kappa+5}$ and let $h\to0$. Choose one $\alpha_\ell\in C_c^\infty(U)$ arbitrarily and the remaining amplitudes to equal one near its support. Since $\psi_\kappa a(\cdot,0)J_0\mathcal C_{r,\kappa}$ is nowhere zero on $U$, it follows that $e_r=0$ on $U$. If $e_r$ is nonzero at some point of $\Gamma_a$, it would be nonzero on a relatively open subset $V\subset\Gamma_a$. By Lemma~\ref{lem:AB}, there is a nonempty $U\Subset V$ on which $\psi_\kappa$ is nowhere zero for some $\kappa\in\{1,3\}$. The previous argument gives $e_r=0$ on $U$, contradicting $e_r\neq0$ on $V$. Hence $e_r=0$ on $\Gamma_a$. Induction on $r$ and the determining property of
$E$ give $e=0$. The local stability conclusion now immediately follows from Theorem~\ref{thm:local-sensitivity}.
\end{proof}

In the end of this section, we give the proof of Theorem~\ref{thm:finite-dimensional-stability}.

\begin{proof}[Proof of Theorem~\ref{thm:finite-dimensional-stability}]
The assertions are immediate if $\mathcal K=\varnothing$ or $E=\{0\}$. Note that the set of potentials satisfying Assumption {\bf (A)} is open in $L^\infty(\Omega)$, and hence is stable under small $L^\infty$ perturbations. By compactness and the
implicit-function theorem, there are $\rho>0$ and an open neighborhood $\mathcal U$ of $\mathcal K$
in $q_{\mathrm{ref}}+E$ such that $\mathscr D_F:\mathcal U\longrightarrow\mathscr Y_\rho$ is well defined and $C^1$. Theorem~\ref{thm:finite-dimensional-uniqueness}, applied with the same $F$ for both potentials, shows that $\mathscr D_F$ is injective on $\mathcal K$.

Define the function
\begin{equation}\label{eq:omega-definition}
\omega(s):=\sup\{\|q_1-q_2\|_{L^\infty(\Omega)}:q_1,q_2\in\mathcal K,\, \|\mathscr D_F(q_1)-\mathscr D_F(q_2)\|_{\mathscr Y_\rho}\leq s\}.
\end{equation}
It is clear that $\omega(s)$ is finite and nondecreasing. If $\omega(s)\not\to0$ as $s\downarrow0$, then there exist
$\varepsilon_0>0$, a sequence $s_j\downarrow0$, and $q_{1,j},q_{2,j}\in\mathcal K$ such that
\begin{equation}
\|\mathscr D_F(q_{1,j})-\mathscr D_F(q_{2,j})\|_{\mathscr Y_\rho} \leq s_j, \quad\|q_{1,j}-q_{2,j}\|_{L^\infty(\Omega)}
\geq\varepsilon_0.
\end{equation}
By compactness, after passing to subsequences, $q_{1,j}\to q_1$ and $q_{2,j}\to q_2$ in $\mathcal K$.
Continuity gives $\mathscr D_F(q_1)=\mathscr D_F(q_2)$, whereas $\|q_1-q_2\|_{L^\infty(\Omega)}\geq\varepsilon_0$,
contradicting injectivity. Thus \eqref{eq:finite-modulus-stability} follows.

To prove (ii), let $N=\dim E$, we choose a basis $e_1,\ldots,e_N$, and write $q_\theta=q_{\mathrm{ref}}+\sum_{j=1}^N\theta_je_j.$ Denote by $K_E\subset\mathbb R^N$ the parameter set corresponding to $\mathcal K$, and set $T(\theta)=\mathscr D_F(q_\theta)$. The map $T$ is $C^1$ on an open neighborhood of $K_E$.
For each $\theta\in K_E$, Corollary~\ref{cor:automatic-sensitivity-rank} gives $\ker DT(\theta)=\{0\}.$
Fix $\theta\in K_E$. Since its range is an $N$-dimensional
space of functions on $\overline B_\rho$, there exist 
$H_{\theta,1},\ldots,H_{\theta,N}\in\overline B_\rho$ such that $v\mapsto\bigl((DT(\theta)v)(H_{\theta,1}),\ldots,(DT(\theta)v)(H_{\theta,N})\bigr)$ is injective on $\mathbb R^N$. For each $\theta_0\in K_E$, continuity of $DT$ gives a
neighborhood $U_{\theta_0}$ and $c_{\theta_0}>0$ such that
\begin{equation}
\sum_{j=1}^{N}|(DT(\eta)v)(H_{\theta_0,j})|^2\geq \big(c_{\theta_0}|v|\big)^2,\quad
\eta\in U_{\theta_0},\quad v\in\mathbb R^N.
\end{equation}
Choose a finite subcover of $K_E$ and collect its associated
inputs as $H_1^0,\ldots,H_{L_0}^0$. Then
\begin{equation}
G_0(\theta):=\bigl(T(\theta)(H_1^0),\ldots,T(\theta)(H_{L_0}^0)\bigr)
\end{equation}
satisfies the estimate
\begin{equation}\label{eq:uniform-differential-lower-bound}
|DG_0(\theta)v|\geq c_0|v|,\quad \theta\in K_E,\quad v\in\mathbb R^N,
\end{equation}
for some $c_0>0$. Since $K_E$ is compactly contained in the $C^1$ domain of $G_0$, we can choose $r>0$ so that, whenever
$\theta,\eta\in K_E$ and $|\theta-\eta|<r$, 
\begin{equation}
\sup_{0\leq t\leq1}\|DG_0(\theta+t(\eta-\theta))-DG_0(\theta)\|\leq \frac{c_0}{2}.
\end{equation}
The fundamental theorem of calculus now gives
\begin{equation}\label{eq:near-parameter-separation}
|G_0(\theta)-G_0(\eta)|
\geq \frac{c_0}{2}|\theta-\eta|.
\end{equation}

If $P_r:=\{(\theta,\eta)\in K_E^2:|\theta-\eta|\geq r\}$ is nonempty, then compactness of $P_r$ and injectivity of $T$ provide inputs $H_1^1,\ldots,H_{L_1}^1\in\overline B_\rho$ such that
\begin{equation}
G_1(\theta):=\bigl(T(\theta)(H_1^1),\ldots,T(\theta)(H_{L_1}^1)\bigr)
\end{equation}
satisfies $\delta_r:=\min_{(\theta,\eta)\in P_r}|G_1(\theta)-G_1(\eta)|>0.$
Set $G=(G_0,G_1)$ and $D=\operatorname{diam}(K_E)$. Combining \eqref{eq:near-parameter-separation} with the definition of $\delta_r$ gives
\begin{equation}
|G(\theta)-G(\eta)|\geq\min\left\{c_0/2,\delta_r/D\right\}|\theta-\eta|,\quad \theta,\eta\in K_E.
\end{equation}
If $P_r$ is empty, take $G=G_0$. Norm equivalence on $E$, together with
$|G(\theta)-G(\eta)|
\leq\sqrt L\,\|T(\theta)-T(\eta)\|_{\mathscr Y_\rho}$,
completes the proof of \eqref{eq:finite-lipschitz-stability}. Therefore, the proof of Theorem~\ref{thm:finite-dimensional-stability} is complete.
\end{proof}

\section{Examples and counterexamples}\label{sec:examples}
In this section, we present some examples that verify the main assumptions of each theorem stated in Subsection~\ref{sec:setting}, along with certain counterexamples that reveal obstructions to recovering the unknown potential $q$.

\subsection{Examples for an admissible pair $(q,F)$}
In each example below, we assume that $q\geq1$. Thus $q$ satisfies  Assumption {\bf (A)}.

\begin{example}
\label{ex:profile-with-boundary-zeros}
Let $\Omega=B(0,1)=\{x\in\mathbb R^n: |x|<1\}$, $n\geq2$, and $x_0=e_n$. Choose a connected relatively open set $\Gamma\Subset\partial\Omega\setminus\{x_0\}$ containing $-e_n$. For any $m\geq2$, set
\begin{equation}
q(x)=1, \quad F(x,u)=\frac{x_1}{m!}u^m+(2+x_2^2)u^{m+1}.
\end{equation}
The first nonzero nonlinear coefficient is $a_F(x)=x_1$. Although $a_F(-e_n)=0$, its restriction to $\Gamma$ is not
identically zero. This pair $(q,F)$ satisfies the assumptions in Theorem~\ref{thm:main}. It also illustrates that the condition $a_j|_\Gamma\not\equiv0$ in Theorem~\ref{thm:boundary-jets} does not require pointwise nonvanishing on $\Gamma$.
\end{example}

\begin{example}\label{ex:admissible}
Let $\Omega\subset\mathbb R^n$ ($n\geq2$) be a bounded, connected and smooth domain. For given $m\geq2$, let $q(x)=1+|x|^2$
and
\begin{equation}
F(x,u)=\frac{e^{x_1}+x_2^2+2}{m!}u^m+(2+|x|^2)u^{m+1}+\sin(x_2)u^{m+2}.
\end{equation}
It is clear that both coefficients are real analytic, and
\begin{equation}
F(x,0)=\partial_uF(x,0)=0,\quad a_F(x):=\partial_u^mF(x,0)=e^{x_1}+x_2^2+2>0.
\end{equation}
Thus this example satisfies the conditions on coefficients needed in Theorem~\ref{thm:main} for every $m\geq2$.
The function $a_F$ also satisfies the pointwise nonvanishing condition in the final assertion of Theorem~\ref{thm:boundary-jets}.
\end{example}

The following example gives an example of $q$ that is nonsmooth at the observation point $x_0$.
\begin{example}\label{ex:rough-admissible}
Let $\Omega=B(0,1)\subset\mathbb R^n$, $n\geq2$, and $x_0=e_1$. Choose a nonempty connected relatively open set
$\Gamma\Subset\partial\Omega\setminus\{x_0\}$, and set
\begin{equation}
q(x)=1+\sqrt{1-x_1},\quad F(x,u)=-(2+e^{x_2})u^2+(1+x_1^2)u^5.
\end{equation}
The potential belongs to $L^\infty(\Omega)$, is real analytic in $\Omega$, and extends real analytically across
$\overline\Gamma$. The derivative of $q$ is unbounded as $x\to x_0$, so it has no $C^1$
extension across $x_0$. Moreover,
\begin{equation}
a_F(x)=\partial_u^2F(x,0)=-2(2+e^{x_2})
\end{equation}
is nowhere zero. This example satisfies the assumptions on coefficients in Theorem~\ref{thm:main}.
\end{example}

The following example presents a nonlinearity \(F\) whose leading coefficient vanishes on the entire boundary $\partial\Omega$.
\begin{example}\label{ex:known-q-vanishing-profile}
Let $\Omega=B(0,1)$, take the known potential $q=1$, and fix $M\geq1$. Set
\begin{equation}
F(x,u)=(1-|x|^2)^M u^2+e^{x_1}u^3.
\end{equation}
Then $F$ is real analytic, $F(x,0)=\partial_uF(x,0)=0$, and
\begin{equation}
a_F(x):=\partial_u^2F(x,0)=2(1-|x|^2)^M.
\end{equation}
The function $a_F$ vanishes to order $M$ on $\partial\Omega$. Nevertheless, the conditions of Corollary~\ref{cor:known-q-one-component} hold. Thus, when $q$ is known, either family of inputs $(f,0)$ or $(0,g)$ determines $F$, even though its leading coefficient vanishes on the entire boundary.
\end{example}

\subsection{Examples of finite-dimensional potential}\label{sec:finite-dimensional-examples}
In the following, Examples~\ref{ex:localized-finite-dimensional-potential} and \ref{ex:finite-dimensional-potential} give nonanalytic potential families satisfying \eqref{eq:jet-visibility}. Example~\ref{ex:interior-sensitivity} shows that
\eqref{eq:jet-visibility} is not necessary for local stability when the boundary inputs are restricted to either
$(f,0)$ or $(0,g)$. Example~\ref{ex:global-rank} gives a compact potential family 
$\mathcal K=\{\theta|x|^2:0\leq\theta\leq\theta_0\}$ for which the fixed boundary input $H_\varepsilon=(\varepsilon,0)$
yields a stable determination result.

\begin{example}\label{ex:localized-finite-dimensional-potential}
Choose nonempty relatively open sets $V_0,V_1\Subset\Gamma$ with disjoint closures. Let $\chi\in C_c^\infty(\Gamma)$ equal one on $\overline V_0$ and vanish near $\overline V_1$, and let $\chi_1\in C_c^\infty([0,t_0))$ equal one near $t=0$.
In inward normal coordinates $(y,t)$, define
\begin{equation}
e_j(y,t)=\chi(y)\chi_1(t)t^{j-1}, \quad 1\leq j\leq N,
\end{equation}
and extend these functions by zero outside the coordinate neighborhood. Set $E_N=\operatorname{span}\{e_1,\ldots,e_N\}$.

For $e=\sum_{j=1}^N c_je_j$ and $y\in V_0$,
\begin{equation}
\partial_\nu^r e(y)=(-1)^r r!c_{r+1},\quad 0\leq r\leq N-1.
\end{equation}
Thus $E_N$ satisfies \eqref{eq:jet-visibility} on every nonempty set $\Gamma_a\Subset V_0$ (in Lemma~\ref{lem:finite-jet-truncation}, one can take $R=N-1$).  By contrast, every element of $E_N$ and all its normal
derivatives vanish on $V_1$. 

Choose a nonzero $r_0\in C_c^\infty(\Omega)$ with $r_0\geq0$ and set $q_{\mathrm{ref}}=1+r_0$. On a sufficiently small
compact parameter set, let
\begin{equation}
q_\theta=q_{\mathrm{ref}}+\sum_{j=1}^N\theta_je_j\geq\frac12.
\end{equation}
Hence every $q_\theta$ satisfies Assumption {\bf (A)}. For example, with
\begin{equation}
F(x,u)=\frac{2+e^{x_1}}{m!}u^m+h(x)u^{m+1}, \quad m\geq2,
\end{equation}
where $h$ is real analytic on $\mathcal O$, the conditions of Theorems~\ref{thm:finite-dimensional-uniqueness} and
\ref{thm:finite-dimensional-stability} hold on $\Gamma_a$. The reference potential is not real analytic in $\Omega$.
\end{example}

\begin{example}\label{ex:finite-dimensional-potential}
Let $r_0\in C_c^\infty(\Omega)$ be nonzero and nonnegative, and let $q_{\mathrm{ref}}=1+r_0$. Let $E_d$ consist of the
restrictions to $\Omega$ of real polynomials of degree at most $d\in\mathbb N$. For $P\in E_d$ and $y\in\Gamma$, the function
$t\mapsto P(y-t\nu(y))$ is a polynomial of degree at most $d$. If
\begin{equation}
\partial_\nu^jP|_\Gamma=0,\quad 0\leq j\leq d,
\end{equation}
then this polynomial vanishes for every $y\in\Gamma$. Thus $P$ vanishes in an open subset of $\Omega$, and hence $P=0$.
Thus, $E_d$ satisfies \eqref{eq:jet-visibility} on any nonempty relatively open subset of $\Gamma$.

Set $\mathcal K:=\left\{q_{\mathrm{ref}}+P:P\in E_d,\ \|P\|_{L^\infty(\Omega)}\leq\frac12\right\}$, which is compact. Moreover, $q_{\mathrm{ref}}=1+r_0\geq1$, so every
$q\in\mathcal K$ satisfies $q\geq\frac12$ and Assumption {\bf (A)}. Taking, for instance, $F(x,u)=e^{x_1}u^2$ verifies the assumptions of Theorems~\ref{thm:finite-dimensional-uniqueness} and
\ref{thm:finite-dimensional-stability}. This example again
allows a reference potential that is not real analytic in
$\Omega$.
\end{example}

\begin{example}\label{ex:interior-sensitivity}
Let $\Omega=B(0,1)$, $q_*=0$, and $F(x,u)=e^{x_1}u^2$. Let $\Psi_0$ be the point adjoint function for $\Delta^2$. By \eqref{eq:transposition}, it is not identically zero. Choose $B\Subset\Omega$ where $\Psi_0$
has a fixed nonzero sign, and take $0\neq e\in C_c^\infty(B)$ with $e\geq0$.
Set
\begin{equation}
H_0^{(0)}=(1,0),\quad v^{(0)}=1,\quad H_0^{(1)}=(0,1),\quad v^{(1)}=\frac{|x|^2-1}{2n}.
\end{equation}
Both $v^{(\sigma)}$ solve $\Delta^2v=0$ with the indicated Navier boundary data. Differentiating
\eqref{eq:sensitivity-formula} at zero boundary input gives
\begin{equation}
D_H\!\big(D_q\mathscr D_F^{(\sigma)}(q_*)[e]\big)(0)[H_0^{(\sigma)}]=-\int_\Omega e\,v^{(\sigma)}\Psi_{q_*}\,dx\neq0,
\quad \sigma=0,1.
\end{equation}
Thus a sufficiently small fixed input $\varepsilon H_0^{(\sigma)}$ gives local Lipschitz stability
on $E=\operatorname{span}\{e\}$, for either $\sigma=0$ or $\sigma=1$. However, the set $E$ fails to satisfy the condition \eqref{eq:jet-visibility}, since $\operatorname{supp}e\Subset\Omega$ and hence $\partial_\nu^r e|_\Gamma=0$ for every $r\geq0$.
\end{example}

\begin{example}\label{ex:global-rank}
Let $\Omega=B(0,1)\subset\mathbb R^n$, let $x_0\in\partial\Omega$, and set
\begin{equation}
q_{\mathrm{ref}}=0,\quad E=\operatorname{span}\{|x|^2\},\quad q_\theta(x)=\theta|x|^2,\quad F(x,u)=e^{x_1}u^2.
\end{equation}
The space $E$ satisfies \eqref{eq:jet-visibility} with $R=0$. Also, $\partial_u^2F(x,0)=2e^{x_1}>0$, and every
$q_\theta$ with $\theta\geq0$ satisfy Assumption {\bf (A)}.
Fix the boundary input $H_\varepsilon=(\varepsilon,0)$ and write
\begin{equation}
g(\theta,\varepsilon):=\Lambda_{q_\theta,F}(H_\varepsilon).
\end{equation}
At $(\theta,\varepsilon)=(0,0)$, the first variation with respect to $\varepsilon$ is $v=1$. Its derivative with respect to $\theta$ is the solution of
\begin{equation}
\Delta^2z=-|x|^2\quad\text{in }\Omega, \quad z=\Delta z=0\quad\text{on }\partial\Omega.
\end{equation}
Direct calculation gives
\begin{equation}
z(x)=\frac{|x|^2-1}{8n(n+2)}-\frac{|x|^6-1}{24(n+2)(n+4)},\quad\partial_\nu z(x_0)=\frac{1}{n(n+2)(n+4)}>0.
\end{equation}
We can check that
\begin{equation}
\partial_\theta\partial_\varepsilon g(0,0)=\frac{1}{n(n+2)(n+4)}.
\end{equation}
Fixing a sufficiently small $\varepsilon>0$ and then choosing $\theta_0>0$ sufficiently small, smooth dependence gives $\partial_\theta g(\theta,\varepsilon)\geq c>0$ for $0\leq\theta\leq\theta_0$. Thus, with  $\mathcal K=\{q_\theta:0\leq\theta\leq\theta_0\}$, we have
\begin{equation}
\|q_\theta-q_\eta\|_{L^\infty(\Omega)}=|\theta-\eta|\leq c^{-1}\left|\Lambda_{q_\theta,F}(H_\varepsilon)-\Lambda_{q_\eta,F}(H_\varepsilon)\right|.
\end{equation}
Thus the stability estimate in Theorem~\ref{thm:finite-dimensional-stability} holds for this family with the single fixed input
$H_\varepsilon$.
\end{example}

\subsection{Counterexamples of recovering the linear potential}
The following proposition and examples give a nonuniqueness construction for the linear problem and clarify the roles of the boundary nonvanishing and analyticity assumptions.

\begin{proposition}[Nonuniqueness for the linear equation]\label{prop:linear-counterexample}
When $F\equiv0$, there are different real potentials $q_0,q_\varepsilon\in C^\infty(\overline\Omega)$ satisfying Assumption {\bf (A)} such that
\begin{equation}
\Lambda_{q_0,0}(H)=\Lambda_{q_\varepsilon,0}(H) \quad\text{for any }H\in\mathcal X_{\mathbb R}.
\end{equation}
\end{proposition}

\begin{proof}
Set $q_0=0$ and denote its point adjoint function by $\Psi_0$. Since $\Psi_0\not\equiv0$ and is smooth in $\Omega$, there
are a ball $B\Subset\Omega$ and $c_0>0$ such that $|\Psi_0|\geq c_0$ on $\overline B$. Choose a real function $w\in C_c^\infty(B)$ with $\Delta^2w\not\equiv0$. For sufficiently small $\varepsilon\neq0$, we define
\begin{equation}
\Psi_\varepsilon=\Psi_0+\varepsilon w,\quad q_\varepsilon(x)=
\begin{cases}
-\dfrac{\varepsilon\Delta^2w(x)}{\Psi_0(x)+\varepsilon w(x)},&x\in B,\\[5pt]
0,&x\notin B.
\end{cases}
\label{eq:counterexample-potential}
\end{equation}
Then $q_\varepsilon\in C_c^\infty(\Omega;\mathbb R)$ is nonzero, $\|q_\varepsilon\|_{L^\infty}=O(|\varepsilon|)$, and $q_\varepsilon\Psi_\varepsilon=-\varepsilon\Delta^2w.$ Assumption {\bf (A)} for $q_\varepsilon$ follows from that for $q_0$ and openness of invertibility under small $L^\infty$ perturbations.

If $Z\in W^{4,p}(\Omega)$ has homogeneous Navier boundary traces, the transposition identity \eqref{eq:transposition} for $\Psi_0$ and the compact support of $w$ give
\begin{equation}
\int_\Omega\Psi_\varepsilon L_{q_\varepsilon}Z\,dx=\int_\Omega\Psi_0\Delta^2Z\,dx+\varepsilon\int_\Omega\bigl(w\Delta^2Z-(\Delta^2w)Z\bigr)\,dx=\partial_\nu Z(x_0).
\end{equation}
Thus $\Psi_\varepsilon$ is the point adjoint function for $q_\varepsilon$.

Let $v_0$ and $v_\varepsilon$ solve the linear problems for $q_0$ and $q_\varepsilon$, respectively, with the same boundary input $H$. Then $Z=v_\varepsilon-v_0$ has homogeneous Navier boundary traces and satisfies the equation $L_{q_\varepsilon}Z=-q_\varepsilon v_0$. Therefore, we have
\begin{equation}
\Lambda_{q_\varepsilon,0}(H)-\Lambda_{q_0,0}(H)=-\int_\Omega\Psi_\varepsilon q_\varepsilon v_0\,dx=\varepsilon\int_\Omega(\Delta^2w)v_0\,dx=\varepsilon\int_\Omega w\Delta^2v_0\,dx=0.
\end{equation}
However, $q_\varepsilon\neq q_0$, and the claim follows.
\end{proof}

\begin{remark}
\label{rem:counterexample-scope}
Proposition~\ref{prop:linear-counterexample} shows
nonuniqueness in the unrestricted class of smooth potentials when $F\equiv0$. It does not address nonzero nonlinearities. The constructed potential $q_\varepsilon$ is nonzero and compactly supported in $\Omega$, so it is not real analytic in $\Omega$. The counterexample also lies in the one-dimensional class $E=\operatorname{span}\{q_\varepsilon\}$, which contains
$q_0=0$ and $q_\varepsilon$. This does not contradict Theorem~\ref{thm:finite-dimensional-uniqueness}: every
element of $E$ has zero normal derivatives on $\Gamma$, so \eqref{eq:jet-visibility} fails. 
\end{remark}

\begin{example}\label{ex:flat-potential}
On $\Omega=B(0,1)$, define
\begin{equation}
r(x)=
\begin{cases}
\exp\!\left(-\dfrac{1}{1-|x|^2}\right),&|x|<1,\\[6pt]
0,&|x|=1.
\end{cases}
\end{equation}
Then $r\in C^\infty(\overline\Omega)$ is real analytic in
$\Omega$, whereas
\begin{equation}
\partial_\nu^k r|_{\partial\Omega}=0\quad\text{for every }k\geq0.
\end{equation}
Thus, interior analyticity and smoothness up to the boundary do not suffice to determine a potential from its boundary normal derivatives of all orders. 
\end{example}

\begin{remark}
Example~\ref{ex:flat-potential} illustrates a limitation of the analytic continuation argument and does not establish nonuniqueness for the nonlinear inverse problem.  Thus its extension to general smooth coefficients remains open.
\end{remark}

\begin{example}
\label{ex:flat-nonlinearity}
Let $H\not\equiv0$ be real analytic in the spatial
variable and set
\begin{equation}
F(x,s)
=H(x)
\begin{cases}
e^{-1/s^2},&s\neq0,\\
0,&s=0.
\end{cases}
\end{equation}
This function is smooth and nonzero, but
\begin{equation}
\partial_s^kF(x,0)=0 \quad\text{for every }k\geq0.
\end{equation}
Hence, for a fixed nonresonant potential $q$, the maps $\Lambda_{q,F}$ and $\Lambda_{q,0}$ have identical Fr\'echet derivatives of every order at zero. Thus, smoothness in $s$ alone does not ensure that $F$ is determined by the Taylor coefficients of the measurement map. This example does not establish equality of the two maps for nonzero boundary inputs.
\end{example}

\section{Numerical experiments}\label{sec:numerical-experiment}

\begin{figure}[htbp]
\centering
\includegraphics[width=0.8\textwidth]{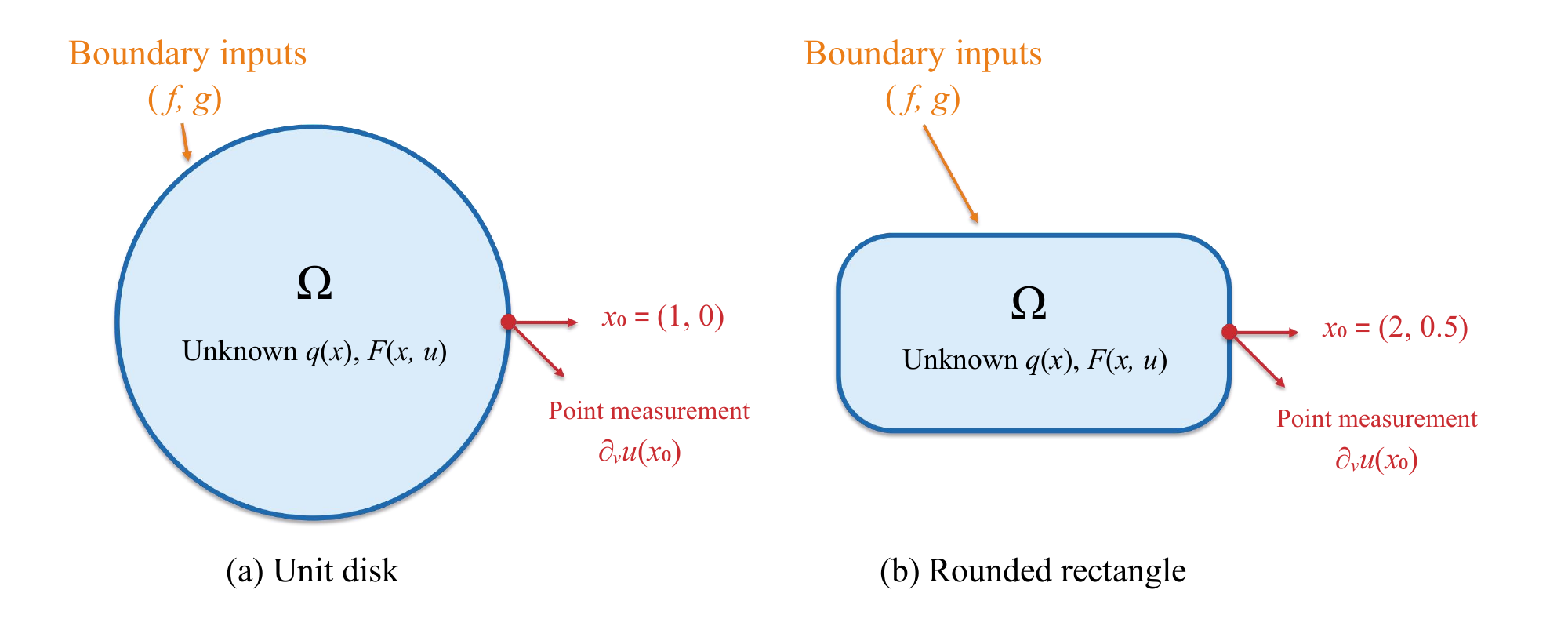}
\caption{\small Experimental setup on the unit disk and the analytic rounded rectangle.
Navier boundary inputs $(f,g)$ are applied on the boundary, and
$\partial_\nu u(x_0)$ is observed at the marked point in each domain.}
\label{fig:num-setup}
\end{figure}

We numerically investigate the simultaneous recovery of $q$ and $F$ from measurements at a fixed boundary point on two domains: a disk and a smoothly rounded rectangle. Four experiments address the recovery in a radial coefficient class, spatially varying coefficients and mixed nonlinearities, identification of the first nonvanishing nonlinear order, and the choice of Navier boundary inputs. Figure~\ref{fig:num-setup} roughly shows the considered domains, boundary inputs, and observation points. The source code and data-generation scripts are available\footnote{https://github.com/Zhikun0416/Recovering-Coefficients-from-Boundary-Point-Measurements}, which contains the full experimental settings.
We solve the nonlinear problem \eqref{eq:model} on the following two domains:
\begin{equation}\label{eq:num-setting}
\begin{aligned}
\text{Case 1:}\quad &\Om=\{x\in\R^2:|x|<1\},&&x_0=(1,0),\\
\text{Case 2:}\quad &\Om=\Phi(\{w\in\mathbb C:|w|<1\}),&&x_0=(2,0.5).
\end{aligned}
\end{equation}
The boundary of Case~1 is clearly analytic. Then we identify $\mathbb R^2$ with $\mathbb C$, the rounded rectangle is defined by the conformal map
\begin{equation}\label{eq:num-rounded-map}
\Phi(w)=1+\frac{i}{2}+s\int_0^w\frac{dt}{\sqrt{1-2\kappa\cos(2\alpha)t^2+\kappa^2t^4}},\quad \kappa=0.9.
\end{equation}
The constants $s>0$ and $0<\alpha<\pi/4$ are chosen so that $\Phi(1)=2+i/2$ and $\Phi(i)=1+i$. Numerically,
$\alpha\approx0.128963$ and $s\approx0.622391$. The map has the symmetries of a rounded rectangle with bounding box $[0,2]\times[0,1]$. It is analytic for $|w|<\kappa^{-1/2}$. Moreover, writing $z_\pm=\kappa e^{\pm2i\alpha}w^2$, we have
\begin{equation}
1+\frac{w\Phi''(w)}{\Phi'(w)}=\frac12\Big(\frac{1+z_+}{1-z_+}+\frac{1+z_-}{1-z_-}\Big).
\end{equation}
Its real part is positive for $|w|<\kappa^{-1/2}$, so $\Phi$ is univalent on $|w|<\kappa^{-1/2}$ by the convexity criterion. In particular, the boundary of Case~2 is also analytic.

We set $(\xi,\eta)=(x_1,x_2)$ in Case 1 and $(\xi,\eta)=(x_1-1,x_2-0.5)$ in Case 2, and write
$\rho^2=\xi^2+\eta^2$.  Each input $H_k=(f_k,g_k)$ gives the scalar measurement data $d_k=\partial_\nu u_{H_k}(x_0)$ at the same observation point. For every reference model considered later, $q^\dagger\geq0.9$ and the first nonzero nonlinear Taylor coefficient $\partial_u^{m_*}F^\dagger(\cdot,0)$ is strictly positive on $\overline\Omega$, where $m_*\in\{2,3\}$ depends on the model. Thus, the homogeneous Navier realization of
$\Delta^2+q^\dagger$ has trivial kernel and is nonresonant. For both geometries, the boundary inputs are represented by finite Fourier series in the angular parameter $\theta$, using the boundary
parametrizations $e^{i\theta}$ and $\Phi(e^{i\theta})$, respectively.
The modal amplitudes decrease with frequency, and a
geometry-dependent normalization places the two Navier components on
comparable scales. Candidate inputs are selected by an A-optimal
criterion for the scaled parameter vector, evaluated at a fixed
nominal model. The same criterion is applied to all three input
strategies in Example~4. Training and validation inputs are selected
from independent candidate sets. In Example~3, both the quadratic and
cubic models are included in the input design.

We impose the Navier boundary data through a biharmonic lifting $u_H$ satisfying
\begin{equation}
\Delta^2u_H=0\quad\text{in }\Omega,\quad (u_H,\Delta u_H)|_{\partial\Omega}=(f,g).
\end{equation}
We solve for $w=u-u_H$, which has homogeneous Navier boundary data, by Picard iteration applied to the equivalent coupled Poisson system. In Case~2, we set $U=u\circ\Phi$, $V=(\Delta u)\circ\Phi$, and $G=|\Phi'|^2$. The
transformed equations are
\begin{equation}
\Delta_w U=GV,\quad \Delta_w V=G\bigl(F(\Phi(w),U)-q(\Phi(w))U\bigr).
\end{equation}
We discretize Case~1 directly on the unit disk. For Case~2, we pull the problem back to the unit disk and use the same polar finite-volume scheme for the transformed equations, retaining the conformal factor $G$. The physical point observation in Case~2 is
\begin{equation}
\partial_\nu u(x_0)=\frac{\partial_rU(1,0)}{|\Phi'(1)|}.
\end{equation}
The radial derivative $\partial_rU(1,0)$ is evaluated using a second-order finite-difference formula at the boundary. In both geometries, manufactured-solution tests exhibit second-order convergence of the numerical solution and the computed point measurement.

The data are generated on a grid finer than both inversion grids. For all fits performed on the refined inversion grid, including noisy-data and candidate-model fits, we use the corrected response
\begin{equation}\label{eq:num-correction}
\widetilde M_h=M_h+\beta(M_h-M_{2h}),\quad \beta=1/3\big(1-h_d^2/h^2\big).
\end{equation}
Here $M_t$ denotes the response computed with grid width $t$, $M_0$ its continuum limit, and $h_d$ the data-grid width. If, along the nested grids, $M_t=M_0+c_2t^2+o(t^2),$ then $\widetilde M_h=M_0+c_2h_d^2+o(h^2).$
The corresponding parameter Jacobian is
\begin{equation}
D_p\widetilde M_h
=(1+\beta)D_pM_h-\beta D_pM_{2h}.
\end{equation}

For the unknown parameter vector $p$, we solve the bound-constrained least-squares problem
\begin{equation}\label{eq:num-minimization}
\min_{\ell\leq p\leq b}\frac{1}{2}\sum_{k=1}^{N}\big|\widetilde M_h(p;H_k)-d_k^\delta\big|^2,
\end{equation}
where the bounds are imposed componentwise. Parameter derivatives are approximated by finite differences in Example~1 and computed by the discrete adjoint method in Examples~2--4. We denote the training Jacobian by $J$.
Let $P=\dim p$. All reference parameter components appearing in the following relative errors are nonzero. We use the following relative errors:
\begin{equation}\label{eq:num-errors}
E_{\rm RMS}=\Big(\frac{1}{P}\sum_{j=1}^{P}\Big|\frac{\widehat p_j-p_j^\dagger}{p_j^\dagger}\Big|^2\Big)^{1/2},\quad E_\infty=\max_{1\leq j\leq P}\frac{|\widehat p_j-p_j^\dagger|}{|p_j^\dagger|},\quad E_2(w)=\frac{\|\widehat w-w^\dagger\|_{L^2(\Omega)}}{\|w^\dagger\|_{L^2(\Omega)}},
\end{equation}
and $P=\dim p.$ The quantity $E_2(w)$ is reported only when $\|w^\dagger\|_{L^2(\Omega)}>0$. All relative errors are expressed as percentages.

Noisy observations are obtained by perturbing the fine-grid data:
\begin{equation}
d_k^\delta=d_k+\delta S Z_k,\quad Z_k\overset{\mathrm{i.i.d.}}{\sim}N(0,1).
\end{equation}
We use $S=1$ in all experiments and generate ten independent noise realizations at each noise level. For comparisons between the two geometries and, in Example~4, among the three input strategies, we use the same standard Gaussian draws for corresponding measurements. Training and validation data are perturbed independently.

\subsection{Simultaneous recovery of radial coefficients}\label{subsec:num-radial}
Example~1 considers simultaneous reconstruction of the potential and
the quadratic nonlinear coefficient in the four-parameter class
\begin{equation}\label{eq:num-radial-model}
q=q_0+q_2\rho^2,\quad F(x,u)=(a_0+a_2\rho^2)u^2.
\end{equation}
The reference coefficients are
\begin{equation}
q^\dagger=1.1+0.8\rho^2,\quad F^\dagger(x,u)=(2.4+1.5\rho^2)u^2.
\end{equation}
The training set contains 48 inputs, including separate and simultaneous excitations of the two Navier boundary components, with positive and negative amplitudes paired. The reference pair satisfies the regularity, nonresonance (Assumption {\bf (A)}), and leading-coefficient nonvanishing condition of Theorem~\ref{thm:main}. The potential space $E=\operatorname{span}\{1,\rho^2\}$ also satisfies the normal-derivative determination condition in Theorem~\ref{thm:finite-dimensional-uniqueness}.

Table~\ref{tab:num-radial} shows the estimates obtained on the refined inversion grid. At the reported precision, the largest relative parameter error rounds to $0.00\%$ in Case~1 and is $0.05\%$ in Case~2. Thus, the selected single-point measurements recover all four parameters accurately in both geometries.

\begin{table}[htbp]
\centering
\small
\caption{\small Example~1: reconstruction on the refined inversion grid from
fine-grid data.}
\label{tab:num-radial}
\begin{tabular}{cccccc}
\toprule
& & \multicolumn{2}{c}{Case 1} & \multicolumn{2}{c}{Case 2}\\
Parameter & True value & Recovered & Error & Recovered & Error\\
\midrule
$q_0$ & 1.1 & 1.1000 & 0.00 & 1.1000 & 0.00\\
$q_2$ & 0.8 & 0.8000 & 0.00 & 0.8002 & 0.03\\
$a_0$ & 2.4 & 2.4000 & 0.00 & 2.3999 & 0.00\\
$a_2$ & 1.5 & 1.5000 & 0.00 & 1.5007 & 0.05\\
\bottomrule
\end{tabular}
\end{table}

\subsection{Spatially varying coefficients and mixed nonlinearities} \label{subsec:num-spatial}
Example~2 introduces nonradial spatial dependence in the potential and the quadratic coefficient, together with a cubic term. The seven unknown coefficients are parametrized by
\begin{equation}\label{eq:num-mixed-model}
q=q_0+q_2\rho^2+q_x\xi,\quad F(x,u)=a(x)u^2+b_3u^3,
\end{equation}
with $a(x)=a_0+a_2\rho^2+a_y\eta.$ The reference model is
\begin{equation}
q^\dagger=1.1+0.8\rho^2+0.2\xi,\quad a^\dagger=2.4+1.5\rho^2+0.3\eta,\quad b_3^\dagger=0.15.
\end{equation}
The training set comprises 96 boundary inputs. The modes in the two Navier boundary components are varied independently, and positive and negative excitation amplitudes are paired. In the theoretical results, the nonlinear Taylor coefficients at $u=0$ are
\begin{equation}
\partial_u^2F(x,0)=2a(x),\quad \partial_u^3F(x,0)=6b_3.
\end{equation}
The tables instead list the parameters of the representation
\eqref{eq:num-mixed-model}.

Table~\ref{tab:num-spatial} compares reconstructions from the same fine-grid data. Relative to the coarse-grid fits, the refined-grid fits with correction~\eqref{eq:num-correction} reduce the maximum relative parameter error from $1.11\%$ to a value rounding to $0.00\%$ in Case~1, and from $26.68\%$ to $0.25\%$ in Case~2. The relative $L^2(\Omega)$ errors $E_2(q)$ and $E_2(a)$ also decrease. 

\begin{table}[htbp]
\centering
\small
\caption{Example~2: reconstruction errors on the base and refined inversion grids.}
\label{tab:num-spatial}
\begin{tabular}{ccccc}
\toprule
Case & Inversion grid & $E_\infty$ (\%) & $E_2(q)$ (\%) & $E_2(a)$ (\%)\\
\midrule
Case 1 & Base    & 1.11  & 0.20 & 0.05\\
Case 1 & Refined & 0.00  & 0.00 & 0.00\\
Case 2 & Base    & 26.68 & 2.88 & 0.71\\
Case 2 & Refined & 0.25  & 0.03 & 0.01\\
\bottomrule
\end{tabular}
\end{table}

\subsection{Identification of the leading nonlinear order}\label{subsec:num-order}
Example~3 examines whether the first nonzero order of the nonlinearity can be recovered from the point measurements. The cubic reference model is
\begin{equation}\label{eq:num-cubic-model}
q^\dagger=1.1+0.8\rho^2+0.2\xi,\quad F^\dagger(x,u)=(2.4+1.5\rho^2+0.3\eta)u^3.
\end{equation}
For the quadratic reference model, $u^3$ is replaced by $u^2$,
with the potential and spatial coefficient unchanged. For each data set,
we fit three candidates: $F=0$, $F=a(x)u^2$, and $F=c(x)u^3$, with
\begin{equation}
q\in\operatorname{span}\{1,\rho^2,\xi\},\quad a,c\in\operatorname{span}\{1,\rho^2,\eta\}.
\end{equation}

We fit the linear model first. Its estimated potential initializes the quadratic and cubic fits, with their nonlinear coefficients initially set to zero. Each model is fitted to 96 training observations using the corrected response \eqref{eq:num-correction}. The validation set consists of 48 boundary inputs not used in fitting. For a fitted parameter vector $\widehat p$, the residual is
\begin{equation}\label{eq:num-validation-rms}
R_{\rm val}(\widehat p)=\Big(\frac{1}{48}\sum_{k=1}^{48}\big|\widetilde M_h(\widehat p;H_k^{\rm val})
-d_k^{\rm val}\big|^2\Big)^{1/2}.
\end{equation}
We select the candidate minimizing $R_{\rm val}$. This procedure identifies the generating model in all four noiseless tests (Table~\ref{tab:num-order-selection}). For the cubic reference model, with training and validation noise generated independently, the cubic candidate is selected in all ten trials at each noise level considered
(Table~\ref{tab:num-order-noise}). 

\begin{table}[htbp]
\centering
\small 
\caption{\small Example~3: validation RMS \eqref{eq:num-validation-rms} for the three candidate models.
The smallest value in each row is in bold.}
\label{tab:num-order-selection}
\begin{tabular}{ccccc}
\toprule
Case & True model & Linear & Quadratic & Cubic\\
\midrule
Case 1 & Quadratic & $6.6226\times10^{-2}$ & $\mathbf{9.0177\times10^{-9}}$ & $6.6046\times10^{-2}$\\
Case 1 & Cubic & $6.9755\times10^{-2}$ & $6.9755\times10^{-2}$ & $\mathbf{9.5387\times10^{-9}}$\\
Case 2 & Quadratic & $1.7579\times10^{-2}$ & $\mathbf{2.7164\times10^{-7}}$ & $1.7578\times10^{-2}$\\
Case 2 & Cubic & $1.2617\times10^{-2}$ & $1.2617\times10^{-2}$ & $\mathbf{1.8140\times10^{-7}}$\\
\bottomrule
\end{tabular}
\end{table}

\begin{table}[htbp]
\centering
\small
\caption{\small Example~3: model-selection counts over ten trials for the cubic reference model.}
\label{tab:num-order-noise}
\begin{tabular}{ccccc}
\toprule
Case & Noise level $\delta$ & Linear & Quadratic & Cubic\\
\midrule
Case 1 & $10^{-5}$ & 0 & 0 & 10\\
Case 1 & $10^{-4}$ & 0 & 0 & 10\\
Case 1 & $10^{-3}$ & 0 & 0 & 10\\
Case 2 & $10^{-5}$ & 0 & 0 & 10\\
Case 2 & $10^{-4}$ & 0 & 0 & 10\\
Case 2 & $10^{-3}$ & 0 & 0 & 10\\
\bottomrule
\end{tabular}
\end{table}

\subsection{Effects of the Navier boundary inputs}\label{subsec:num-inputs}
For $m_*=2$, the proof of Theorem~\ref{thm:boundary-jets} uses variations in both Navier boundary components to distinguish the potential from the leading quadratic coefficient. Example~4 compares reconstruction using $f$-only, $g$-only, and mixed $(f,g)$ inputs for
\begin{equation}\label{eq:num-input-model}
q^\dagger=1.1+0.8\rho^2+0.2\xi, \quad F^\dagger(x,u)=(2.4+1.5\rho^2+0.3\eta)u^2.
\end{equation}
In each domain, every strategy uses 96 training measurements to estimate the six coefficients, subject to the same weighted boundary-energy budget. For $H=(f,g)$, the energy is
\begin{equation}\label{eq:num-input-energy}
\|H\|_w^2=\frac{1}{|\partial\Omega|}\int_{\partial\Omega}\Big(|f|^2+\frac{|g|^2}{s_\Omega^2}\Big)\,ds.
\end{equation}
The geometry-dependent factor $s_\Omega$ is the ratio of the
Frobenius norms of the discrete biharmonic lifting matrices for the
normalized $f$- and $g$-boundary modes. For each strategy, the inputs
are selected separately at the same nominal model using the same
A-optimal design criterion. In the mixed-input strategy, the design
also determines the energy allocation and relative sign of the two
nonzero boundary components.

All six discrete training Jacobians in Table~\ref{tab:num-inputs} have full column rank. Among the selected designs, mixed inputs yield the largest $\sigma_{\min}(J)$ in each geometry, indicating the greatest worst-direction sensitivity of the linearized discrete parameter-to-data map. This local comparison does not imply the smallest reconstruction error: in Case~2, the $f$-only strategy gives the smallest noiseless error.

\begin{table}[htbp]
\centering
\small
\caption{Example~4: base-grid sensitivities and refined-grid errors under equal measurement and input-energy budgets.}
\label{tab:num-inputs}
\begin{tabular}{cccccc}
\toprule
Case & Inputs & $\sigma_{\min}(J)$ & $\operatorname{cond}(J)$
& $E_{\rm RMS}$ (\%) & $E_\infty$ (\%)\\
\midrule
Case 1 & only $f$ & $1.4354\times10^{-3}$ & 109.1314 & 0.00 & 0.00\\
Case 1 & only $g$ & $5.8266\times10^{-4}$ & 544.8580 & 0.00 & 0.00\\
Case 1 & both $f,g$ & $8.0233\times10^{-3}$ & 41.6280 & 0.00 & 0.00\\
Case 2 & only $f$ & $1.1275\times10^{-3}$ & 37.1002 & 0.02 & 0.02\\
Case 2 & only $g$ & $1.0132\times10^{-3}$ & 69.6849 & 0.07 & 0.15\\
Case 2 & both $f,g$ & $1.7990\times10^{-3}$ & 36.1479 & 0.11 & 0.25\\
\bottomrule
\end{tabular}
\end{table}

\begin{figure}[ht!]
\centering
\includegraphics[width=0.88\textwidth]{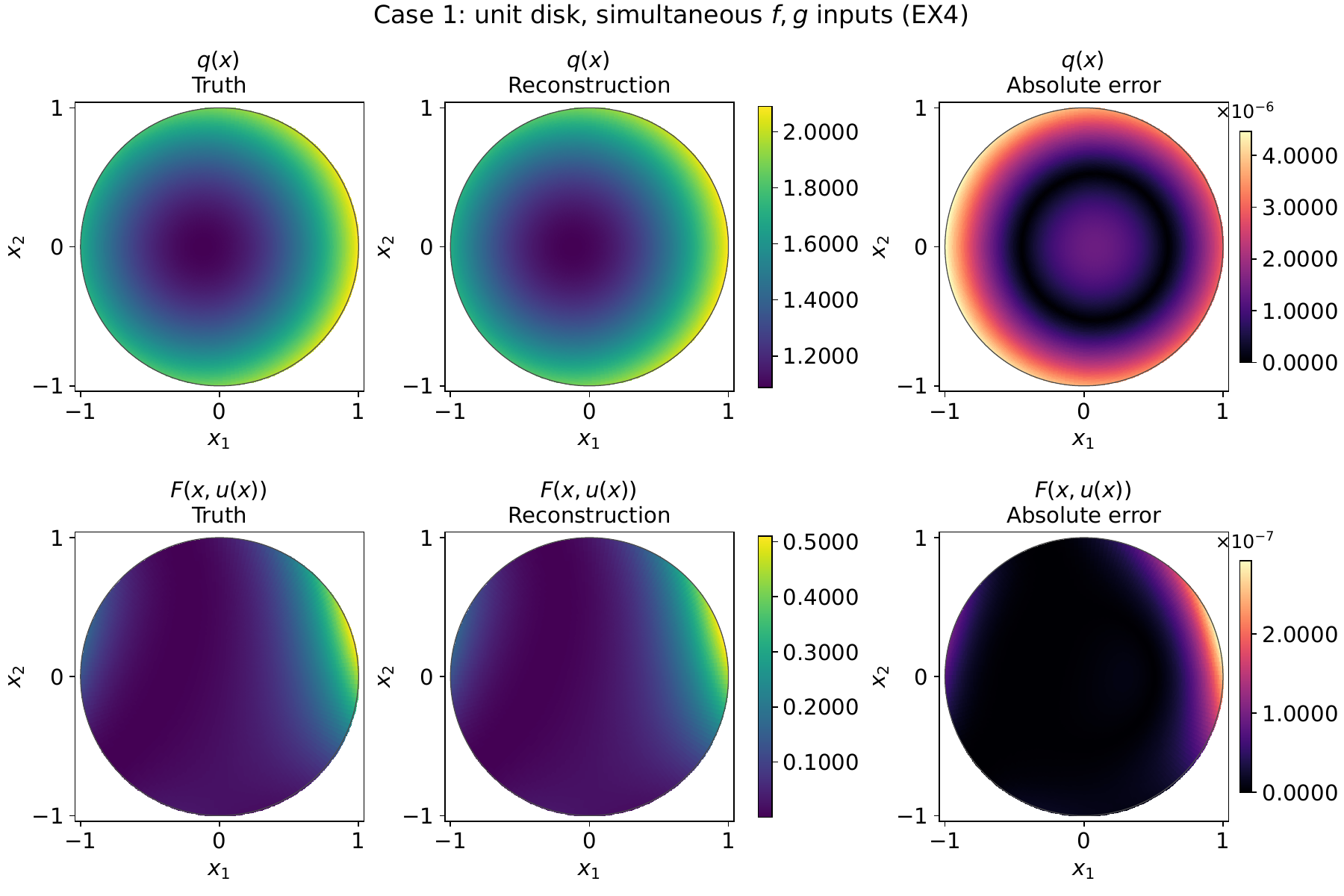}\\

\includegraphics[width=0.88\textwidth]{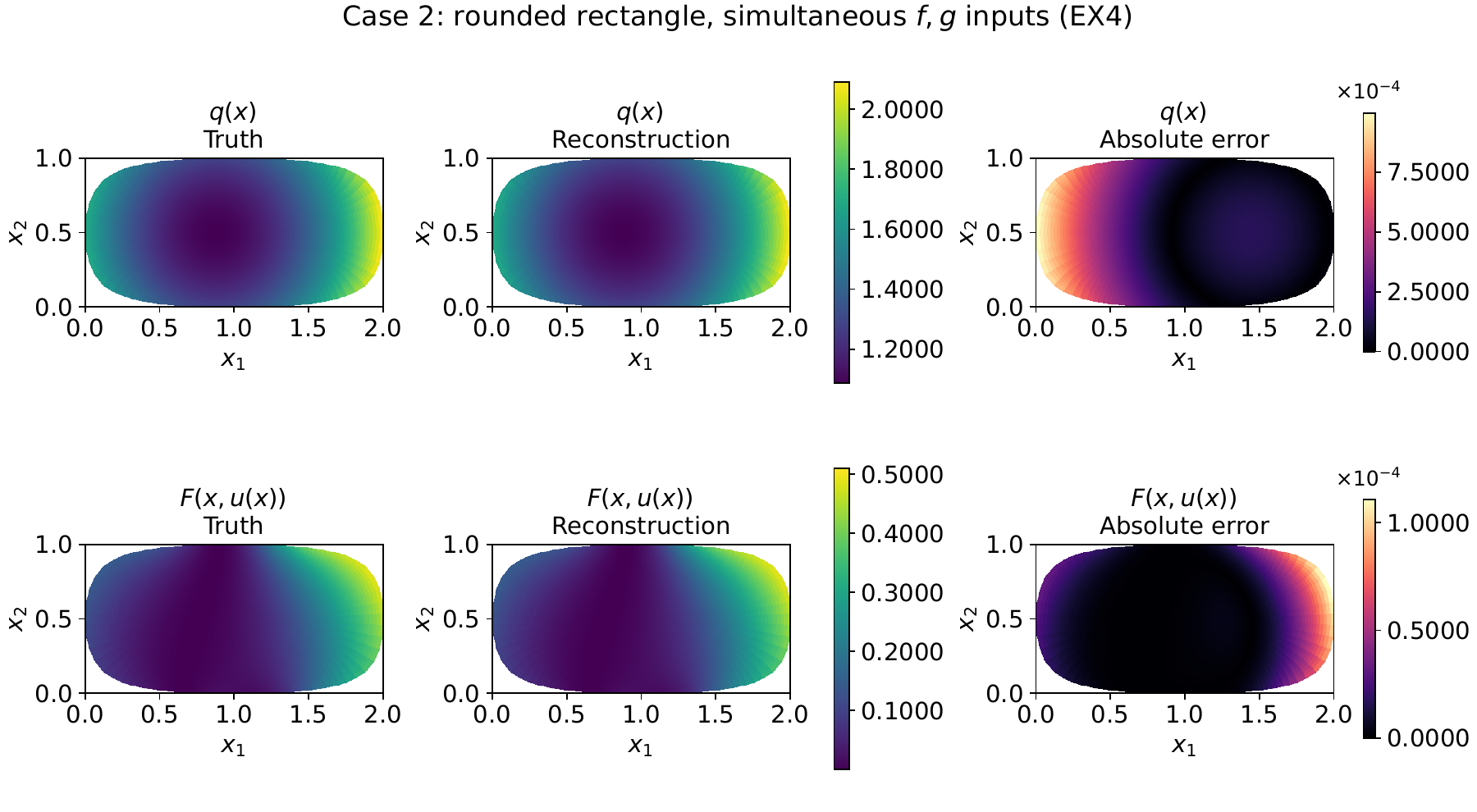}
\caption{\small Example~4 with mixed $(f,g)$ inputs: Case~1 (top) and
Case~2 (bottom). Within each case, the rows show $q$ and the
nonlinear term $F(x,u(x))$; the columns show the reference field,
reconstruction, and absolute difference. The difference in
$F(x,u(x))$ reflects both coefficient and state errors. Each
difference panel uses the indicated color scale.}
\label{fig:num-dual}
\end{figure}

All strategies are evaluated with the same absolute noise standard deviation $\delta$ and the corrected response
\eqref{eq:num-correction}. Table~\ref{tab:num-input-noise} reports the sample mean of $E_{\rm RMS}$ over ten paired noise realizations. The mixed-input strategy has the smallest mean at every tested noise level in both geometries. Nevertheless, at $\delta=10^{-3}$, its mean error reaches $11.37\%$ in Case~1 and $52.77\%$ in Case~2. These values indicate substantial noise sensitivity even for the
mixed-input strategy.
\begin{table}[h]
\centering
\small
\caption{Example~4: mean parameter relative RMS error (\%) over ten noisy trials.
}
\label{tab:num-input-noise}
\begin{tabular}{ccccc}
\toprule
Case & Noise level $\delta$ & $f$ only & $g$ only & Both $f,g$\\
\midrule
Case 1 & $10^{-5}$ & 0.29 & 0.88 & 0.09\\
Case 1 & $10^{-4}$ & 3.40 & 10.78 & 0.83\\
Case 1 & $10^{-3}$ & 38.09 & 64.66 & 11.37\\
Case 2 & $10^{-5}$ & 0.92 & 1.26 & 0.75\\
Case 2 & $10^{-4}$ & 8.16 & 11.70 & 6.11\\
Case 2 & $10^{-3}$ & 95.00 & 110.73 & 52.77\\
\bottomrule
\end{tabular}
\end{table}

Figure~\ref{fig:num-dual} displays the Example~4 mixed-input reconstructions obtained on the refined inversion grids. For a fixed spatially varying boundary input $H=(f,g)$ with both components nonzero, the reference and reconstructed states, $u^\dagger$ and $\widehat u$, are evaluated on a common display grid. The nonlinear terms shown are $F^\dagger(x,u^\dagger(x))$ and $\widehat F(x,\widehat u(x))$. Their difference reflects errors in both the coefficient and the state, rather than coefficient error alone. 

The numerical experiments indicate that, within the tested finite-dimensional classes, simultaneous reconstruction of the potential and the nonlinear coefficients is achievable from finitely many measurements at a single boundary point. In the noiseless tests, refined-grid fits with response correction yield substantially smaller errors than base-grid fits. Mixed Navier boundary inputs produce the largest minimum singular values of the training Jacobian and the lowest mean parameter errors at the tested noise levels, although the reconstruction accuracy deteriorates substantially at the highest noise level.

\section{Conclusions and possible extensions}\label{sec:conclusion}
In this final section, we establish uniqueness results for nonlinearities of the form $F(x,u)=\sum_{k=2}^{M}\frac{a_k(x)}{k!}u^k$ ($M\ge 2$) without assuming spatial analyticity of the coefficients $a_k$, and then discuss possible extensions.

\subsection{Recovery of a nonanalytic coefficient class}

The spatial analyticity of $F$ can be dispensed with when the potential is fixed and the full Navier boundary data $(f,g)$ is available. We will apply the CGO solutions for $\Delta^2$ with a zeroth-order perturbation. The absence of a first-order term gives a stronger remainder estimate than the one needed for general perturbations. 

\begin{lemma}\label{lem:nonanalytic-cgo-product}
Let $3\le n\le6$, let $\Omega\Subset\widetilde\Omega\subset\mathbb R^n$ be smooth bounded domains, and let $q\in L^\infty(\widetilde\Omega)$. Set $p_n=2$ for $3\le n\le5$ and $p_6=5/2$. For every $\xi\in\mathbb R^n$, there are solutions
$v_{1,h},v_{2,h}$ of $(\Delta^2+q)v=0$ in $\widetilde\Omega$ such that
\begin{equation}\label{eq:nonanalytic-cgo-pair}
v_{\ell,h}=e^{ix\cdot\zeta_{\ell,h}/h}(1+s_{\ell,h}),\quad \zeta_{\ell,h}\cdot\zeta_{\ell,h}=0,\quad \zeta_{1,h}+\zeta_{2,h}=h\xi,
\end{equation}
and
\begin{equation}\label{eq:nonanalytic-cgo-convergence}
\|s_{\ell,h}\|_{L^{2p_n}(\Omega)}\rightarrow0, \quad v_{1,h}v_{2,h}\rightarrow e^{ix\cdot\xi}
\quad\text{in }L^{p_n}(\Omega).
\end{equation}
\end{lemma}

\begin{proof}
The phases in \eqref{eq:nonanalytic-cgo-pair} can be chosen as in \cite[Section~3]{KrupchykLassasUhlmann2014}. With
\begin{equation}
P_{\zeta,h}:=-h^2\Delta-2ih\zeta\cdot\nabla,
\end{equation}
the equation for the constant-amplitude remainder is
\begin{equation}\label{eq:nonanalytic-conjugated-remainder}
\bigl(P_{\zeta,h}^2+h^4q\bigr)s_h=-h^4q.
\end{equation}
The solvability estimate of
\cite[Proposition~2.3]{KrupchykLassasUhlmann2014} gives
\begin{equation}\label{eq:nonanalytic-H1-remainder}
\|s_h\|_{H^1_{\mathrm{scl}}(\widetilde\Omega)}\le Ch^2.
\end{equation}
In particular, \eqref{eq:nonanalytic-conjugated-remainder} implies
$\|P_{\zeta,h}^2s_h\|_{L^2(\widetilde\Omega)}=O(h^4)$.
For bounded $\zeta$, the symbol of $P_{\zeta,h}^2$ satisfies
\begin{equation}
\langle\eta\rangle^4
\le C\left(1+\bigl|(|\eta|^2+2\zeta\cdot\eta)^2\bigr|\right).
\end{equation}
The corresponding Fourier estimate, localized by nested cutoffs between $\Omega$ and $\widetilde\Omega$, gives the uniform interior
semiclassical elliptic estimate
\begin{equation}
\|s_h\|_{H^4_{\mathrm{scl}}(\Omega)}
\le C\bigl(\|P_{\zeta,h}^2s_h\|_{L^2(\widetilde\Omega)}
+\|s_h\|_{L^2(\widetilde\Omega)}\bigr)
=O(h^2).
\end{equation}
Semiclassical Sobolev embedding therefore yields
\begin{equation}\label{eq:nonanalytic-L2p-remainder}
\|s_h\|_{L^{2p_n}(\Omega)}\le C h^{-n(1/2-1/(2p_n))}\|s_h\|_{H^4_{\mathrm{scl}}(\Omega)}=O\!\big(h^{2-\frac n2+\frac n{2p_n}}\big).
\end{equation}
The exponent is positive for the stated choices of $p_n$. The exponential factors cancel, and the second convergence in
\eqref{eq:nonanalytic-cgo-convergence} follows from H\"older's inequality.
\end{proof}
The following theorem gives the recovery of a polynomial nonlinearity without spatial analyticity.
\begin{theorem}\label{thm:nonanalytic-polynomial-known-q}
Let $3\le n\le6$, and let $q\in L^\infty(\Omega;\mathbb R)$ satisfy Assumption {\bf (A)}. For $j=1,2$, let
\begin{equation}\label{eq:nonanalytic-polynomial-class}
F_j(x,u)=\sum_{k=2}^{M}\frac{a_{j,k}(x)}{k!}u^k,\quad a_{j,k}\in C^\alpha(\overline\Omega;\mathbb R), \quad 0<\alpha<1,
\end{equation}
where $M\ge2$ is fixed. If $\Lambda_{q,F_1}(H)=\Lambda_{q,F_2}(H)$ for any sufficiently small $H=(f,g)\in\mathcal X_{\mathbb R}$, then
\begin{equation}
a_{1,k}=a_{2,k}\quad {\rm in}\,\, \Omega,\quad 2\le k\le M.
\end{equation}
\end{theorem}

\begin{proof}
For $p=p_n$ as in Lemma~\ref{lem:nonanalytic-cgo-product}, the polynomial Nemytskii map $u\mapsto F_j(\cdot,u)$ is analytic from
$W^{4,p}(\Omega)$ to $L^p(\Omega)$. The analytic implicit-function theorem for the equation $L_qu=F_j(x,u)$ yields real-analytic local measurement maps. We prove $a_{1,k}=a_{2,k}$ by induction. Suppose that $a_{1,k}=a_{2,k}$ in $\Omega$ for  $k=2,\ldots,N-1$. The mixed variations of orders less than $N$ consequently satisfy the same Navier boundary value problems and therefore coincide. Thus, the condition $\Lambda_{q,F_1}(H)=\Lambda_{q,F_2}(H)$ and the transposition identity \eqref{eq:transposition} yield
\begin{equation}\label{eq:nonanalytic-product-identity}
\int_\Omega \Psi_q c_N v_1\cdots v_N\,dx=0, \quad c_N=a_{1,N}-a_{2,N},
\end{equation}
for all solutions $L_qv_\ell=0$ obtained from Navier boundary data.
The identity also holds for complex solutions by complex multilinear
extension. Fix a nontrivial such solution $w$ and set
$v_3=\cdots=v_N=w$ when $N>2$.

Fix $\xi\in\mathbb R^n$. Choose a smooth bounded domain $\Omega'\supset\overline\Omega$, extend $q$ by zero to $\Omega'$, and
apply Lemma~\ref{lem:nonanalytic-cgo-product} to the extended potential. Since $\overline\Omega\subset\Omega'$, interior regularity gives $v\in H^4(\Omega)$. For $n\leq6$, the embedding $H^4(\Omega)\hookrightarrow L^\infty(\Omega)$ and the equation
$\Delta^2v=-qv$ then imply, by interior $L^p$ regularity on $\Omega'$, that $v\in W^{4,p}(\Omega)$ for every finite $p$.
Thus its Navier boundary traces belong to $\mathcal X_{\mathbb C}$. Applying Lemma~\ref{lem:point-adjoint}, with $p=p_n$, we have $\Psi_q\in L^{p_n'}(\Omega)$. 
As $c_Nw\in L^\infty(\Omega)$, H\"older's
inequality and Lemma~\ref{lem:nonanalytic-cgo-product} allow passage
to the limit in \eqref{eq:nonanalytic-product-identity}, yielding
\begin{equation}
\int_\Omega \Psi_q c_N w^{N-2}e^{ix\cdot\xi}\,dx=0 \quad \xi\in\mathbb R^n.
\end{equation}
Since $\Psi_q c_N w^{N-2}\in L^1(\Omega)$, the vanishing of its Fourier transform for every $\xi\in\mathbb R^n$ implies $\Psi_q c_N w^{N-2}=0$ almost everywhere in $\Omega$. Recall that the point adjoint function $\Psi_q$ is nontrivial, and $w$ was chosen nontrivial. Interior elliptic regularity makes both functions continuous, and unique continuation for the system associated with
$L_q$ shows that neither vanishes on an open subset of $\Omega$. Their nonzero sets are therefore dense. Since $c_N$ is continuous, we have $c_N=0$ in $\Omega$, closing the induction.
\end{proof}

\begin{remark}
The class \eqref{eq:nonanalytic-polynomial-class} includes genuinely nonanalytic spatial coefficients. For instance, if
$0\ne b\in C_c^\infty(\Omega)$, then $F(x,u)=(1+b(x))u^2+b(x)u^3$ is covered by the theorem although $b$ is not real analytic in $\Omega$.
\end{remark}

The preceding theorem can also be used when the potential is initially unknown but belongs to a class determined by its normal boundary
derivatives. We state the following corollary.

\begin{corollary}\label{cor:nonanalytic-polynomial-finite-q}
Let $3\le n\le6$. Fix nonempty relatively open sets $\Gamma_a\Subset\Gamma\Subset\partial\Omega\setminus\{x_0\}$. Suppose that $q_{\mathrm{ref}}\in L^\infty(\Omega;\mathbb R)$ and every element of a finite-dimensional space $E\subset L^\infty(\Omega;\mathbb R)$ be smooth up to the boundary near $\overline\Gamma$. Assume that $E$ is determined by its normal derivatives on $\Gamma_a$ in the sense of Definition~\ref{def:jet-visible}. For $j=1,2$, let $q_j\in q_{\mathrm{ref}}+E$ satisfy Assumption {\bf (A)}, and let
\begin{equation}
F_j(x,u)=\sum_{k=m_*}^{M}\frac{a_{j,k}(x)}{k!}u^k,\quad a_{j,k}\in C^\infty(\overline\Omega;\mathbb R), \quad 2\le m_*\le M,
\end{equation}
with $a_{j,m_*}(y)\ne0$ on $y\in\Gamma_a$. If $\Lambda_{q,F_1}(H)=\Lambda_{q,F_2}(H)$ for any sufficiently small $H=(f,g)\in\mathcal X_{\mathbb R}$, then $q_1=q_2$ in $\Omega$ and $F_1=F_2$ on $\Omega\times\mathbb R.$
\end{corollary}

\begin{proof}
Equality of the measurement maps gives
\begin{equation}
D^k\Lambda_{q_1,F_1}(0)=D^k\Lambda_{q_2,F_2}(0), \quad k=1,m_*.
\end{equation}
Since the leading coefficients are nowhere zero on $\Gamma_a$, Theorem~\ref{thm:boundary-jets} yields
\begin{equation}
\partial_\nu^r(q_1-q_2)|_{\Gamma_a}=0, \quad r\geq0.
\end{equation}
As $q_1-q_2\in E$, Definition~\ref{def:jet-visible} implies $q_1=q_2:=q$ in $\Omega$. Theorem~\ref{thm:nonanalytic-polynomial-known-q},
now applied with this common potential, gives $a_{1,k}=a_{2,k}$ for $m_*\leq k\leq M$, and hence $F_1=F_2$.
\end{proof}

\begin{remark}
We take $q_{\mathrm{ref}}=1$, $E=\operatorname{span}\{1\}$, and $q_\theta=1+\theta$ with $|\theta|<1/2$ for example. Then $E$ is determined by its zeroth boundary derivative and every $q_\theta$ is nonresonant. If $0\ne b\in C_c^\infty(\Omega)$ and
$F(x,u)=(1+b(x))u^2+b(x)u^3$, the leading coefficient is nonzero on $\Gamma_a$ although $F$ is not spatially real analytic.


\end{remark}

When $q=0$, the first variations corresponding to the boundary inputs $(f,0)$ are harmonic. This yields a separate partial-input result in all dimensions $n\ge2$.

\begin{proposition}\label{prop:nonanalytic-polynomial-zero-q}
Let $\Omega\subset\mathbb R^n$, $n\ge2$, be a bounded connected domain with smooth boundary, let $x_0\in\partial\Omega$, and let
$\Gamma\subset\partial\Omega$ be nonempty and relatively open. Fix an integer $M\ge2$. For
$j=1,2$, let
\begin{equation}
F_j(x,u)=\sum_{k=2}^{M}\frac{a_{j,k}(x)}{k!}u^k,\quad a_{j,k}\in C^\alpha(\overline\Omega;\mathbb R),
\quad 0<\alpha<1.
\end{equation}
If $\Lambda_{0,F_1}(f,0)=\Lambda_{0,F_2}(f,0)$ for any $f\in C^\infty(\partial\Omega)$ with $\operatorname{supp}f\subset\Gamma$ and sufficiently small $C^{4,\alpha}(\partial\Omega)$ norm, then $a_{1,k}=a_{2,k}$ in $\Omega$ for $2\le k\le M$.
\end{proposition}

\begin{proof}
Since the Navier realization of $\Delta^2$ is invertible, the implicit function theorem yields smooth solution maps
$f\mapsto u_{j,f}\in C^{4,\alpha}(\overline\Omega)$, $j=1,2$, for $f$ in a neighborhood of zero in $C^{4,\alpha}(\partial\Omega)$. We also argue by induction on $N=2,\ldots,M$. At each order, equality of the coefficients below $N$ and of the point measurements gives
\begin{equation}\label{eq:nonanalytic-zero-q-identity}
\int_\Omega \Psi_0(a_{1,N}-a_{2,N})v_1\cdots v_N\,dx=0,
\end{equation}
where $v_j$ is harmonic with smooth boundary trace supported in $\Gamma$. Indeed, the first variations satisfy $\Delta^2v_j=0$ and $\Delta v_j|_{\partial\Omega}=0$, and hence $\Delta v_j=0$ in $\Omega$.

Let $G_D$ and $P_D$ be the Green and Poisson kernels of the Dirichlet Laplacian $-\Delta_D$. The factorization of the Navier inverse as
$(-\Delta_D)^{-2}$ gives
\begin{equation}\label{eq:nonanalytic-zero-q-adjoint}
\Psi_0(x)=-\int_\Omega G_D(x,y)P_D(y,x_0)\,dy<0 \quad x\in\Omega.
\end{equation}
Moreover, $P_D(\cdot,x_0)\in L^1(\Omega)$ by the standard Poisson-kernel bound $P_D(y,x_0)\le C\operatorname{dist}(y,\partial\Omega)
|y-x_0|^{-n}$. If $T=(-\Delta_D)^{-1}1$, then $T\in L^\infty(\Omega)$ and Tonelli's
theorem gives
\begin{equation}
\|\Psi_0\|_{L^1(\Omega)}=\int_\Omega T(y)P_D(y,x_0)\,dy<\infty.
\end{equation}
Choose $0\le h_+\in C_c^\infty(\Gamma)$, $h_+\not\equiv0$, and let $v_+$ be its harmonic extension. Then $v_+>0$ in $\Omega$. Taking
$v_3=\cdots=v_N=v_+$ in \eqref{eq:nonanalytic-zero-q-identity}, with the empty product understood as $1$ when $N=2$, we have
\begin{equation}
\int_\Omega B_N v_1v_2\,dx=0,\quad B_N:=\Psi_0(a_{1,N}-a_{2,N})v_+^{N-2}\in L^1(\Omega).
\end{equation}
The completeness theorem for products of harmonic functions \cite[Theorem~1.3]{SaloTzou2023} implies $B_N=0$ almost everywhere.
By \eqref{eq:nonanalytic-zero-q-adjoint} and $v_+>0$, the difference $a_{1,N}-a_{2,N}$ vanishes almost everywhere, hence in $\Omega$ by continuity.
\end{proof}

\begin{remark}
Unlike Theorem~\ref{thm:nonanalytic-polynomial-known-q}, Proposition~\ref{prop:nonanalytic-polynomial-zero-q} uses only one
Navier input component on $\Gamma$, but requires $q=0$. The restriction $n\leq6$ in Theorem~\ref{thm:nonanalytic-polynomial-known-q} arises from the integrability estimates used to pair the CGO products with the point adjoint function. It does not indicate nonuniqueness of recovering coefficients in higher dimensions.
\end{remark}

The above argument can be directly extended to more general nonlinearities with $q=0$. Suppose that $u\mapsto F_j(\cdot,u)$, $j=1,2$, is real analytic near zero as a $C^\alpha(\overline\Omega)$-valued map and satisfies $F_j(\cdot,0)=\partial_uF_j(\cdot,0)=0$. Assume that the corresponding point-measurement maps agree for all sufficiently small smooth inputs $(f,0)$ with $\operatorname{supp}f\subset\Gamma$. Writing $a_{j,N}=\partial_u^NF_j(\cdot,0)$, suppose inductively that
$a_{1,k}=a_{2,k}$ for $2\le k<N$. By some similar arguments, we can derive
\begin{equation}
\int_\Omega\Psi_0\bigl(a_{1,N}-a_{2,N}\bigr)v_1\cdots v_N\,dx=0.
\end{equation}
The same completeness argument gives $a_{1,N}=a_{2,N}$. Starting at $N=2$, the induction process identifies all Taylor coefficients $\partial_u^kF(x,0)$. The analyticity of $F_j$ in $u$ then gives $F_1(x,u)=F_2(x,u)$ for $x\in\Omega$ and $u$ in a common neighborhood of zero.

However, for the simultaneous recovery of a general smooth potential and nonlinearity, both the first-order linearized solutions and the point-adjoint function depend on the unknown potential. The previous reduction to products of harmonic functions is therefore unavailable.

\subsection{Possible extensins}

At the end of this section, we briefly raise some questions that we leave for future work. In particular, it is more challenge to establish simultaneous global recovery for a general nonlinearity $F$ and an unrestricted smooth potential $q$ without analytic continuation or finite-dimensional assumptions. If both coefficients $q,F$ are unknown, the same question remains open if the Navier boundary inputs are available on an arbitrary relatively open subset of \(\partial\Omega\), and only one component boundary data $(f,0)$ or $(0,g)$ is available. The linear counterexample \eqref{prop:linear-counterexample}, which corresponds to $F\equiv0$, settles neither issue. Moreover, the analytic continuation argument does not yield a quantitative stability estimate for the infinite-dimensional class of potentials. By contrast, under the assumptions of Theorem~\ref{thm:finite-dimensional-stability}, compact finite-dimensional classes are shown to admit a stability modulus and a stability estimate based on finitely many fixed inputs. Therefore, it is also of great interest to study the stable determination of $F$ and $q$ in some infinite-dimensional space from the point measurement map $\Lambda_{q,F}$.


\appendix
\section{Proof of Lemma~\ref{lem:local-boundary-regularity} and Lemma~\ref{lem:boundary-ucp-system}}

\begin{proof}[Proof of Lemma~\ref{lem:local-boundary-regularity}]
We set
\begin{equation}
D_j:=\mathcal O_j\cap\Omega, \quad \Sigma_j:=\mathcal O_j\cap\partial\Omega, \quad j=1,2.
\end{equation}
The boundary operators $\gamma_0u=u|_{\partial\Omega}$ and $\gamma_2u=\Delta u|_{\partial\Omega}$ form the standard Navier elliptic boundary system for $\Delta^2$ and satisfy the Lopatinskii-Shapiro complementing condition.  Since $q$
is a zeroth-order term, it does not affect either the principal symbol or the principal boundary system.  Thus the same complementing condition holds for $L_q=\Delta^2+q.$

Choose an open set $\mathcal V_0$ such that $\mathcal O_1\Subset\mathcal V_0\Subset\mathcal O_2.$
The local $L^s$ boundary estimate for fourth-order elliptic problems \cite[Theorem~15.3, pp.~706]{ADN1959}, gives
\begin{equation}\label{eq:lemma42-basic-ADN}
\|u\|_{W^{4,s}(\mathcal V_0\cap\Omega)}\le C\bigl(\|L_qu\|_{L^s(D_2)}+\|u\|_{L^s(D_2)}+\|\gamma_0u\|_{W^{4-1/s,s}(\Sigma_2)}
+\|\gamma_2u\|_{W^{2-1/s,s}(\Sigma_2)}\bigr).
\end{equation}
The compact set $\overline{\mathcal V}_0\cap\overline\Omega$ can be covered by finitely many coordinate neighborhoods, using boundary charts near $\partial\Omega$ and interior neighborhoods away from the boundary. The Agmon-Douglis-Nirenberg boundary estimate is applied near $\partial\Omega$, while standard interior elliptic estimates are used away from the boundary.
Using the condition $ L_qu=0$ in $D_2$ and $\gamma_0u=\gamma_2u=0$ on $\Sigma_2$, \eqref{eq:lemma42-basic-ADN} reduces to
\begin{equation}
\|u\|_{W^{4,s}(\mathcal V_0\cap\Omega)}\le C_0\|u\|_{L^s(D_2)}.
 \label{eq:lemma42-W4-base}
\end{equation}

It remains to improve the regularity.  Set $w:=\Delta u.$ We still decompose $L_qu=0$ into the coupled second-order
system
\begin{equation}
 \begin{cases}
  \Delta w=-qu,\\
  \Delta u=w,
 \end{cases}
 \quad\text{in }D_2,
 \label{eq:lemma42-second-order-system}
\end{equation}
with homogeneous Dirichlet conditions $u=w=0$ on $\Sigma_2$. We next use the standard local $W^{m,s}$-regularity estimate for the Dirichlet Laplacian. Suppose that $\mathcal V'\Subset\mathcal V\Subset\mathcal O_2$, and that $z\in W^{2,s}_{\mathrm{loc}}$ up to $\mathcal V\cap\partial\Omega$ satisfies
\begin{equation}
\Delta z=f\quad\text{in }\mathcal V\cap\Omega,\quad z=0\quad\text{on }\mathcal V\cap\partial\Omega.
\end{equation}
By the standard local $W^{m+2,s}$-regularity estimate for the Dirichlet-Laplacian problem (see, for instance, \cite[Chapter~2]{Grisvard1985}), if $f\in W^{m,s}(\mathcal V\cap\Omega),$
then
\begin{equation}\label{eq:lemma42-local-Dirichlet}
\|z\|_{W^{m+2,s}(\mathcal V'\cap\Omega)} \le  C\left( \|f\|_{W^{m,s}(\mathcal V\cap\Omega)} +\|z\|_{L^s(\mathcal V\cap\Omega)}\right).
\end{equation}

We next apply the bootstrap arguments.  Suppose that $m\ge4$ and that, on some boundary neighborhood $\mathcal V\Subset\mathcal O_2$, one has
\begin{equation}\label{eq:lemma42-bootstrap-hypothesis}
\|u\|_{W^{m,s}(\mathcal V\cap\Omega)}\le C_m\|u\|_{L^s(D_2)}.
\end{equation}
Choose two smaller sets $\mathcal V''\Subset\mathcal V'\Subset\mathcal V.$
Since $q$ is smooth up to the boundary in $\mathcal O_2$,  we have
\begin{equation}\label{eq:lemma42-q-multiplier}
\|qu\|_{W^{m,s}(\mathcal V\cap\Omega)}\le C_m \|q\|_{W^{m,\infty}(\mathcal V\cap\Omega)} \|u\|_{W^{m,s}(\mathcal V\cap\Omega)}.
\end{equation}
Applying \eqref{eq:lemma42-local-Dirichlet} to the first equation in
\eqref{eq:lemma42-second-order-system}, using $ w=\Delta u\in W^{m-2,s}(\mathcal V\cap\Omega)$ and $m\ge4$, we have
$w\in W^{2,s}(\mathcal V\cap\Omega)$.  Moreover, $w=0$ on the boundary portion $\mathcal V\cap\partial\Omega$.  Hence
\begin{equation}\label{eq:lemma42-w-estimate}
\|w\|_{W^{m+2,s}(\mathcal V'\cap\Omega)} \le C\bigl(\|qu\|_{W^{m,s}(\mathcal V\cap\Omega)} +\|w\|_{L^s(\mathcal V\cap\Omega)}
\bigr).
\end{equation}
Note that
\begin{equation}
\|w\|_{L^s(\mathcal V\cap\Omega)}\le C\|u\|_{W^{m,s}(\mathcal V\cap\Omega)}.
\end{equation}
Combining this estimate with \eqref{eq:lemma42-bootstrap-hypothesis} and \eqref{eq:lemma42-q-multiplier}, we can obtain
\begin{equation}
\|w\|_{W^{m+2,s}(\mathcal V'\cap\Omega)}\le C_m'\|u\|_{L^s(D_2)}.
 \label{eq:lemma42-w-higher}
\end{equation}
Since $u=0$ on $\mathcal V'\cap\partial\Omega$, using $\mathcal V''\Subset\mathcal V'$ and \eqref{eq:lemma42-w-higher}, applying \eqref{eq:lemma42-local-Dirichlet} to the second equation in \eqref{eq:lemma42-second-order-system}, we can get
\begin{equation} \label{eq:lemma42-bootstrap-step}
\|u\|_{W^{m+4,s}(\mathcal V''\cap\Omega)}\le C\left(\|w\|_{W^{m+2,s}(\mathcal V'\cap\Omega)}+\|u\|_{L^s(\mathcal V'\cap\Omega)}\right)\le C_{m+4}\|u\|_{L^s(D_2)}.
\end{equation}
Thus the local regularity improves from $W^{m,s}$ to $W^{m+4,s}$ on a smaller boundary neighborhood.

Fix $N\ge1$ and choose sets
\begin{equation}
\mathcal O_1\Subset\mathcal V_{2N-2}\Subset\mathcal V_{2N-3}\Subset\cdots\Subset\mathcal V_1\Subset\mathcal V_0
\Subset\mathcal O_2.
\end{equation}
When $N=1$, this is understood simply as $\mathcal O_1\Subset\mathcal V_0\Subset\mathcal O_2.$ The estimate \eqref{eq:lemma42-W4-base} gives the initial $W^{4,s}$-bound on $\mathcal V_0$.  For $r=0,\ldots,N-2$, applying the bootstrap step with
\begin{equation}
m=4(r+1),\quad \mathcal V=\mathcal V_{2r},\quad \mathcal V'=\mathcal V_{2r+1},\quad \mathcal V''=\mathcal V_{2r+2},
\end{equation}
the induction argument gives
\begin{equation}
\|u\|_{W^{4N,s}(\mathcal V_{2N-2}\cap\Omega)}\le C_{N,s}\|u\|_{L^s(D_2)}.
\end{equation}
Since $D_1\subset\mathcal V_{2N-2}\cap\Omega,$ we obtain
\begin{equation}\label{eq:lemma42-W4N}
\|u\|_{W^{4N,s}(D_1)}\le C_{N,s}\|u\|_{L^s(D_2)}.
\end{equation}
Given an arbitrary integer $k\ge0$, choose $N\ge1$ such that $4N\ge k.$ By the embedding $W^{4N,s}(D_1)\hookrightarrow W^{k,s}(D_1),$ it follows from \eqref{eq:lemma42-W4N} that
\begin{equation}
\|u\|_{W^{k,s}(D_1)}\le C_{k,s}\|u\|_{L^s(D_2)}.
\end{equation}
This proves the desired estimate.

For each fixed $k$, only finitely many derivatives of $q$ are involved in the previous argument. Consequently, the constant $C_{k,s}$ depends only on $k$, $s$, the neighborhoods $\mathcal O_1$ and $\mathcal O_2$, the local geometry of $\partial\Omega$, and finitely many $C^j$-norms of $q$ on $\mathcal O_2\cap\overline\Omega$.
\end{proof}

\begin{proof}[Proof of Lemma~\ref{lem:boundary-ucp-system}]
Fix $x_*\in U_*$ and choose a ball $B$ centered at $x_*$ such that $B\cap\partial\Omega\Subset U_*.$ After decreasing the radius of $B$ if necessary, the functions $u$ and $w$ belong to $H^2(B\cap\Omega)$. For $v=u$ or $v=w$, the assumptions give $v=\partial_\nu v=0$ on $B\cap\partial\Omega.$ The standard local zero-extension criterion for $H^2$ functions on a
smooth domain therefore implies that the extension
\begin{equation}
\widetilde v(x):=
\begin{cases}
v(x),&x\in B\cap\Omega,\\
0,&x\in B\setminus\overline\Omega,
\end{cases}
\end{equation}
belongs to $H^2_{\mathrm{loc}}(B)$.  Indeed, the vanishing of both the Dirichlet and normal traces prevents the first and second distributional derivatives of $\widetilde v$ from acquiring terms supported on $B\cap\partial\Omega$.

Extend $c$ to a function $\widetilde c\in L^\infty(B)$, for instance by setting it equal to zero in $B\setminus\overline\Omega$.  We claim that
\begin{equation}\label{eq:extended-ucp-system}
\Delta\widetilde u=\widetilde w,\quad\Delta\widetilde w=-\widetilde c\,\widetilde u \quad\text{in }\mathcal D'(B).
\end{equation}
To verify the first identity, let $\phi\in C_c^\infty(B)$.  Green's
formula and the vanishing Cauchy data of $u$ give the following identity
\begin{align}
\langle\Delta\widetilde u,\phi\rangle & =\int_{B\cap\Omega}u\,\Delta\phi\,dx=\int_{B\cap\Omega}(\Delta u)\phi\,dx+\int_{B\cap\partial\Omega}\bigl(u\,\partial_\nu\phi-\phi\,\partial_\nu u\bigr)\,dS \notag\\
&=\int_{B\cap\Omega}w\phi\,dx=\int_B\widetilde w\phi\,dx.
\end{align}
The second identity follows in the same way from $w=\partial_\nu w=0$ and $\Delta w=-cu$.  Thus, with
\begin{equation}
\widetilde{\mathbf U}:=
\begin{pmatrix}
\widetilde u\\
\widetilde w
\end{pmatrix},
\quad
\widetilde C(x):=
\begin{pmatrix}
0&-1\\
\widetilde c(x)&0
\end{pmatrix},
\end{equation}
equation \eqref{eq:extended-ucp-system} becomes
\begin{equation}\label{eq:extended-vector-equation}
\Delta\widetilde{\mathbf U}+\widetilde C(x)\widetilde{\mathbf U}=0 \quad\text{in }\mathcal D'(B),\quad \widetilde C\in L^\infty(B;\mathbb R^{2\times2}).
\end{equation}
The principal part of this system is the diagonal strongly elliptic operator.  Hence the local unique-continuation theorem for
second-order elliptic systems \cite[Theorem~3.1]{HondaLinNakamuraSasayama2022} applies.  Although that theorem is stated for $C^2$ solutions, the Carleman estimate used in its proof is an $L^2$ estimate and extends, by approximation, to $H^2_{\mathrm{loc}}$ solutions. Since $ \widetilde{\mathbf U}=0$ in the nonempty open set $B\setminus\overline\Omega,$
unique continuation yields $\widetilde{\mathbf U}=0$ in $B$. This implies $u=w=0$ in $B\cap\Omega,$ which proves the local assertion \eqref{local-vanishing-u}. Finally, suppose that $c\in L^\infty(\Omega)$ and that
\eqref{eq:ucp-system} holds in the connected domain $\Omega$. The local conclusion provides a nonempty open set $\mathcal O_0\subset\Omega$ on which $\mathbf U=(u,w)^T$ vanishes. Considering the equation
\begin{equation}
\Delta\mathbf U+
\begin{pmatrix}
0&-1\\
c&0
\end{pmatrix}
\mathbf U=0
\quad\text{in }\Omega,
\end{equation}
successive application of the same local result along a finite chain of overlapping balls yields $\mathbf U\equiv0$ in $\Omega$. Therefore the proof is finished.
\end{proof}

\subsection*{ Acknowledgments}  

S. Fu is supported by the Fundamental Research Funds for the Central Universities, NPU, under grant number D5000250416, the Natural Science Basic Research Program of Shaanxi under grant number 2026JC-YBQN-0003, and the key program of the National Natural Science Foundation of China under grant number 62433020. H. Liu is supported by the Hong Kong RGC General Research Funds (projects 11311122,  11304224, and 11303125)

\subsection*{Data availability statement}
Data available on request from the authors.
	
\subsection*{Conflict of Interests} Hereby we declare there are no conflicts of interest.

\bibliography{references} 
	
\bibliographystyle{alpha}

\end{document}